%% file: main.tex
\documentclass{amsart}

\input{preamble}
\input{macros}

\providecommand{\figalttext}[2][]{}

\title[Equivariant Quantisation via Symplectic Cuts]{Equivariant Quantisation of Hypertoric Varieties via Symplectic Cuts}

\author{Benjamin C. W. Brown}
\address{School of Mathematics, The University of Edinburgh, James Clerk Maxwell Building, Peter Guthrie Tait Road, Edinburgh EH9 3FD, United Kingdom}

\urladdr{\url{https://orcid.org/0000-0001-8699-2677}}


\hypersetup{pdftitle={Equivariant Quantisation of Hypertoric Varieties via Symplectic Cuts},pdfauthor={Benjamin C. W. Brown}}

\begin{document}

\begin{abstract}
	We study the equivariant geometric quantisation of a non-compact hypertoric orbifold of complex dimension $2n$ by compact symplectic cuts.
	When the defining arrangement has a non-empty bounded base region and the fixed lift is integral, we prove that, at every positive regular integral level $d$, the space of sections on the cut is the sum of the fibre-rotation weight spaces of weight at most $d$ in the original quantisation.
	The higher cohomology of the cut pre-quantum line orbibundle vanishes, so its equivariant index is the character of its actual section space.
	Using the moment polyptychs introduced in a companion paper, we give an orbifold Kawasaki--Riemann--Roch localisation formula for this character.
	We also compute it as a weighted lattice-point sum: a residual weight $x$ occurs at level $d$ if and only if the arrangement slack function satisfies $\mcS(x)\leq d$, with multiplicity $\binom{d - \mcS(x) + n}{n}$.
	The examples compare these calculations in rank one, including a cut with an orbifold point, and in rank two, for smooth cuts that are not toric.
\end{abstract}

\maketitle

\input{sections/1-introduction}

\input{sections/2-geometric-input}

\input{sections/3-localisation-formulae}

\input{sections/4-hypertoric-equivariant-quantisation}

\input{sections/5-algebraic-cutting-vanishing}

\input{sections/6-combinatorial-quantisation}

\input{sections/7-conclusion}

\printbibliography

\end{document}

%% file: preamble.tex
\usepackage[T1]{fontenc}
\usepackage[utf8]{inputenc}
\usepackage[lining]{ebgaramond}
\usepackage{microtype}
\usepackage{dsfont}

\usepackage{amsmath,amssymb,amsthm}
\usepackage{stmaryrd} % provides \sslash
\usepackage{braket}   % provides \Set

\usepackage{geometry}
\usepackage[skip=0pt,indent=20pt]{parskip}

\usepackage[shortlabels]{enumitem}

\usepackage{graphicx}
\usepackage{float}
\usepackage{tikz}
\usepackage{tikz-cd}
\usetikzlibrary{calc,arrows.meta}

\usepackage{subcaption}

\usepackage[
	backend=biber,
	style=alphabetic,
	sorting=nyt,
	maxalphanames=4,
	minalphanames=4,
	maxbibnames=99
]{biblatex}

\usepackage{hyperref}

\hypersetup{
	hidelinks
}

\tikzset{
	polyboundary/.style={line width=0.9pt},
	polyedge/.style={line width=0.45pt},
	vertex/.style={circle,fill,inner sep=1.35pt},
	region/.style={font=\scriptsize,inner sep=1pt},
}

\theoremstyle{plain}
\newtheorem{theorem}{Theorem}[section]
\newtheorem{proposition}[theorem]{Proposition}
\newtheorem{corollary}[theorem]{Corollary}
\newtheorem{lemma}[theorem]{Lemma}

\theoremstyle{definition}
\newtheorem{definition}[theorem]{Definition}
\newtheorem{example}[theorem]{Example}

\theoremstyle{remark}
\newtheorem{remark}[theorem]{Remark}
\newtheorem{notation}[theorem]{Notation}

%% file: macros.tex
\newcommand{\ra}{\rightarrow}
\newcommand{\lra}{\longrightarrow}

\newcommand{\w}{\omega}

\newcommand{\PP}{\mathbb{P}}
\newcommand{\RR}{\mathbb{R}}
\newcommand{\CC}{\mathbb{C}}
\newcommand{\QQ}{\mathbb{Q}}

\newcommand{\ZZ}{\mathbb{Z}}

\newcommand{\one}{\mathds{1}}

\newcommand{\mcA}{\mathcal{A}}

\newcommand{\mcC}{\mathcal{C}}
\newcommand{\mcE}{\mathcal{E}}

\newcommand{\mcG}{\mathcal{G}}
\newcommand{\mcI}{\mathcal{I}}
\newcommand{\mcIv}{\mathcal{I}_{v}}

\newcommand{\mcL}{\mathcal{L}}

\newcommand{\mcN}{\mathcal{N}}
\newcommand{\mcO}{\mathcal{O}}
\newcommand{\mcP}{\mathcal{P}}

\newcommand{\mcS}{\mathcal{S}}

\newcommand{\mcZ}{\mathcal{Z}}

\newcommand{\mfk}{\mathfrak{k}}
\newcommand{\mft}{\mathfrak{t}}

\makeatletter
\newcount\wb@n
\newdimen\wb@delta
\newcommand{\wb@xht}[1]{%
	\ifx#1\displaystyle\fontdimen5\textfont2 \else
	\ifx#1\textstyle\fontdimen5\textfont2 \else
	\ifx#1\scriptstyle\fontdimen5\scriptfont2 \else
	\fontdimen5\scriptscriptfont2 \fi\fi\fi}
\newsavebox{\wb@nuc}
\newsavebox{\wb@acc}
\def\wb@rep{%
	\ifnum\wb@n>\z@
	\copy\wb@acc \kern-\wd\wb@acc \kern\wb@delta
	\advance\wb@n\m@ne
	\expandafter\wb@rep
	\fi}
\newcommand{\wb@draw}[2]{%
	\sbox\wb@nuc{$\m@th#1#2$}%
	\sbox\wb@acc{$\m@th#1\mathchar"0016$}%
	\@tempdimb=\ht\wb@acc
	\dimen@=\ht\wb@nuc \advance\dimen@-\wb@xht{#1}\relax
	\ifdim\dimen@>\z@ \advance\@tempdimb\dimen@ \fi
	\advance\@tempdimb-\ht\wb@acc
	\ifdim\wd\wb@nuc>\wd\wb@acc
	\wb@delta=.25\wd\wb@acc                % overlap generously
	\wb@n=\@ne \@tempdimc=\wd\wb@acc
	\@whiledim\@tempdimc<\wd\wb@nuc\do{%
		\advance\@tempdimc\wb@delta \advance\wb@n\@ne}%
	\wb@delta=\dimexpr(\wd\wb@nuc-\wd\wb@acc)/\numexpr\wb@n-1\relax\relax
	\advance\wb@n\m@ne                     % that many gaps
	\else
	\wb@n=\z@ \wb@delta=\z@
	\fi
	\mathord{\rlap{\copy\wb@nuc}%
		\raise\@tempdimb\hbox to\wd\wb@nuc{\hfil\wb@rep\copy\wb@acc\hfil}}}
\newcommand{\widebar}[1]{\mathpalette\wb@draw{#1}}
\makeatother

\newcommand{\Cout}{\mathfrak{C}_{\mathrm{out}}}

\newcommand{\diff}{\mathrm{diff}}
\newcommand{\aux}{\mathrm{aux}}

\newcommand{\alg}{\mathrm{alg}}

\renewcommand{\aux}{\mathrm{aux}}

\newcommand{\aff}{\mathrm{aff}}
\newcommand{\tang}{\mathrm{tan}}
\newcommand{\fib}{\mathrm{fib}}

\DeclareMathOperator{\CritVal}{CritVal}

\DeclareMathOperator{\supp}{supp}
\DeclareMathOperator{\mult}{mult}
\DeclareMathOperator{\Conv}{Conv}

\DeclareMathOperator{\Verts}{Vert}

\DeclareMathOperator{\Identity}{Id}

\DeclareMathOperator{\Image}{im}
\DeclareMathOperator{\Span}{Span}

\DeclareMathOperator{\Td}{Td}
\DeclareMathOperator{\Ch}{Ch}

\DeclareMathOperator{\rk}{rk}

\DeclareMathOperator{\codim}{codim}

\DeclareMathOperator{\HK}{HK}
\DeclareMathOperator{\cut}{cut}
\DeclareMathOperator{\Proj}{Proj}
\DeclareMathOperator{\Spec}{Spec}
\DeclareMathOperator{\sstable}{ss}

\DeclareMathOperator{\diag}{diag}
\DeclareMathOperator{\wt}{wt}

\DeclareMathOperator{\SU}{SU}

%% file: sections/1-introduction.tex
% !TeX root = ../main.tex
\section{Introduction}
\label{section:section-1-introduction}

For a compact K\"ahler manifold equipped with a pre-quantum line bundle, Riemann--Roch computes the holomorphic Euler characteristic of the line bundle, which is the quantisation.
On a non-compact space, the section space may be infinite-dimensional, but a compatible circle action can provide a grading with finite-dimensional weight spaces.
We study this situation for hypertoric varieties: compact symplectic cuts realise finite pieces of their fibre-rotation grading, and orbifold localisation and weighted lattice-point counts compute the same characters.

One motivation comes from the equivariant Verlinde formula for Higgs-bundle moduli, proposed by Gukov and Pei \cite{GukovPei2017}.
It was proved by Andersen, Gukov, and Pei \cite{AndersenGukovPei2016} and, independently and algebro-geometrically, by Halpern-Leistner \cite{HalpernLeistner2016}.
A related geometric model is Weitsman's real polarisation \cite{Weitsman1992} of the moduli space of flat $\SU(2)$-connections: Jeffrey and Weitsman identify the Bohr--Sommerfeld fibres \cite{JeffreyWeitsman1992}, whose number is given by the Verlinde formula \cite{Verlinde1988}, and describe a dense open toric subset and a toric compactification \cite{JeffreyWeitsman1994}.
The latter has a different complex structure from the moduli space, so its lattice-point count cannot be compared directly with the algebraic section spaces \cite[p.~153]{JeffreyWeitsman1994}.
For hypertoric cuts, on the other hand, a weighted lattice-point count and the algebraic section spaces can be compared directly.

Two direct antecedents are Martens's use of compact cuts to compute equivariant integrals on non-compact symplectic and hyperk\"ahler quotients \cite{Martens2008}, and Hausel's compactification of Higgs-bundle moduli \cite{Hausel1998} by the symplectic cut of Lerman \cite{Lerman1995}.
Hausel suggests applying this construction to toric hyperk\"ahler manifolds \cite[\S1]{Hausel1998}, introduced by Bielawski and Dancer \cite{BielawskiDancer2000}.
The companion paper \cite{BrownI} develops these hypertoric cuts and their moment polyptychs.
Here, we treat the cuts as individual finite stages of quantisation and compare their sections with those of the original non-compact variety.

Let $M_{\nu}$ be a non-empty hypertoric orbifold of real dimension $4n$ at a generic parameter, with residual torus $T^{n}$.
The primitive normal vectors $u_{i} \in \mft_{\ZZ}^{n}$ of $M_{\nu}$ and an integral lift $\lambda$ of $\nu$ give the affine-linear functionals $\ell_{i}(x) = \langle x ,\, u_{i} \rangle + \lambda_{i}$ and the arrangement of hyperplanes $H_{i} = \{ \ell_{i} = 0 \}$.
We coorient the arrangement so that its distinguished base region is non-empty and bounded.
Fibre rotation then has a proper moment map $\Phi$ \cite[Proposition~1.3 and the following remark]{HaradaProudfoot2004}.
We write $M_{\nu}^{\leq \delta}$ for the compact K\"ahler orbifold cut at a positive regular level $\delta$; when the level is an integer $d$, the cut $M_{\nu}^{\leq d}$ carries a pre-quantum line orbibundle $\mcL_{\nu}^{\leq d}$.

Our main comparison holds at every positive regular integral level: in the outer cut chamber $d > c_{\max}$, above the largest critical value of $\Phi$, and in the inner chambers below it.
With algebraic sections in complex structure $J_{1}$ and the pullback action of Notation~\ref{notation:representation-convention}, fibre rotation defines the exhaustive filtration
\[
	F_{d}H^{0}(M_{\nu} ; \mcL_{\nu}) :=
	\bigoplus\nolimits_{i \leq d}H^{0}(M_{\nu} ; \mcL_{\nu})_{i}, \qquad d \in \ZZ_{\geq 0}.
\]
Algebraic cutting, following Edidin and Graham \cite{EdidinGraham1998}, identifies these pieces with invariant sections on $M_{\nu} \times \CC$.
Restriction and descent give
\[
	F_{d}H^{0}(M_{\nu} ; \mcL_{\nu})
	\overset{\sim}{\longrightarrow}
	H^{0}(M_{\nu}^{\leq d} ; \mcL_{\nu}^{\leq d}),
	\qquad
	H^{q}(M_{\nu}^{\leq d} ; \mcL_{\nu}^{\leq d}) = 0 \quad ( q > 0 )
\]
for every positive regular integral level $d$; this is Theorem~\ref{theorem:cut-cohomology-every-level}, which more generally compares the cohomology in each fibre-rotation weight.
Its proof applies Teleman's quantisation theorem \cite{Teleman2000} in the form of Halpern-Leistner \cite[Theorem~3.29]{HalpernLeistner2015}, whose hypotheses we verify by showing that the unstable strata of the cut are closed blades of fibre-rotation cells (Proposition~\ref{prop:cells-are-blades}) and that the pre-quantum weights along the fixed components have the required signs (Lemma~\ref{lemma:stratum-weights}).
The cut section spaces are therefore finite stages of the non-compact quantisation; successive filtration quotients recover the individual weight spaces.

The arrangement determines the weights and multiplicities of these finite stages.
Let $\mcS(x) = \sum_{i = 1}^{N} \max\{ 0 , -\ell_{i}(x) \}$ be the slack function.
The moment polyptych $\Delta_{\nu}^{\leq d}$ is assembled from truncated arrangement regions, and its support is $\{ \mcS \leq d \}$ \cite[Proposition~4.6]{BrownI}.
Every $F_{d}$ is finite-dimensional, and a residual weight $x$ occurs with multiplicity $\binom{d - \mcS(x) + n}{n}$ when $\mcS(x) \leq d$ and with multiplicity zero otherwise; this is Theorem~\ref{theorem:combinatorial-multiplicity}.
Therefore, at a positive regular integral level,
\[
	\operatorname{char}_{T^{n}} H^{0}(M_{\nu}^{\leq d} ; \mcL_{\nu}^{\leq d})(e^{\xi}) = \sum_{x \in |\Delta_{\nu}^{\leq d}| \cap ( \mft_{\ZZ}^{n} )^{\ast}} \binom{d - \mcS(x) + n}{n} e^{\langle x ,\, \xi \rangle}.
\]
The character of $H^{0}( M_{\nu} ; \mcL_{\nu} )$ for $T^{n}$ and the fibre-rotation action agrees, after a change of grading, with the abelian monopole formula with background flux \cite[\S5(ii)]{Nakajima2016}; this relation, and the role of the semi-invariant section comparison in it, are explained in Remark~\ref{remark:monopole-formula}.
The monopole comparison uses the identification of the sections of $\mcL_{\nu}^{\otimes m}$ with the semi-invariants of weight $m\nu$ on $\phi_{\CC}^{-1}(0)$, which is proved here in full as Proposition~\ref{prop:positive-ray-semi-invariants}.

Since the cut realises each finite stage as the quantisation of a compact orbifold, the character of $H^{0}( M_{\nu}^{\leq d} ; \mcL_{\nu}^{\leq d} )$ also has a local geometric expression.
The fixed-point and tangent data of~\cite{BrownI}, together with the pre-quantum fibre characters, give a Kawasaki localisation formula for it, which is Theorem~\ref{theorem:character-hypertoric-orbifold}.
Even for a smooth hypertoric manifold, the orbifold terms can be needed, since a reducible core forces singular outer-chamber cuts \cite[Theorem~4.18]{BrownI}.
The cohomology comparison identifies the localised index with the actual section-space character, so localisation and the weighted lattice sum compute the same representation.

Four examples from~\cite{BrownI} illustrate this comparison.
In rank one, the cuts of $T^{\ast}\CC\PP^{1}$ are smooth toric four-manifolds, so their characters also follow from the lattice-point description of toric quantisation \cite{Hamilton2008}.
The outer-chamber cuts of the $A_{2}$-resolution have a fixed point with local group $\ZZ / 2\ZZ$ and, for one parity of the level, a non-trivial fibre character.
In rank two, the cuts of $T^{\ast}\CC\PP^{2}$ and $T^{\ast}( \CC\PP^{1} \times \CC\PP^{1} )$ are smooth but not toric, so they test Theorem~\ref{theorem:character-hypertoric-orbifold} beyond the toric case.

\subsection*{Organisation}
\label{subsection:subsection-1.1-organisation}

Section~\ref{section:section-2-conventions-geometric-input} fixes the conventions and records the geometric inputs from the companion paper.
Section~\ref{section:section-3-localisation-formulae} recalls orbifold localisation and constructs the cut pre-quantum line orbibundle, and Section~\ref{section:section-4-equivariant-quantisation-cut-spaces} gives the cut localisation formula.
Section~\ref{section:section-5-algebraic-cutting-vanishing} proves the cohomology comparison at every positive regular integral level, and Section~\ref{section:section-6-combinatorial-quantisation} derives the multiplicities and compares the two character calculations in the examples.
Section~\ref{section:section-7-conclusion} discusses two questions arising from these results.

%% file: sections/2-geometric-input.tex
% !TeX root = ../main.tex
\section{Conventions and Geometric Input}
\label{section:section-2-conventions-geometric-input}

We begin by fixing the quotient and line-bundle conventions used throughout the paper, then recalling the cut geometry and local data from the companion paper~\cite{BrownI}.

\subsection{Hypertoric Quotients}
\label{subsection:subsection-2.1-hypertoric-quotients}

We use the hyperk\"ahler quotient construction of \cite{HitchinKarlhedeLindstromRocek1987}, with the conventions of \cite[Section~2]{BrownI}.
For a torus $T = \mft / \mft_{\ZZ}$, we use the period-one convention: a weight $\alpha \in \mft_{\ZZ}^{\ast}$ has the character $\chi_{\alpha}( [ \xi ] ) = \exp( 2\pi\sqrt{-1} \langle \alpha ,\, \xi \rangle )$.
Moment maps satisfy $d\langle \mu ,\, \xi \rangle = \w( - ,\, X_{\xi} )$.
Primitive vectors $u_{1} , \ldots , u_{N}$ generating $\mft_{\ZZ}^{n}$ give an exact cocharacter sequence
\begin{equation}
	\label{eqn:lie-algebra-seq}
	\begin{tikzcd}[ampersand replacement=\&]
		0 \ar[r] \& \mfk_{\ZZ} \ar[r, hook, "\imath_{\ZZ}"] \& \mft_{\ZZ}^{N} \ar[r, "\pi_{\ZZ}"] \& \mft_{\ZZ}^{n} \ar[r] \& 0,
	\end{tikzcd}
\end{equation}
with $\pi_{\ZZ}(e_{i}) = u_{i}$, and hence a subtorus $K \subset T^{N}$ with residual torus $T^{n} = T^{N}/K$.
Let $\imath : \mfk \hookrightarrow \mft^{N}$ be the inclusion, $\pi : \mft^{N} \ra \mft^{n}$ the projection, $\epsilon_{i}$ the basis dual to $e_{i}$, and $\alpha_{i} = \imath_{\ZZ}^{\ast}(\epsilon_{i})$.
The dual sequence of \eqref{eqn:lie-algebra-seq} identifies $( \mft_{\ZZ}^{n} )^{\ast}$ with $\ker \imath_{\ZZ}^{\ast}$, and shows that every integral $\nu \in \mfk_{\ZZ}^{\ast}$ has an integral lift $\lambda = \sum_{i} \lambda_{i}\epsilon_{i}$.

On $T^{\ast}\CC^{N}$, the action is $t \cdot ( z , w ) = ( t_{i} z_{i} , t_{i}^{-1}w_{i} )_{i = 1}^{N}$, and our forms are
\[
	\w_{\RR} = \frac{\sqrt{-1}}{4\pi} \sum\nolimits_{i} ( dz_{i} \wedge d\bar{z}_{i} + dw_{i} \wedge d\bar{w}_{i} ),
	\qquad
	\w_{\CC} = \frac{\sqrt{-1}}{2\pi} \sum\nolimits_{i} dz_{i} \wedge dw_{i}.
\]
The $K$-moment maps are $\phi_{\RR} = \tfrac{1}{2} \sum_{i} ( |z_{i}|^{2} - |w_{i}|^{2} ) \alpha_{i}$ and $\phi_{\CC} = \sum_{i} (z_{i}w_{i}) \alpha_{i}$.
Fix a parameter $\nu \in \mfk^{\ast}$ and a lift $\lambda \in (\mft^{N})^{\ast}$, i.e., $\nu = \imath^{\ast}(\lambda) \in \mfk^{\ast}$, and assume $\nu$ is generic, which means that $( \nu , 0 )$ is a regular value of $\phi_{\HK} = ( \phi_{\RR} , \phi_{\CC} )$.
Then $M_{\nu} := \phi_{\HK}^{-1}(\nu , 0) / K$ is a non-empty hyperk\"ahler orbifold of real dimension $4n$; finite stabilisers are allowed throughout.
In its complex structure $J_{1}$, and writing $[z, w]$ for the class of $(z, w) \in \phi_{\HK}^{-1}(\nu, 0)$, the residual $T^{n}$-moment maps are determined by
\[
	\pi^{\ast} \circ \mu_{\RR}[z , w] = \tfrac{1}{2} \sum\nolimits_{i} ( |z_{i}|^{2} - |w_{i}|^{2} ) \epsilon_{i} - \lambda,
	\qquad
	\pi_{\CC}^{\ast} \circ \mu_{\CC}[z , w] = \sum\nolimits_{i} (z_{i}w_{i}) \epsilon_{i}.
\]

\subsection{Algebraic Quotients and Pre-Quantum Conventions}
\label{subsection:subsection-2.2-algebraic-quotients-prequantum}

For the algebraic and pre-quantum constructions, take $\nu$ and the fixed lift $\lambda$ to be integral.
Let $\mcN := \phi_{\CC}^{-1}(0)$, and let $K_{\CC}$ and $T_{\CC}^{N}$ act by complexifying the actions above, with exact sequence $1 \ra K_{\CC} \ra T_{\CC}^{N} \ra T_{\CC}^{n} \ra 1$.
The trivial line bundle $L_{\nu} := \mcN \times \CC$ has the geometric linearisation
\begin{equation}
	\label{eqn:positive-character-linearisation}
	k \cdot ( x , v ) = ( k \cdot x , \chi_{\nu}(k) v ),
	\qquad k \in K_{\CC},
\end{equation}
extended to $T_{\CC}^{N}$ by the character $\chi_{\lambda}$.

\begin{notation}[Representation Convention]
	\label{notation:representation-convention}
	For a geometric action on a line bundle $L$, with fibre map $t_{L}$, sections carry the pullback action $( t \cdot s )(x) := t_{L}^{-1} s( t \cdot x )$.
	Thus, the $K_{\CC}$-invariant sections of $L_{\nu}^{\otimes m}$ are the semi-invariants $\CC[\mcN]_{m\nu} := \{ f \mid f( k \cdot x ) = \chi_{\nu}(k)^{m}f(x) \}$.
	If a monomial $z^{a}w^{b}$ satisfies $\imath_{\ZZ}^{\ast}( a - b ) = \nu$, write $a - b = \lambda + \pi_{\ZZ}^{\ast}(x)$; its descended section has residual weight $+x$.
\end{notation}

\begin{definition}
	\label{def:GIT-hypertoric-variety}
	The hypertoric variety at $\nu$ as an algebraic quotient and its semi-invariant ring are
	\begin{equation}
		\label{eqn:GIT-hypertoric-definition}
		M_{\nu}^{\alg} := \mcN^{\sstable}(\nu) \sslash K_{\CC} \cong \Proj R_{\nu},
		\qquad R_{\nu} := \bigoplus\nolimits_{m \geq 0} \CC[\mcN]_{m\nu}.
	\end{equation}
\end{definition}

For integral $\nu$, genericity also gives equality of the stable and semistable loci: by \cite[Proposition~2.1]{Konno2000}, $\nu$ lies in no span of fewer than $N - n$ of the $\alpha_{i}$, so it is a regular value of $\phi_{\RR}$ on $T^{\ast}\CC^{N}$, and \cite[Theorem~2.3]{Proudfoot2004} applies.
An explicit description of the common locus $\mcN^{\sstable}(\nu)$ in the hypertoric coordinates is given in \cite[Lemma~3.4]{Konno2008}.

By \cite[Theorem~6.3]{HauselSturmfels2002}, see also \cite[Theorem~2.3]{Proudfoot2004} and \cite[Proposition~3.3]{Konno2008}, the inclusion of the hyperk\"ahler level set into the semistable locus induces the residual $T^{n}$-equivariant Kempf--Ness biholomorphism:
\begin{equation}
	\label{eqn:hypertoric-kempf-ness}
	( \phi_{\HK}^{-1}(\nu , 0)/K ,\, J_{1} ) \overset{\sim}{\lra} \mcN^{\sstable}(\nu) / K_{\CC}.
\end{equation}
We denote the common complex orbifold by $M_{\nu}$; since the stable and semistable loci coincide, $M_{\nu}^{\alg}$ is its coarse space.

The canonical morphism $M_{\nu}^{\alg} \ra \Spec R_{\nu, 0} = \Spec \CC[\mcN]^{K_{\CC}}$ is projective, and we call its target the affinisation.
Moreover, $\mcN$ is a reduced irreducible complete intersection because every $u_{i}$ is non-zero \cite[Lemma~4.7]{BellamyKuwabara2012}.
Its non-empty semistable locus and quotient are irreducible, and the quotient is normal because it has only finite quotient singularities \cite[Proposition~6.2]{HauselSturmfels2002}.
Thus, $M_{\nu}$ and $M_{\nu}^{\alg}$ are connected, irreducible, and normal.

Give $L_{\nu}$ the Gaussian metric $\|( ( z , w ) , v )\|^{2} = |v|^{2} \exp( -( \|z\|^{2} + \|w\|^{2} )/2 )$.
Its Chern connection has curvature $( 2\pi/\sqrt{-1} ) \w_{\RR}|_{\mcN}$, with our period-one convention; compare \cite[\S3.1]{Konno2008}.
Reduction gives the $T^{n}$-equivariant holomorphic pre-quantum line orbibundle
\begin{equation}
	\label{eqn:prequantum-orbibundle}
	\mcL_{\nu} := L_{\nu}|_{\phi_{\HK}^{-1}(\nu,0)}/K \longrightarrow M_{\nu},
	\qquad F_{\nabla_{\nu}} = ( 2\pi/\sqrt{-1} ) \w_{\RR},
\end{equation}
which agrees under \eqref{eqn:hypertoric-kempf-ness} with the descent of \eqref{eqn:positive-character-linearisation}.
It is relatively ample over the affinisation; a sufficiently divisible positive power descends to a relatively ample line bundle on the coarse quotient, whereas $\mcL_{\nu}$ itself need not descend.

\begin{notation}[Weight Conventions]
	\label{notation:weight-conventions}
	We use the geometric $+\chi_{\lambda}$ lift and the pullback action on sections.
	At a residual fixed point $p$, the geometric pre-quantum fibre weight is $-\mu_{\RR}(p)$, and the pullback fibre weight is $+\mu_{\RR}(p)$ \cite[p.~241, (7)]{Meinrenken1998}.
\end{notation}

\subsection{Arrangement, Core, and Cuts}
\label{subsection:subsection-2.3-arrangement-core-cuts}

Set $\ell_{i}(x) := \langle x ,\, u_{i} \rangle + \lambda_{i}$, $H_{i} := \{ \ell_{i} = 0 \}$, and $H_{i}^{\pm} := \{ \pm\ell_{i} \geq 0 \}$.
The arrangement $\mcA = \{ H_{1} , \ldots , H_{N} \}$ is simple \cite[Section~2]{BrownI}, and for each subset $A \subseteq \{1, \ldots, N\}$, the closed region $\Delta_{A}$ is defined by $\ell_{i} \leq 0$ for $i \in A$ and $\ell_{i} \geq 0$ otherwise.
We assume that a bounded region exists, and coorient the arrangement so that it is $\Delta_{\emptyset}$.
Write $\mcI_{x} := \{ i ~|~ \ell_{i}(x) = 0 \}$, $H_{\mcI_{x}} := \cap_{i \in \mcI_{x}} H_{i}$, and $\rk(\mcI_{x}) := \codim H_{\mcI_{x}}$; the arrangement vertices $\Verts(\mcA)$ are the points with $\rk(\mcI_{x}) = n$.

The extended core $\mcE_{\nu} := \mu_{\CC}^{-1}(0)$ has toric components
\[
	\mcE_{A} = \{ [ z , w ] \in \mcE_{\nu} ~|~ w_{i} = 0 \text{ if } i \not\in A,~z_{i} = 0 \text{ if } i \in A \}
\]
indexed by the non-empty regions, each a symplectic toric orbifold for the residual $T^{n}$-action \cite[Theorem~6.5(iii)]{BielawskiDancer2000} and \cite[\S2]{HaradaProudfoot2004}.
Each restriction $\mu_{\RR} : \mcE_{A} \ra \Delta_{A}$ is a proper surjection \cite[Lemma~3.3]{BrownI}.
The union of the components with bounded regions is the compact core $\mcC$, and $X_{\nu} := \mcE_{\emptyset}$ is its distinguished compact toric component \cite[\S2]{HaradaProudfoot2004}.
The fibre-rotation action $\tau \cdot [ z , w ] = [ z , \tau w ]$, where $\tau \in U_{1}$, preserves $\w_{\RR}$, scales $\w_{\CC}$ by $\tau$, and has a proper moment map $\Phi[z , w] = \tfrac{1}{2}\|w\|^{2}$ with zero set $X_{\nu}$ \cite[Proposition~1.3 and the following remark]{HaradaProudfoot2004}.
On $\mcE_{A}$, this action is the residual one-parameter subgroup with generator $u_{A}$, and $\Phi$ restricts as follows:
\begin{equation}
	\label{eqn:U1-action-generator}
	u_{A} := -\sum\nolimits_{i \in A} u_{i},
	\qquad \lambda_{A} := \sum\nolimits_{i \in A} \lambda_{i},
	\qquad \Phi|_{\mcE_{A}} = \langle \mu_{\RR} ,\, u_{A} \rangle - \lambda_{A}.
\end{equation}
The fibre-rotation action extends to $\CC^{\ast}$; since $\Phi$ is proper and bounded below, its fixed locus is compact and $\lim_{\tau \ra 0} \tau \cdot p$ exists for every $p$, so $M_{\nu}^{\alg}$ is semiprojective in the sense of \cite[Definition~1.1.1]{HauselRodriguezVillegas2015}.

On $M_{\nu} \times \CC_{\zeta}$, use the diagonal action $\tau \cdot ( p , \zeta ) = ( \tau \cdot p , \tau \zeta )$ and the product K\"ahler form $\w_{\RR} + \w_{\aux}$, where $\w_{\aux} = ( \sqrt{-1}/4\pi ) d\zeta \wedge d\bar{\zeta}$.
Its moment map is $\rho(p , \zeta) = \Phi(p) + \tfrac{1}{2} |\zeta|^{2}$, with $\CritVal(\rho) = \CritVal(\Phi)$.
For a positive regular level $\delta$, the symplectic cut of Lerman \cite[\S1.1, p.~249]{Lerman1995}, in the orbifold formulation of \cite[\S4]{Meinrenken1998}, is
\[
	M_{\nu}^{\leq \delta} := \rho^{-1}(\delta)/U_{1},
	\qquad Z_{\nu}^{\delta} := \Phi^{-1}(\delta)/U_{1}.
\]
The cut is a K\"ahler orbifold of complex dimension $2n$ carrying the residual $T^{n}$-action; $Z_{\nu}^{\delta}$ is its cut divisor, and the complement $M_{\nu}^{< \delta}$ of $Z_{\nu}^{\delta}$ is equivariantly symplectomorphic to $\{ \Phi < \delta \}$.
We also use $\mu_{\RR}$ and $\Phi$ for the functions descended to the cut.

\begin{lemma}[{\cite[Lemma~4.1]{BrownI}}]
	\label{lemma:proper-moment-map}
	The cutting moment map $\rho$ is proper, so $M_{\nu}^{\leq \delta}$ is compact at every positive regular level $\delta$.
\end{lemma}

\subsection{Slack and Moment Polyptychs}
\label{subsection:subsection-2.4-slack-moment-polyptychs}

The piecewise-linear slack function is
\begin{equation}
	\label{eqn:slack-function}
	\mcS : (\mft^{n})^{\ast} \lra \RR_{\geq 0}, \qquad
	\mcS(x) := \sum_{i = 1}^{N} \max\{ 0 , -\ell_{i}(x) \}.
\end{equation}
It satisfies $\mcS( \mu_{\RR}(p) ) \leq \Phi(p)$, with equality if and only if $p \in \mcE_{\nu}$ \cite[Lemma~4.3]{BrownI}.
On $\Delta_{A}$, it is given by $\mcS(x) = \langle x ,\, u_{A} \rangle - \lambda_{A}$ \cite[Lemma~4.2]{BrownI}.
Cutting $\mcE_{A}$ gives a closed toric suborbifold $\mcE_{A}^{\leq \delta} \subset M_{\nu}^{\leq \delta}$ with moment image $\Delta_{A}^{\leq \delta} := \Delta_{A} \cap \{ \mcS \leq \delta \}$ and cut face $\Pi_{A}^{\delta} := \Delta_{A} \cap \{ \mcS = \delta \}$ \cite[Proposition~4.4]{BrownI}.
The non-empty truncations $\Delta_{A}^{\leq \delta}$ are compact $n$-polytopes, rational when $\lambda$ and $\delta$ are rational, and their faces form the moment polyptych $\Delta_{\nu}^{\leq \delta}$.
Its boundary subcomplex $\Pi_{\nu}^{\delta}$ consists of the cut faces and their faces.
The supports of these complexes are
\[
	|\Delta_{\nu}^{\leq \delta}| =
	\bigcup\nolimits_{A} \Delta_{A}^{\leq \delta}
	= \mu_{\RR}(M_{\nu}^{\leq \delta}) = \{ \mcS \leq \delta \},
	\qquad |\Pi_{\nu}^{\delta}| = \{ \mcS = \delta \},
\]
where empty truncations contribute nothing \cite[Proposition~4.6]{BrownI}.
We call a vertex of $\Delta_{\nu}^{\leq \delta}$ a boundary vertex if it lies in $\Pi_{\nu}^{\delta}$, and an interior vertex otherwise.
We also put $\mcE_{\nu}^{\leq \delta} := \cup_{A} \mcE_{A}^{\leq \delta}$.

The critical values are $\CritVal(\Phi) = \{ \mcS(v) ~|~ v \in \Verts(\mcA) \}$ \cite[Proposition~4.7]{BrownI}.
The cut chambers are the components of $( 0 , \infty ) \setminus \CritVal(\Phi)$.
The outer chamber is
\begin{equation}
	\label{eqn:outer-cut-chamber}
	\Cout = ( c_{\max} , \infty ), \qquad c_{\max} := \max \CritVal(\Phi).
\end{equation}
The critical-value set is finite, and $\Phi \leq c_{\max}$ on $\mcC$; hence an outer-chamber cut divisor is disjoint from the compact core.

\subsection{Fixed-Point and Lattice Data}
\label{subsection:subsection-2.5-fixed-point-lattice-data}

At every positive regular level, $\mu_{\RR}$ identifies the residual fixed points with the vertices of $\Delta_{\nu}^{\leq \delta}$, preserving the interior/boundary distinction \cite[Lemma~5.1]{BrownI}.
For the fixed point $p$ over $v$, the set $\mcI_{v}$ has $n$ elements at an interior vertex and $n - 1$ at a boundary vertex.
For each $j \in \mcI_{v}$, there are two incident non-cutting edges, $\varepsilon_{v, j}$ in the direction where $\ell_{j}$ increases and $\varphi_{v, j}$ where it decreases; a boundary vertex has one further cutting edge $\varkappa_{v}$, tangent to $H_{\mcI_{v}}$ and directed into the cut \cite[Theorem~5.2 and Notation~5.3]{BrownI}.
The tangent weights along these edges are the rational covectors $\varrho_{p, j} \in T_{v}\varepsilon_{v, j}$ and $\varsigma_{p, j} \in T_{v}\varphi_{v, j}$ normalised by
\begin{equation}
	\label{eqn:noncutting-weight-normalisation}
	\langle \varrho_{p , j} ,\, u_{j} \rangle = 1, \qquad
	\langle \varsigma_{p , j} ,\, u_{j} \rangle = -1.
\end{equation}
At a boundary vertex, choose a cut face $\Pi_{A}^{\delta}$ through $v$ and normalise the cutting weight $\vartheta_{p} \in T_{v}\varkappa_{v}$ by
\begin{equation}
	\label{eqn:cutting-weight-normalisation}
	\langle \vartheta_{p} ,\, u_{j} \rangle = 0 \quad ( j \in \mcI_{v} ),
	\qquad \langle \vartheta_{p} ,\, u_{A} \rangle = -1.
\end{equation}
These conditions determine the weights uniquely.

For $I \subseteq \{ 1 , \ldots , N \}$, let $U_{I} := \sum_{i \in I} \ZZ u_{i}$ denote the integer span, without taking its saturation.
The local normal lattice is
\begin{equation}
	\label{eqn:local-normal-lattice}
	\Lambda_{v} := \begin{cases}
		U_{\mcI_{v}},                             & \text{for } v \text{ an interior vertex},                      \\
		U_{\mcI_{v}} + \ZZ u_{A},     & \text{for } v \in \Pi_{A}^{\delta} \text{ a boundary vertex}.
	\end{cases}
\end{equation}
The full labelled normal $u_{A}$ is used, even when it is non-primitive.
If two cut faces meet at $v$, their generators differ by an element of $U_{\mcI_{v}}$ \cite[Lemma~4.5]{BrownI}; thus both $\Lambda_{v}$ and $\vartheta_{p}$ are independent of the chosen face.
The lattice $\Lambda_{v}$ has full rank, and the local group and the local covering torus are
\begin{equation}
	\label{eqn:local-orbifold-group}
	\Gamma_{p} \cong \mft_{\ZZ}^{n}/\Lambda_{v},
	\qquad \widetilde{T}_{v}^{n} := \mft^{n}/\Lambda_{v} \longrightarrow T^{n}.
\end{equation}
The character lattice of $\widetilde{T}_{v}^{n}$ is
\begin{equation}
	\label{eqn:local-dual-lattice}
	\Lambda_{v}^{\vee} := \{ \alpha \in ( \mft^{n} )^{\ast} \mid
	\langle \alpha , \eta \rangle \in \ZZ \text{ for all } \eta \in \Lambda_{v} \}.
\end{equation}
These are the full local groups of the cut, inherited by each toric cut component through $p$; see \cite[Section~5]{BrownI} and \cite[Lemma~6.6]{LermanTolman1997}.

\begin{proposition}[{\cite[Proposition~5.5 and the paragraph after Lemma~5.6]{BrownI}}]
	\label{prop:isotropy-weights}
	The weights above lie in $\Lambda_{v}^{\vee}$, and the tangent isotropy representation of the local covering torus is
	\[
		T_{p}M_{\nu}^{\leq \delta} \cong
		\begin{cases}
			\displaystyle \bigoplus\nolimits_{j \in \mcI_{v}} ( V_{\varrho_{p, j}} \oplus V_{\varsigma_{p,j}} ), & \text{if } p \in M_{\nu}^{< \delta}, \\[1mm]
			\displaystyle \bigoplus\nolimits_{j \in \mcI_{v}} ( V_{\varrho_{p, j}} \oplus V_{\varsigma_{p,j}} )
			\oplus V_{\vartheta_{p}}^{\oplus 2}, & \text{if } p \in Z_{\nu}^{\delta}.
		\end{cases}
	\]
\end{proposition}
In particular, the cutting edge has one geometric direction, but its weight occurs twice in the full tangent space.

\begin{lemma}[Local Group Action on a Weight Space {\cite[Lemma~5.6]{BrownI}}]
	\label{lemma:local-cover-kernel-character}
	Let $\Gamma_{p} \cong \mft_{\ZZ}^{n} / \Lambda_{v}$ be the deck group of the covering $\widetilde{T}_{v}^{n} \ra T^{n}$.
	If $g \in \Gamma_{p}$ is represented by $\gamma_{g} \in \mft_{\ZZ}^{n}$, then its geometric eigenvalue on a weight line $V_{\alpha}$, for $\alpha \in \Lambda_{v}^{\vee}$, is $\exp( 2\pi\sqrt{-1}\langle \alpha ,\, \gamma_{g} \rangle )$.
	This is independent of the representative.
\end{lemma}
Thus, the geometric tangent eigenvalues used below are
\begin{equation}
	\label{eqn:group-characters}
	\rho_{p , j}(g) := e^{2\pi\sqrt{-1} \langle \varrho_{p , j} ,\, \gamma_{g} \rangle}, \qquad
	\sigma_{p , j}(g) := e^{2\pi\sqrt{-1} \langle \varsigma_{p , j} ,\, \gamma_{g} \rangle},
\end{equation}
and, at a boundary fixed point,
\begin{equation}
	\label{eqn:group-character-boundary}
	\theta_{p}(g) := e^{2\pi\sqrt{-1} \langle \vartheta_{p} ,\, \gamma_{g} \rangle}
\end{equation}
on each of the two cutting-weight summands.

%% file: sections/3-localisation-formulae.tex
% !TeX root = ../main.tex
\section{Equivariant Localisation Formulae}
\label{section:section-3-localisation-formulae}

We now turn from the geometry of the compact cuts to the calculation of their equivariant indices.
The equivariant orbifold localisation formula reduces this global calculation to one involving local data at the fixed points: their tangent representations, finite orbifold structure groups, and the actions on the fibres of the pre-quantum line orbibundle.
We first fix the equivariant conventions and recall the localisation formula, and then determine the fibre data of the cut, which complete the local inputs recalled in Subsection~\ref{subsection:subsection-2.5-fixed-point-lattice-data}.

\subsection{Equivariant Localisation Formulae}
\label{subsection:subsection-3.1-equivariant-localisation-formulae}

Let $M$ be a compact K\"ahler manifold with K\"ahler form $\w$, and let $\mcL \ra M$ be a holomorphic pre-quantum line bundle, normalised so that $F_{\nabla} = ( 2 \pi / \sqrt{-1} )\, \w$ and $c_{1}(\mcL) = [ \w ]$.
Write $\chi(M ; \mcL) = \sum_{q \geq 0} ( -1 )^{q} \dim_{\CC} H^{q}(M ; \mcL)$ for its holomorphic Euler characteristic, or Riemann--Roch number.
The Hirzebruch--Riemann--Roch formula expresses this number as $\chi(M ; \mcL) = \int_{M} \Td(TM) \wedge \Ch(\mcL)$ and, when the higher cohomology vanishes, the alternating sum reduces to $\chi(M ; \mcL) = \dim_{\CC}H^{0}(M ; \mcL)$.

Suppose that $\dim_{\CC}M = D$, and that a compact torus $T$ acts holomorphically on $M$ with a moment map $\mu : M \ra \mft^{\ast}$ and a lift to $\mcL$ preserving the pre-quantum connection.
Since $T$ is compact and connected, this action also preserves the Hermitian structure of $\mcL$.

We use the geometric linearisation fixed in Section~\ref{section:section-2-conventions-geometric-input} and the pullback action on sections given in Notation~\ref{notation:representation-convention}, i.e., a geometric fibre character is inverted when regarded as a character in the induced representation on sections.
In the equivariant formulae below, we use a logarithmic variable $\xi \in \mft_{\CC}$: the element $e^{\xi} \in T_{\CC}$ is characterised by $\chi_{\alpha}(e^{\xi}) = e^{\langle \alpha ,\, \xi \rangle}$ for every $\alpha \in \mft_{\ZZ}^{\ast}$.
Thus, for $\eta \in \mft$, the element $[ \eta ] \in T$ of Section~\ref{section:section-2-conventions-geometric-input} is written as $e^{2\pi\sqrt{-1}\eta}$ in this notation.

For a local tangent line bundle $V_{\alpha}$ of geometric weight $\alpha$, the equivariant first Chern class is
\begin{equation}
	\label{eqn:equivariant-first-Chern-class-convention}
	c_{1}^{T}(V_{\alpha})(\xi)
	= -\langle \alpha ,\, \xi \rangle,
\end{equation}
so that the corresponding factor of the fixed-point denominator is $1 - e^{\langle \alpha , \xi \rangle}$.
For a pre-quantum lift normalised by the moment map $\mu$, the geometric action on the fibre over a fixed point $p$ has weight $-\mu(p)$.
We write $\beta_{p} = \mu(p)$ for the corresponding pullback fibre weight, which enters the equivariant character.

Each $H^{q}(M ; \mcL)$ is a finite-dimensional representation of $T$.
We define the \textbf{equivariant holomorphic Euler characteristic} by
\begin{equation}
	\label{eqn:holomorphic-Euler-characteristic}
	\chi_{T}(M ; \mcL) := \sum\nolimits_{q \geq 0}( -1 )^{q} [ H^{q}(M ; \mcL) ] \in R(T).
\end{equation}
Suppose that $M^{T}$ is finite and write, at each fixed point $p$, $T_{p}M \cong \bigoplus_{j = 1}^{D} V_{\alpha_{p , j}}$ with $\alpha_{p , j} \neq 0$.
The equivariant index theorem expresses the character of $\chi_{T}(M ; \mcL)$ as the integral of the equivariant Riemann--Roch class $\Td_{T}(TM) \wedge \Ch_{T}(\mcL)$.
Atiyah--Bott--Berline--Vergne localisation evaluates this integral as a sum over the fixed points, dividing each contribution by the equivariant Euler class of $T_{p}M$.
With \eqref{eqn:equivariant-first-Chern-class-convention}, this division cancels the linear numerators in the Todd class, leaving one factor $( 1 - e^{\langle \alpha_{p, j} , \xi \rangle} )^{-1}$ for each tangent line bundle, and the equivariant Chern character of the pre-quantum fibre contributes $e^{\langle \mu(p) , \xi \rangle}$.
We therefore obtain the holomorphic fixed-point formula; compare \cite[Corollary~I.10]{GuilleminGinzburgKarshon2002}:
\begin{equation}
	\label{eqn:equivariant-HRR}
	\chi_{T}(M ; \mcL)(e^{\xi}) = \sum_{p \in M^{T}}
	\frac{
		e^{\langle \mu(p) ,\, \xi \rangle}
	}{
		\prod_{j = 1}^{D}
		\left( 1 - e^{\langle \alpha_{p , j} ,\, \xi \rangle} \right)
	},
\end{equation}
for $\xi \in \mft_{\CC}$ away from the denominator hyperplanes where $e^{\langle \alpha_{p , j} , \xi \rangle} = 1$.
Although the individual summands may have poles, their sum is the character of a finite-dimensional virtual representation and is regular at the identity.
Furthermore, if the higher cohomology vanishes, then this character is the character of $H^{0}(M ; \mcL)$.

Of course, the cut spaces considered below need not be smooth manifolds.
Even when $M_{\nu}$ is smooth, \cite[Theorem~4.18]{BrownI} shows that a reducible core $\mcC$ forces every cut $M_{\nu}^{\leq \delta}$ with $\delta \in \Cout$ to have non-trivial cyclic isotropy at the cut boundary and a singular coarse local quotient there.
We therefore need the orbifold version of the previous formula.
At an orbifold fixed point $p$, the $T$-action on a local uniformising chart may lift only to a finite covering torus $\widetilde{T}_{p} \ra T$.
This makes the local tangent and fibre weights characters of $\widetilde{T}_{p}$; that is, they may be regarded as rational weights with respect to the global lattice $\mft_{\ZZ}^{\ast}$.

\begin{theorem}[Equivariant Orbifold Kawasaki--Riemann--Roch Formula]
	\label{thm:equivariant-KRR-index-theorem}
	Let $M$ be a compact K\"ahler orbifold of complex dimension $D$, equipped with a holomorphic action of a torus $T$, and let $\mcL \ra M$ be a $T$-equivariant holomorphic line orbibundle.
	Suppose that the fixed-point set $M^{T}$ is finite.

	For each $p \in M^{T}$, let $\Gamma_{p}$ be the local orbifold structure group.
	Assume that, on a local uniformising chart, the $T$-action lifts to a finite covering torus $\widetilde{T}_{p} \ra T$ whose action commutes with that of $\Gamma_{p}$, and that the tangent representation splits into common $\widetilde{T}_{p}$- and $\Gamma_{p}$-weight spaces $T_{p}M \cong \bigoplus_{j = 1}^{D} V_{\alpha_{p , j}}$ with $\alpha_{p , j} \neq 0$.
	For $g \in \Gamma_{p}$, let $\chi_{p , j}^{\tang}(g) \in U_{1}$ be the \emph{geometric} eigenvalue of $g$ on the tangent line bundle $V_{\alpha_{p , j}}$, let $\beta_{p}$ be the pullback torus weight of the local fibre, and let $\chi_{p}^{\fib}(g) \in U_{1}$ be its local \emph{pullback} character; i.e., $\chi_{p}^{\fib}(g)$ is the inverse of the geometric eigenvalue of $g$ on the local fibre $\mcL_{p}$.
	We say that $\xi \in \mft_{\CC}$ is \textbf{generic} if
	\begin{equation}
		\label{eqn:orbifold-generic-xi}
		\chi_{p, j}^{\tang}(g) e^{\langle \alpha_{p , j} ,\, \xi \rangle} \neq 1
	\end{equation}
	for every $p \in M^{T}$, $g \in \Gamma_{p}$, and $j = 1 , \ldots , D$.

	Then, for all generic $\xi \in \mft_{\CC}$, the \emph{Kawasaki localisation formula} gives
	\begin{equation}
		\label{eqn:equivariant-KRR-character-formula}
		\chi_{T}(M ; \mcL)(e^{\xi}) = \sum_{p \in M^{T}}
		\frac{1}{|\Gamma_{p}|}
		\sum_{g \in \Gamma_{p}}
		\frac{
			\chi_{p}^{\fib}(g) e^{\langle \beta_{p} ,\, \xi \rangle}
		}{
			\prod\nolimits_{j = 1}^{D}
			( 1 - \chi_{p, j}^{\tang}(g) e^{\langle \alpha_{p , j} ,\, \xi \rangle} )
		}.
	\end{equation}
\end{theorem}

Theorem~\ref{thm:equivariant-KRR-index-theorem} is the isolated-fixed-point specialisation of the equivariant orbifold index theorem of Vergne \cite{Vergne1996}, which builds on Kawasaki's orbifold Riemann--Roch theorem \cite{Kawasaki1979}; see also \cite[Theorems~3.4--3.5]{Meinrenken1998}.
The finite-group average over the local structure group is written explicitly in the isolated-fixed-point case in \cite[\S10.1 and Remark~10.4]{CannasDaSilva1996}; compare also \cite[Corollary~I.10 and (I.13)]{GuilleminGinzburgKarshon2002}.
The cited formulae are initially stated either near the identity or under a stronger genericity hypothesis; both sides are meromorphic functions of $\xi$, so the equality extends by meromorphic continuation to every $\xi$ satisfying \eqref{eqn:orbifold-generic-xi}.
We call the equivariant index character $\chi_{T}(M ; \mcL)(e^{\xi})$ computed by \eqref{eqn:equivariant-KRR-character-formula} the \textbf{Kawasaki character}.

For the cut space of a hypertoric variety, $\Gamma_{p}$ is the kernel of the local covering homomorphism $\widetilde{T}_{v}^{n} \ra T^{n}$, and its action on the chart is the restriction of the covering-torus action.
The decomposition into one-dimensional $\widetilde{T}_{v}^{n}$-weight spaces is therefore also a decomposition into common $\widetilde{T}_{v}^{n}$- and $\Gamma_{p}$-eigenspaces, as required in Theorem~\ref{thm:equivariant-KRR-index-theorem}.
When applying the Kawasaki localisation formula to $M_{\nu}^{\leq \delta}$, we use the full local group $\Gamma_{p}$ at each fixed point and its representations on the full tangent space and the pre-quantum fibre of the cut.
On the compact cut spaces, algebraic and holomorphic orbifold sheaf cohomology agree, so we use $H^{q}(M ; \mcL)$ for both of these groups.

The eigenvalues \eqref{eqn:group-characters} and \eqref{eqn:group-character-boundary}, with their multiplicities, are the geometric tangent eigenvalues $\chi_{p , j}^{\tang}(g)$ that appear in Theorem~\ref{thm:equivariant-KRR-index-theorem}.
Each one-dimensional summand in Proposition~\ref{prop:isotropy-weights} contributes a factor $1 - \chi_{p , j}^{\tang}(g) e^{\langle \alpha_{p , j} ,\, \xi \rangle}$ to the denominator, so the two copies of $V_{\vartheta_{p}}$ at a boundary fixed point contribute $( 1 - \theta_{p}(g) e^{\langle \vartheta_{p} ,\, \xi \rangle} )^{2}$.
In particular, the genericity condition \eqref{eqn:orbifold-generic-xi} is stronger than $e^{\langle \alpha_{p , j} ,\, \xi \rangle} \neq 1$.

\subsection{Pre-Quantum Line Orbibundles}
\label{subsection:subsection-3.2-pre-quantum-line-orbibundles}

To complete the local data for the Kawasaki localisation formula, we now construct the $T^{n}$-equivariant pre-quantum line orbibundle on the cut and determine its fibre characters.
Fix a positive regular integral level $d$, and retain the integral lift $\lambda$ from Section~\ref{section:section-2-conventions-geometric-input}, together with the pre-quantum line orbibundle $\mcL_{\nu} \ra M_{\nu}$ of \eqref{eqn:prequantum-orbibundle}, whose Chern curvature is $F_{\nabla_{\nu}} = ( 2\pi / \sqrt{-1} ) \w_{\RR}$.
The fibre-rotation action lifts to the trivial line bundle $L_{\nu} = \mcN \times \CC \ra \mcN$ via $\tau \cdot ( ( z , w ) , a ) = ( ( z , \tau w ) , a )$ for $\tau \in U_{1}$, which descends to $\mcL_{\nu} \ra M_{\nu}$ since it commutes with the $K$-action.

The auxiliary line $\CC_{\zeta}$ carries the K\"ahler form $\w_{\aux} = ( \sqrt{-1} / 4\pi ) d\zeta \wedge d\bar\zeta$ from Section~\ref{section:section-2-conventions-geometric-input}.
Since $d$ is integral, $\varpi_{d}(\tau) = \tau^{d}$ defines a character $\varpi_{d} : U_{1} \ra U_{1}$, which we use to twist the geometric fibre action on the auxiliary line bundle $\CC_{\zeta} \times \CC_{Y} \ra \CC_{\zeta}$.
We equip this line bundle with the Hermitian metric $\| ( \zeta , Y ) \|^{2} = |Y|^{2} \exp( -|\zeta|^{2} / 2 )$ and the lift $\tau \cdot ( \zeta , Y ) = ( \tau \zeta , \tau^{d} Y )$ for $\tau \in U_{1}$.
We denote the resulting pre-quantum line bundle by $\mcO_{\CC}(\varpi_{d}) \ra \CC_{\zeta}$; its Chern connection $\nabla_{\aux}$ has curvature $F_{\nabla_{\aux}} = ( 2\pi / \sqrt{-1} ) \w_{\aux}$.

Taking the product of these lifts gives the diagonal cutting-circle action on $\mcL_{\nu} \boxtimes \mcO_{\CC}(\varpi_{d}) \ra M_{\nu} \times \CC_{\zeta}$, which commutes with the residual $T^{n}$-action.
The untwisted product lift is compatible with the cutting moment map $\rho : M_{\nu} \times \CC_{\zeta} \ra \RR$, given by $\rho( [ z , w ] , \zeta ) = \Phi[z , w] + \tfrac{1}{2}|\zeta|^{2}$.
Twisting the auxiliary fibre action by $\varpi_{d}$ changes the compatible moment-map normalisation from $\rho$ to $\rho - d$ by the Kostant formula.

Since $( \rho - d )|_{\rho^{-1}(d)} = 0$, the product connection restricted to this level set is basic for the cutting circle and descends under reduction.
Furthermore, the regularity of $d$ ensures that the cutting $U_{1}$-action is locally free on $\rho^{-1}(d)$, so the quotient line bundle
\begin{equation}
	\label{eqn:cut-prequantum-orbibundle}
	\mcL_{\nu}^{\leq d} := \left.
	\left( \left. \left( \mcL_{\nu} \boxtimes \mcO_{\CC}(\varpi_{d}) \right) \middle|_{\rho^{-1}(d)} \right. \right) \middle/ U_{1}
	\right. \longrightarrow M_{\nu}^{\leq d}
\end{equation}
is a $T^{n}$-equivariant holomorphic pre-quantum line orbibundle.
The product Hermitian metric also descends, and the descended connection $\nabla^{\leq d}$ satisfies $F_{\nabla^{\leq d}} = ( 2\pi / \sqrt{-1} ) \w_{\RR}^{\leq d}$; compare the pre-quantum cutting constructions in \cite[\S3]{DGMW1995} and \cite[\S4]{Meinrenken1998}.

\begin{lemma}
	\label{lemma:cut-bundle-locality}

	Under the canonical $T^{n}$-equivariant symplectomorphism $M_{\nu}^{< d} \cong \{ q \in M_{\nu} ~|~ \Phi(q) < d \}$, there is a canonical $T^{n}$-equivariant isomorphism of Hermitian pre-quantum line orbibundles with connection $\mcL_{\nu}^{\leq d}|_{M_{\nu}^{< d}} \cong \mcL_{\nu}|_{\{ \Phi < d \}}$.
\end{lemma}

\begin{proof}
	This is the locality property of the pre-quantum symplectic-cut construction \cite[Theorem~4.5]{Meinrenken1998}.
\end{proof}

Let us now determine the local fibre characters of $\mcL_{\nu}^{\leq d}$ at a residual fixed point $p \in ( M_{\nu}^{\leq d} )^{T^{n}}$, with $v = \mu_{\RR}(p) \in ( \mft^{n} )^{\ast}$.
Recall from \eqref{eqn:local-dual-lattice} that $\Lambda_{v}^{\vee}$ is the character lattice of the local covering torus $\widetilde{T}_{v}^{n} = \mft^{n} / \Lambda_{v}$.
The integrality of $\lambda$ and $d$ gives $v \in \Lambda_{v}^{\vee}$.
Indeed, for every $j \in \mcI_{v}$, the equality $\ell_{j}(v) = 0$ gives $\langle v ,\, u_{j} \rangle = -\lambda_{j} \in \ZZ$.
At a boundary vertex $v \in \Pi_{A}^{d}$, the cutting equality also gives $\langle v ,\, u_{A} \rangle = d + \lambda_{A} \in \ZZ$.
The vectors $u_{j}$ for $j \in \mcI_{v}$, together with $u_{A}$ at a boundary vertex, generate the local normal lattice $\Lambda_{v}$ of \eqref{eqn:local-normal-lattice}, and hence $v$ pairs integrally with $\Lambda_{v}$.

Pre-quantum reduction preserves the normalisation of the residual moment map $\mu_{\RR}$, so the Kostant formula implies that the local covering torus $\widetilde{T}_{v}^{n}$ acts on the local fibre of $\mcL_{\nu}^{\leq d}$ at $p$ with geometric weight $-v$.
Under the pullback convention of Notation~\ref{notation:representation-convention}, the weight entering the Kawasaki localisation formula is therefore $\beta_{p} = v \in \Lambda_{v}^{\vee}$.
This applies to both interior and boundary fixed points; for interior ones, it also follows from Lemma~\ref{lemma:cut-bundle-locality}.

For $g \in \Gamma_{p}$ represented by $\gamma_{g} \in \mft_{\ZZ}^{n}$, Lemma~\ref{lemma:local-cover-kernel-character}, applied to the weight $-v$, shows that $g$ acts geometrically on this local fibre by multiplication by $\exp( -2 \pi \sqrt{-1} \langle v ,\, \gamma_{g} \rangle )$.
As before, under the pullback action on sections, the corresponding value of the fibre character in \eqref{eqn:equivariant-KRR-character-formula} is the inverse, namely
\begin{equation}
	\label{eqn:local-cut-fibre-character}
	\chi_{p}^{\fib}(g) = \exp( 2 \pi \sqrt{-1} \langle \mu_{\RR}(p) ,\, \gamma_{g} \rangle ),
\end{equation}
and this value is independent of the representative $\gamma_{g}$ because $v \in \Lambda_{v}^{\vee}$.
Together with the tangent decomposition in Proposition~\ref{prop:isotropy-weights} and the local group $\Gamma_{p}$ in \eqref{eqn:local-orbifold-group}, these fibre characters complete the local data required by Theorem~\ref{thm:equivariant-KRR-index-theorem} for $( M_{\nu}^{\leq d} , \mcL_{\nu}^{\leq d} )$.

%% file: sections/4-hypertoric-equivariant-quantisation.tex
% !TeX root = ../main.tex
\section{Equivariant Quantisation of Cut Hypertoric Varieties}
\label{section:section-4-equivariant-quantisation-cut-spaces}

We now compute the equivariant indices of the compact hypertoric cuts from the local data of Section~\ref{section:section-3-localisation-formulae}.
Although the full cut is not toric for the residual $T^{n}$-action, the moment polyptych determines its fixed points and the local representations that enter the orbifold Kawasaki--Riemann--Roch formula.

First, we assemble the localisation formula, which Subsection~\ref{subsection:subsection-6.2-worked-examples} evaluates in the examples of~\cite{BrownI}.
Section~\ref{section:section-5-algebraic-cutting-vanishing} relates these indices to the circle-graded quantisation of the original non-compact hypertoric variety and proves that the higher cohomology groups on the cuts vanish, which identifies these indices with characters of section spaces, and Section~\ref{section:section-6-combinatorial-quantisation} describes the resulting filtration combinatorially.

\subsection{Equivariant Localisation Formula}
\label{subsection:subsection-4.1-equivariant-localisation-formula}

At a positive regular integral level $d$, the cut space $M_{\nu}^{\leq d}$ is a compact K\"ahler orbifold of complex dimension $2n$, with K\"ahler form $\w_{\RR}^{\leq d}$.
It is equipped with the residual $T^{n}$-action and the $T^{n}$-equivariant pre-quantum line orbibundle $\mcL_{\nu}^{\leq d} \ra M_{\nu}^{\leq d}$.
We write $\chi_{T^{n}}(M_{\nu}^{\leq d} ; \mcL_{\nu}^{\leq d}) \in R(T^{n})$ for its equivariant holomorphic Euler characteristic, as in \eqref{eqn:holomorphic-Euler-characteristic}.
Although this is a priori a virtual $T^{n}$-representation, Theorem~\ref{theorem:cut-cohomology-every-level} identifies it with the representation on $H^{0}(M_{\nu}^{\leq d} ; \mcL_{\nu}^{\leq d})$.

The following theorem expresses $\chi_{T^{n}}(M_{\nu}^{\leq d} ; \mcL_{\nu}^{\leq d})$ in terms of the moment polyptych $\Delta_{\nu}^{\leq d}$ and its associated lattice data, combining the fixed-point, tangent, and fibre calculations of Sections~\ref{section:section-2-conventions-geometric-input} and \ref{section:section-3-localisation-formulae}.

\begin{theorem}[Equivariant Localisation Formula for Orbifold Cut Spaces]
	\label{theorem:character-hypertoric-orbifold}
	Let $M_{\nu}$ be the non-empty hypertoric orbifold constructed in Section~\ref{section:section-2-conventions-geometric-input} at a generic parameter $\nu$, with the arrangement cooriented so that the distinguished base region $\Delta_{\emptyset}$ is non-empty and bounded as in Subsection~\ref{subsection:subsection-2.3-arrangement-core-cuts}.
	Assume that the fixed lift $\lambda$ is integral, $\lambda \in ( \mft_{\ZZ}^{N} )^{\ast}$, so that $\nu = \imath_{\ZZ}^{\ast}(\lambda) \in \mfk_{\ZZ}^{\ast}$.
	Let $d$ be a positive regular integral level, let $\mcL_{\nu}^{\leq d} \to M_{\nu}^{\leq d}$ be the cut pre-quantum line orbibundle from \eqref{eqn:cut-prequantum-orbibundle}, and let $\xi \in \mft_{\CC}^{n}$ be generic in the sense of \eqref{eqn:orbifold-generic-xi}, so that every denominator factor displayed below is non-zero.

	Then the equivariant Kawasaki character is
	\begin{equation}
		\label{eqn:character-hypertoric-orbifold}
		\chi_{T^{n}}(M_{\nu}^{\leq d} ; \mcL_{\nu}^{\leq d})(e^{\xi}) = \sum\nolimits_{p \in ( M_{\nu}^{\leq d} )^{T^{n}}} \chi_{p}(e^{\xi}),
	\end{equation}
	where $\chi_{p}(e^{\xi})$ denotes the local orbifold contribution of the fixed point $p$.
	For a fixed point $p$, write $v := \mu_{\RR}(p)$ so that $\Gamma_{p} \cong \mft_{\ZZ}^{n} / \Lambda_{v}$ and, for each $g \in \Gamma_{p}$, let $\gamma_{g} \in \mft_{\ZZ}^{n}$ be a representative.
	If $p \in M_{\nu}^{< d}$ is an interior fixed point, then
	\begin{equation}
		\label{eqn:cut-localisation-interior}
		\chi_{p}(e^{\xi}) = \frac{1}{|\Gamma_{p}|}
		\sum_{g \in \Gamma_{p}}
		\frac{
			e^{2 \pi \sqrt{-1}\langle v ,\, \gamma_{g} \rangle}
			e^{\langle v ,\, \xi \rangle}
		}{
			\prod\nolimits_{j \in \mcI_{v}}
			\left( 1 - \rho_{p , j}(g) e^{\langle \varrho_{p , j} ,\, \xi \rangle} \right)
			\left( 1 - \sigma_{p , j}(g) e^{\langle \varsigma_{p , j} ,\, \xi \rangle} \right)
		}
	\end{equation}
	and, if $p \in Z_{\nu}^{d}$ is a boundary fixed point, then
	\begin{equation}
		\label{eqn:cut-localisation-boundary}
		\chi_{p}(e^{\xi}) = \frac{1}{|\Gamma_{p}|}
		\sum_{g \in \Gamma_{p}}
		\frac{
			e^{2 \pi \sqrt{-1} \langle v ,\, \gamma_{g} \rangle}
			e^{\langle v ,\, \xi \rangle}
		}{
			\left( 1 - \theta_{p}(g) e^{\langle \vartheta_{p} ,\, \xi \rangle} \right)^{2}
			\prod\nolimits_{j \in \mcIv}
			\left( 1 - \rho_{p , j}(g) e^{\langle \varrho_{p , j} ,\, \xi \rangle} \right)
			\left( 1 - \sigma_{p , j}(g) e^{\langle \varsigma_{p , j} ,\, \xi \rangle} \right)
		}.
	\end{equation}
	These expressions are independent of the choice of $\gamma_{g} \in \mft_{\ZZ}^{n}$.
	Here, $\rho_{p , j}(g)$, $\sigma_{p , j}(g)$, and $\theta_{p}(g)$ are the geometric tangent eigenvalues from \eqref{eqn:group-characters} and \eqref{eqn:group-character-boundary}.
	Moreover, under the pullback action on sections, the fibre factor $e^{2\pi\sqrt{-1} \langle v ,\, \gamma_{g} \rangle}$ is the character value $\chi_{p}^{\fib}(g)$ from \eqref{eqn:local-cut-fibre-character}.
\end{theorem}

\begin{proof}

	By \cite[Lemma~5.1]{BrownI}, the residual fixed points are isolated and correspond to the vertices $v$ of $\Delta_{\nu}^{\leq d}$.
	At the point $p$ above $v$, the local group is $\Gamma_{p} \cong \mft_{\ZZ}^{n}/\Lambda_{v}$ by \eqref{eqn:local-orbifold-group}.
	Lemma~\ref{lemma:local-cover-kernel-character} identifies $\Gamma_{p}$ with the kernel of the local covering $\widetilde{T}_{v}^{n} = \mft^{n} / \Lambda_{v} \ra T^{n}$.
	The action of $\Gamma_{p}$ on the local chart is obtained by restricting the covering-torus action, so each element of $\Gamma_{p}$ acts by a scalar on every isotropy weight space in Proposition~\ref{prop:isotropy-weights}, which is the decomposition required by Theorem~\ref{thm:equivariant-KRR-index-theorem}.

	The tangent summands are given in Proposition~\ref{prop:isotropy-weights}.
	Together with the geometric eigenvalues in \eqref{eqn:group-characters} and \eqref{eqn:group-character-boundary}, they determine the denominator factors.
	At a boundary fixed point, the summand $V_{\vartheta_{p}}^{\oplus2}$ gives the squared factor $( 1 - \theta_{p}(g) e^{\langle \vartheta_{p} ,\, \xi \rangle} )^{2}$.
	At every fixed point, the pre-quantum normalisation gives the pullback fibre weight $\beta_{p} = v = \mu_{\RR}(p)$, and \eqref{eqn:local-cut-fibre-character} gives $\chi_{p}^{\fib}(g) = e^{2\pi\sqrt{-1} \langle v,\, \gamma_{g} \rangle}$.

	Substituting these data into \eqref{eqn:equivariant-KRR-character-formula} proves \eqref{eqn:cut-localisation-interior} and \eqref{eqn:cut-localisation-boundary}.
	The expressions are independent of the representatives $\gamma_{g}$ because all the tangent and fibre weights involved lie in $\Lambda_{v}^{\vee}$.
\end{proof}

For the smooth cut spaces characterised by \cite[Corollary~4.20]{BrownI}, the local groups $\Gamma_{p}$ are trivial, and the previous theorem simplifies to the following localisation formula.

\begin{corollary}[Equivariant Localisation Formula for Smooth Cut Spaces]
	\label{corollary:character-smooth-hypertoric-cut}
	Under the hypotheses of Theorem~\ref{theorem:character-hypertoric-orbifold}, suppose that $M_{\nu}^{\leq d}$ is a smooth manifold with trivial local orbifold groups.
	Then $\mcL_{\nu}^{\leq d} \ra M_{\nu}^{\leq d}$ is an ordinary pre-quantum line bundle, and the contribution of a residual fixed point $p$, with $v = \mu_{\RR}(p)$, to \eqref{eqn:character-hypertoric-orbifold} is
	\begin{equation}
		\label{eqn:smooth-cut-local-contribution}
		\chi_{p}(e^{\xi}) =
		\begin{cases}
			\displaystyle
			\frac{
				e^{\langle v ,\, \xi \rangle}
			}{
				\prod\nolimits_{j \in \mcIv}
				\left( 1 - e^{\langle \varrho_{p , j} ,\, \xi \rangle} \right)
				\left( 1 - e^{\langle \varsigma_{p , j},\, \xi \rangle} \right)
			},
			\qquad & \text{if } p \in M_{\nu}^{< d}, \\[5mm]
			\displaystyle
			\frac{
				e^{\langle v ,\, \xi \rangle}
			}{
				\left( 1 - e^{\langle \vartheta_{p} ,\, \xi \rangle} \right)^{2}
				\prod\nolimits_{j \in \mcIv}
				\left( 1 - e^{\langle \varrho_{p , j},\, \xi \rangle} \right)
				\left( 1 - e^{\langle \varsigma_{p , j} ,\, \xi \rangle} \right)
			},
			\qquad & \text{if } p \in Z_{\nu}^{d}.
		\end{cases}
	\end{equation}
	Moreover, $\Lambda_{v} = \mft_{\ZZ}^{n}$ at every vertex, so the local covering torus is the residual torus $T^{n}$.
	In particular, this formula applies for every positive integral level $d$ whenever the equivalent conditions of \cite[Corollary~4.20]{BrownI} hold.
\end{corollary}

\begin{proof}
	Since $\Gamma_{p} = \{ 1 \}$, equation~\eqref{eqn:local-orbifold-group} gives $\Lambda_{v} = \mft_{\ZZ}^{n}$.
	The unique element of $\Gamma_{p}$ is represented by $0 \in \mft_{\ZZ}^{n}$, so all finite-group tangent and fibre character values in \eqref{eqn:cut-localisation-interior} and \eqref{eqn:cut-localisation-boundary} equal $1$, and the averages reduce to the expressions in \eqref{eqn:smooth-cut-local-contribution}.
	The final claim follows from \cite[Corollary~4.20]{BrownI}.
\end{proof}

We now return to the boundary contributions in the orbifold localisation formula \eqref{eqn:cut-localisation-boundary}.
When $n = 1$, the tangent representation at a boundary fixed point consists of two copies of the cutting weight, so each summand in the finite-group average has a squared denominator.
The local group is cyclic, of order the multiplicity $k$ of the labelled normal $u_{A}$.

Lemma~\ref{lemma:boundary-vertex-closed-form} evaluates this average as a rational function in the residual character variable $t \in T^{1}$, determined by $k$ and the coordinate $v$ of the boundary vertex.
This expression allows us to calculate the boundary contributions in the worked examples of Subsection~\ref{subsection:subsection-6.2-worked-examples} without repeating the finite-group summation.

To state the lemma, which uses Theorem~\ref{theorem:character-hypertoric-orbifold}, we continue to assume that the lift is integral and fix a positive regular integral level $d$.
For a boundary vertex $v \in \Pi_{A}^{d}$, write $p_{v} \in Z_{\nu}^{d}$ for the corresponding residual fixed point, with $\mu_{\RR}(p_{v}) = v$.
Henceforth, including in the worked examples, we index its local data by $v$, writing $\Gamma_{v}$, $\vartheta_{v}$, $\theta_{v}$, $\chi_{v}$, $\chi_{v}^{\fib}$, and $\chi_{v , j}^{\tang}$.

\begin{lemma}[Rank-One Boundary Contribution]
	\label{lemma:boundary-vertex-closed-form}
	Suppose that $n = 1$, and let $v \in \Pi_{A}^{d}$ be a boundary vertex.
	Choose the primitive lattice generator $\widebar{u} \in \mft_{\ZZ}^{1}$ such that $u_{A} = k\widebar{u}$, with $k \in \ZZ_{> 0}$.
	Identify $( \mft^{1} )^{\ast}$ with $\RR$ via the linear isomorphism $x \mapsto \langle x ,\, \widebar{u} \rangle$, and continue to write $v$ for the real coordinate $\langle v ,\, \widebar{u} \rangle$ of the vertex.
	Let $t = e^{\xi}$ be the $T^{1}$-character variable for this choice of generator.
	Write
	\begin{equation}
		\label{eqn:boundary-vertex-division}
		d + \lambda_{A} = ak + b, \qquad a , b \in \ZZ, \quad 0 \leq b < k,
	\end{equation}
	so that $v = a + b/k$ and $a = \lfloor v \rfloor$.
	Then $\Gamma_{v} \cong \ZZ / k\ZZ$, and
	\begin{equation}
		\label{eqn:boundary-vertex-closed-form}
		\chi_{v}(t) = \sum\nolimits_{\ell \geq 0} ( k\ell + b + 1 )\, t^{a - \ell} = t^{a}\, \frac{( b + 1 ) + ( k - b - 1 )\, t^{-1}}{( 1 - t^{-1} )^{2}}.
	\end{equation}
	The contribution $\chi_{v}(t)$ contains only integral powers of $t$: the fractional powers appearing in the individual terms cancel in the average.
\end{lemma}

\begin{proof}

	By \cite[Theorem~5.2(ii)]{BrownI}, we have $\mcI_{v} = \varnothing$, so $\Lambda_{v} = \ZZ u_{A}$ and $\Gamma_{v} = \mft_{\ZZ}^{1} / \ZZ u_{A} \cong \ZZ / k\ZZ$ is generated by the class of $\widebar{u}$.
	For $j = 0 , \ldots , k - 1$, let $g_{j} := [ j\widebar{u} ] \in \Gamma_{v}$ and choose $\gamma_{j} := j\widebar{u} \in \mft_{\ZZ}^{1}$ as its representative.
	The normalisation \eqref{eqn:cutting-weight-normalisation} gives $\langle \vartheta_{v} ,\, u_{A} \rangle = -1$, so $\vartheta_{v} = -1/k$ relative to the chosen coordinate, and hence $t^{\vartheta_{v}} = t^{-1/k}$.
	By \eqref{eqn:U1-action-generator}, $\lambda_{A} = \sum_{i \in A} \lambda_{i}$.
	Since the fixed lift $\lambda$ is integral, $\lambda_{A}$ is an integer, so \eqref{eqn:boundary-vertex-division} is a Euclidean division of integers.
	Moreover, $v \in \Pi_{A}^{d}$ gives $kv = \langle v ,\, u_{A} \rangle = d + \lambda_{A}$, so $v = a + b/k$.
	It follows that
	\[
		\theta_{v}(g_{j}) = e^{-2\pi\sqrt{-1}\, j/k}, \qquad
		\chi_{v}^{\fib}(g_{j}) = e^{2\pi\sqrt{-1} \langle v ,\, \gamma_{j} \rangle} = e^{2\pi\sqrt{-1}\, jb/k},
	\]
	so substituting these character values into the boundary localisation formula \eqref{eqn:cut-localisation-boundary} of Theorem~\ref{theorem:character-hypertoric-orbifold}, and using $\mcI_{v} = \varnothing$, gives
	\[
		\chi_{v}(t) = \frac{t^{v}}{|\Gamma_{v}|} \sum_{g \in \Gamma_{v}} \frac{\chi_{v}^{\fib}(g)}{( 1 - \theta_{v}(g) e^{\langle \vartheta_{v} , \xi \rangle} )^{2}} = \frac{t^{\, a + b/k}}{k} \sum_{j = 0}^{k - 1} \frac{e^{2 \pi \sqrt{-1}\, jb/k}}{\big( 1 - e^{-2 \pi \sqrt{-1}\, j/k}\, t^{-1/k} \big)^{2}}.
	\]
	Expanding $( 1 - z )^{-2} = \sum_{m \geq 0} ( m + 1 ) z^{m}$ and interchanging the finite sum with the formal power series, we find by character orthogonality for $\ZZ / k\ZZ$ that 	the inner average vanishes unless $m \equiv b \pmod{k}$.
	Writing the surviving indices as $m = b + k\ell$ with $\ell \geq 0$, and using $t^{\, a + b/k}\, t^{-( b + k\ell ) / k} = t^{\, a - \ell}$, gives the first equality in \eqref{eqn:boundary-vertex-closed-form}; summing $\sum_{\ell \geq 0} t^{-\ell}$ and $\sum_{\ell \geq 0} \ell\, t^{-\ell}$ gives the second.
\end{proof}

\begin{remark}
	\label{remark:boundary-vertex-specialisations}
	For $k = 1$, we have $b = 0$ and $a = v$, so \eqref{eqn:boundary-vertex-closed-form} reduces to $\chi_{v}(t) = t^{v}( 1 - t^{-1} )^{-2}$, the contribution of a boundary vertex with trivial local group.
	More generally, $b = 0$ if and only if $v$ is integral, that is, if and only if $\chi_{v}^{\fib}$ is trivial.
	The residue $b$ therefore measures the non-integrality of the vertex: the averaged contribution has integral powers for every $b$, with leading exponent $\lfloor v \rfloor$, equal to $v$ when $b = 0$.
\end{remark}

%% file: sections/5-algebraic-cutting-vanishing.tex
% !TeX root = ../main.tex
\section{Algebraic Cutting and Vanishing of Higher Cohomology}
\label{section:section-5-algebraic-cutting-vanishing}

At every positive regular integral level $d$, the cut quantisation admits three descriptions: the localisation index of Section~\ref{section:section-4-equivariant-quantisation-cut-spaces}, an algebraic description, and a combinatorial description.
Here, we consider the algebraic description: we compare sections with polynomial semi-invariants, realise the cut as a GIT quotient, identify its section space with a finite filtration piece of $H^{0}(M_{\nu} ; \mcL_{\nu})$, and then prove that its higher cohomology groups vanish, so that the index is the character of that space.
Section~\ref{section:section-6-combinatorial-quantisation} then computes the residual-weight multiplicities from the coordinate ring, obtaining a weighted lattice-point formula for the same equivariant character computed by the localisation formula of Theorem~\ref{theorem:character-hypertoric-orbifold}.

\subsection{Sections and Semi-Invariants}
\label{subsection:subsection-5.1-sections-semi-invariants}

We first compare the global sections of the tensor powers of $\mcL_{\nu}$ with polynomial semi-invariants.
The comparison is used throughout this section and the next; its proof uses the ambient Lawrence toric quotient and does not depend on the algebraic-cut construction.

\begin{proposition}[Semi-Invariant Comparison]
	\label{prop:positive-ray-semi-invariants}
	Let $P$ be the coordinate ring of $T^{\ast}\CC^{N}$, and let $J_{\CC}P$ be the ideal generated by the components of the complex moment map:
	\[
		P := \CC[z_{1} , \ldots , z_{N} , w_{1} , \ldots , w_{N}], \qquad
		J_{\CC}P := \langle \langle \phi_{\CC} ,\, \kappa \rangle ~|~ \kappa \in \mfk_{\CC} \rangle_{P}.
	\]
	The ideal $J_{\CC}P$ is the vanishing ideal $I( \phi_{\CC}^{-1}(0) )$, as $\mcN = \phi_{\CC}^{-1}(0)$ is a reduced complete intersection (Subsection~\ref{subsection:subsection-2.2-algebraic-quotients-prequantum}).
	Hence, $\CC[\mcN] \cong P / J_{\CC}P$.

	For every $m \geq 0$, let $(P / J_{\CC}P )_{m\nu}$ be the space of $\bar{f}$ with $k^{\ast} \bar{f} = \chi_{\nu}(k)^{m} \bar{f}$ for every $k \in K_{\CC}$; under the isomorphism above, it is the space $\CC[\mcN]_{m\nu}$ of Notation~\ref{notation:representation-convention}.
	Then there exists a natural $T^{n}$-equivariant isomorphism
	\begin{equation}
		\label{eqn:positive-ray-semi-invariants}
		( P / J_{\CC}P )_{m\nu} \overset{\sim}{\lra}
		H^{0}(M_{\nu} ; \mcL_{\nu}^{\otimes m}).
	\end{equation}
	Furthermore, these isomorphisms are compatible with multiplication, so they assemble into an isomorphism of graded $\CC$-algebras $\bigoplus\nolimits_{m \geq 0} ( P / J_{\CC}P )_{m\nu} \overset{\sim}{\lra} \bigoplus\nolimits_{m \geq 0} H^{0}(M_{\nu} ; \mcL_{\nu}^{\otimes m})$.
\end{proposition}

\begin{proof}
	The unstable locus in the complex moment fibre $\mcN$ may have codimension one, as in \cite[Example~1.2]{BellamyCrawSchedler2024}, so we cannot argue that sections extend directly across that locus.
	We therefore first prove that the complement of $( T^{\ast}\CC^{N} )^{\sstable}(\nu)$ in $T^{\ast}\CC^{N}$ has codimension at least two, and then work on the ambient Lawrence quotient.

	The Hilbert--Mumford criterion \cite[Lemmas~4.2--4.3]{BellamyKuwabara2012} describes instability in terms of the cone generated by the weights of the non-zero coordinates.
	Since the unstable locus is a union of coordinate subspaces, it suffices to show that the generic point of each coordinate hyperplane is semistable.
	We treat the hyperplane $\{ z_{i} = 0 \}$; the case $\{ w_{i} = 0 \}$ is identical.
	At its generic point, all coordinates $z_{j}$ and $w_{j}$ with $j \neq i$ are non-zero, so the weights $\pm\alpha_{j}$ for $j \neq i$ occur.
	The weights $\{ \alpha_{j} ~|~ j \neq i \}$ still span $\mfk^{\ast}$.
	Indeed, if $\kappa \in \mfk$ is annihilated by every $\alpha_{j}$ with $j \neq i$, then $\imath(\kappa) = ce_{i}$ for some $c \in \RR$.
	Since $\pi \circ \imath = 0$, we obtain $cu_{i} = 0$, which forces $c = 0$ because $u_{i} \neq 0$, and hence $\kappa = 0$.
	Since both signs $\pm\alpha_{j}$ occur among the remaining cotangent coordinates, these weights generate $\mfk^{\ast}$ as a convex cone.
	This cone contains the linearisation character $\nu$, so the generic point of each coordinate hyperplane is semistable.
	Hence, the unstable locus has codimension at least two.

	Let $W_{\nu} := [ (T^{\ast}\CC^{N})^{\sstable}(\nu) / K_{\CC} ]$ be the ambient Lawrence toric orbifold, with coarse space
	\[
		W_{\nu}^{\alg} := T^{\ast}\CC^{N} \sslash_{\nu} K_{\CC} = \Proj ( \oplus_{m \geq 0} P_{m\nu} ),
	\]
	and let $\mcA_{\nu} \ra W_{\nu}$ be the line orbibundle obtained by descending the trivial line bundle on $( T^{\ast}\CC^{N} )^{\sstable}(\nu)$, with geometric fibre action given by $\chi_{\nu}$.
	Recall that descent identifies sections of $\mcA_{\nu}^{\otimes m}$ with semi-invariants of character $m\nu$ on this semistable locus.
	Since its complement has codimension at least two in the normal affine space $T^{\ast}\CC^{N}$, those semi-invariants extend uniquely to polynomial semi-invariants, so $H^{0}(W_{\nu} ; \mcA_{\nu}^{\otimes m}) \cong P_{m\nu}$ for every $m \geq 0$.

	We also need the higher cohomology to vanish for these tensor powers on $W_{\nu}$.
	For each $m \geq 0$, the tensor power $\mcA_{\nu}^{\otimes m}$ is a line orbibundle on the canonical toric orbifold $W_{\nu}$, and its coarse pushforward is the rank-one reflexive sheaf $\mcO(E_{m})$, where $E_{m}$ is an integral torus-invariant Weil divisor.
	The divisor $E_{m}$ is $\QQ$-Cartier because the coarse Lawrence toric variety is simplicial.
	The orbifold is canonical because its coordinate divisors have trivial generic stabiliser: at a generic point of $\{ z_{i} = 0 \}$, every $w_{j}$ is non-zero, and $K_{\CC} \subseteq T_{\CC}^{N}$.
	Since the local stabilisers are finite and we work in characteristic zero, taking invariants is exact, so the coarse and orbifold cohomology groups agree.

	Let $\pi_{\aff} : W_{\nu}^{\alg} \ra \Spec P^{K_{\CC}}$ denote the projective toric morphism to the affinisation, which is surjective since it is projective and dominant.
	Under the description of $W_{\nu}^{\alg}$ as a $\Proj$, the sheaf $\mcO(E_{m})$ is the twisting sheaf $\mcO_{W_{\nu}^{\alg}}(m)$, and since $\oplus_{m \geq 0} P_{m\nu}$ is finitely generated over $P^{K_{\CC}}$, there is a $k > 0$ for which $\mcO(E_{k})$ is invertible and $\pi_{\aff}$-ample \cite[Proposition~4.6.18]{Grothendieck1961}.
	As $kE_{m}$ is linearly equivalent to $mE_{k}$, every $E_{m}$ with $m \geq 0$ is $\pi_{\aff}$-nef.
	Hence, Fujino's relative toric vanishing theorem \cite[Theorem~5.3]{Fujino2007} gives $R^{q}(\pi_{\aff})_{\ast} \mcO(E_{m}) = 0$ for every $q > 0$.
	Since $\Spec P^{K_{\CC}}$ is affine, the Leray spectral sequence and the previous cohomology comparison give $H^{q}(W_{\nu} ; \mcA_{\nu}^{\otimes m}) = 0$ for every $q > 0$ and $m \geq 0$.

	Reductivity of $K_{\CC}$ makes the restriction $P_{m\nu} \ra ( P/J_{\CC}P )_{m\nu}$ surjective for every $m \geq 0$, so $\mcN^{\sstable}(\nu) = \mcN \cap ( T^{\ast}\CC^{N} )^{\sstable}(\nu)$.
	The inclusion of these semistable loci induces $M_{\nu} \hookrightarrow W_{\nu}$.
	Furthermore, since $\mcA_{\nu}$ and $\mcL_{\nu}$ both descend via $\chi_{\nu}$, we have $\mcL_{\nu} = \mcA_{\nu}|_{M_{\nu}}$.

	Write $r := \dim_{\CC}K_{\CC} = N - n$.
	After choosing a basis of $\mfk_{\CC}$, let $f_{1} , \ldots , f_{r}$ be the corresponding components of $\phi_{\CC}$.
	Each $f_{a}$ is a linear combination of the products $z_{i}w_{i}$, and therefore has $K_{\CC}$-character zero.
	Hausel and Sturmfels show that these moment-map equations generate the homogeneous prime complete-intersection ideal denoted by $\operatorname{Circ}(B)$ in \cite[\S6]{HauselSturmfels2002}.
	The regular-sequence property is also recorded in \cite[Lemma~4.7]{BellamyKuwabara2012}, whose non-zero-row hypothesis is our condition $u_{i} \neq 0$.
	It follows that $f_{1} , \ldots , f_{r}$ is a regular sequence, and $M_{\nu} \hookrightarrow W_{\nu}$ is the complete intersection cut out by their descended sections.

	We regard the descended moment map as a section of the trivial vector bundle with fibre $\mfk_{\CC}^{\ast}$, so its Koszul complex uses the exterior powers of $( \mfk_{\CC}^{\ast} )^{\vee} = \mfk_{\CC}$.
	Tensoring the Koszul resolution of $\mcO_{M_{\nu}}$ with $\mcA_{\nu}^{\otimes m}$ then gives
	\[
	\begin{tikzcd}[ampersand replacement=\&]%, column sep=large]
		0 \ar[r] \&
		\mcA_{\nu}^{\otimes m} \otimes \bigwedge^{r} \mfk_{\CC} \ar[r] \&
		\ldots \ar[r] \&
		\mcA_{\nu}^{\otimes m} \otimes \mfk_{\CC} \ar[r] \&
		\mcA_{\nu}^{\otimes m} \ar[r] \&
		\mcL_{\nu}^{\otimes m} \ar[r] \&
		0.
	\end{tikzcd}
	\]
	All the terms involve the same tensor power $\mcA_{\nu}^{\otimes m}$ because each $f_{a}$ has trivial $K_{\CC}$-character.
	The previous toric vanishing therefore makes the global-section functor exact on this resolution.
	Hence, the restriction map $H^{0}(W_{\nu} ; \mcA_{\nu}^{\otimes m}) \ra H^{0}(M_{\nu} ; \mcL_{\nu}^{\otimes m})$ is surjective with kernel $( J_{\CC}P )_{m\nu}$, so $H^{0}(M_{\nu} ; \mcL_{\nu}^{\otimes m}) \cong P_{m\nu} / ( J_{\CC}P )_{m\nu} = ( P / J_{\CC}P )_{m\nu}$.

	The restriction and descent maps commute with the residual $T^{n}$-action and with the fibre-rotation action, taken with the lift from Subsection~\ref{subsection:subsection-3.2-pre-quantum-line-orbibundles} that leaves the fibre coordinate of $L_{\nu} = \mcN \times \CC$ unchanged.
	Thus, \eqref{eqn:positive-ray-semi-invariants} is $T^{n}$-equivariant and preserves fibre-rotation weights.
	Moreover, descent and restriction preserve multiplication: the product of two semi-invariants of characters $m\nu$ and $m^{\prime}\nu$ corresponds to the product of the corresponding sections of $\mcL_{\nu}^{\otimes m}$ and $\mcL_{\nu}^{\otimes m^{\prime}}$.
	The isomorphisms therefore assemble into the stated graded-algebra isomorphism.
\end{proof}

\begin{remark}
	Proposition~\ref{prop:positive-ray-semi-invariants} concerns the characters $m\nu$ with $m \in \ZZ_{\geq 0}$.
	For other characters of $K_{\CC}$, the semi-invariants on the semistable locus need not extend to polynomial semi-invariants on the level set $\mcN = \phi_{\CC}^{-1}(0)$.
	Indeed, already for the hypertoric presentation of $T^{\ast}\CC\PP^{1}$, the unstable locus can have codimension one, and its removal can enlarge the semi-invariant ring; see \cite[Example~1.2]{BellamyCrawSchedler2024}.
\end{remark}

\subsection{Algebraic Cutting}
\label{subsection:subsection-5.2-algebraic-cutting}

The localisation formula computes the equivariant index of each compact cut.
To relate these cuts to the quantisation of the original non-compact hypertoric variety $M_{\nu}$, we decompose $H^{0}(M_{\nu} ; \mcL_{\nu})$ into its fibre-rotation weight spaces, which we call the \textbf{fibre-rotation grading}.
The algebraic cutting method used here expresses the cut level $d \in \ZZ_{> 0}$ as a condition on these weights, producing a filtration whose pieces are identified with $H^{0}(M_{\nu}^{\leq d} ; \mcL_{\nu}^{\leq d})$ at every positive regular integral level $d$.

In \cite{EdidinGraham1998}, Edidin and Graham formulate an algebraic counterpart to Lerman's symplectic cut \cite{Lerman1995}, which they call algebraic cutting.
In their construction, one forms a GIT quotient of $M \times \mathbb{A}^{1}$ by an anti-diagonal one-dimensional torus action, with the cut level encoded by a choice of linearisation.
Here, we use the diagonal cutting convention of Subsection~\ref{subsection:subsection-2.3-arrangement-core-cuts} (i.e., opposite to theirs), so that the sub-level region $\{ \Phi \leq d \}$ is retained.
Since $M_{\nu}$ is already realised as a GIT quotient and the fibre-rotation action admits a holomorphic $\CC^{\ast}$-complexification, its algebraic cut can be formulated using graded section rings.
This description allows us to compare the cutting-invariant sections on $M_{\nu} \times \CC$ with the sections of the pre-quantum line orbibundle on the cut space.

Fix a generic integral parameter $\nu \in \mfk_{\ZZ}^{\ast}$, retaining the fixed integral lift $\lambda \in (\mft_{\ZZ}^{N})^{\ast}$ and the pre-quantum line orbibundle $\mcL_{\nu} \ra M_{\nu}$ from Section~\ref{section:section-2-conventions-geometric-input}.
We consider its section ring
\[
	R := \bigoplus\nolimits_{m \geq 0} H^{0}(M_{\nu} ; \mcL_{\nu}^{\otimes m}),
\]
whose degree $m$ records the tensor power.
The hypertoric variety $M_{\nu}^{\alg}$ is semiprojective in the sense fixed in Section~\ref{section:section-2-conventions-geometric-input}, and $\mcL_{\nu}$ is relatively ample.
By Proposition~\ref{prop:positive-ray-semi-invariants}, $R$ is isomorphic as a graded $\CC$-algebra to the semi-invariant ring $R_{\nu}$ of \eqref{eqn:GIT-hypertoric-definition}.
It follows that $R$ is a finitely generated graded $\CC$-algebra and $M_{\nu}^{\alg} \cong \Proj R$ as a projective scheme over $\Spec R_{0}$.

The fibre-rotation action extends holomorphically to $\CC^{\ast}$ and lifts to $\mcL_{\nu}$ by the complexification of the $U_{1}$-lift fixed in Subsection~\ref{subsection:subsection-3.2-pre-quantum-line-orbibundles}.
On the trivial line bundle $L_{\nu} = \mcN \times \CC \ra \mcN$, the lift acts by $\tau \cdot ( ( z , w ) , a ) = ( ( z , \tau w ) , a )$ for $\tau \in \CC^{\ast}$, leaving the fibre coordinate $a$ unchanged.
After descent, it acts on each $R_{m} = H^{0}(M_{\nu} ; \mcL_{\nu}^{\otimes m})$ by the pullback convention of Notation~\ref{notation:representation-convention}, and commutes with the residual $T^{n}$-action.

On the auxiliary affine line $\CC_{\zeta}$, the geometric action is $\tau \cdot \zeta = \tau\zeta$.
For an integer $d$, let $\varpi_{d} : \CC^{\ast} \ra \CC^{\ast}$ be the character $\varpi_{d}(\tau) = \tau^{d}$, and equip the trivial line bundle $\mcO_{\CC}(\varpi_{d}) \ra \CC_{\zeta}$ with geometric fibre weight $+d$, as in the auxiliary construction of Subsection~\ref{subsection:subsection-3.2-pre-quantum-line-orbibundles}.
Write $y(\zeta) = \zeta$ for the coordinate function and $\one_{d}$ for the nowhere-vanishing constant section of $\mcO_{\CC}(\varpi_{d})$ in its standard trivialisation; under the pullback action of Notation~\ref{notation:representation-convention}, they have weights $+1$ and $-d$, respectively.
The section ring of this line bundle is
\[
	S := \bigoplus\nolimits_{m \geq 0} H^{0}( \CC ; \mcO_{\CC}(\varpi_{d})^{\otimes m} ),
\]
and the section $\one_{d}^{\otimes m}$ trivialises the $m$-th tensor power, so $S_{m} \cong \CC[y] \cdot \one_{d}^{\otimes m}$.
Multiplication is induced by tensor product, making $S$ a polynomial algebra over $\CC[y]$ generated by $\one_{d}$ in section-ring degree one; note that the coefficient $y$ has section-ring degree zero.
In particular, $\Proj S \cong \Spec S_{0} = \Spec \CC[y] \cong \CC$.

To form the section ring of the product line orbibundle $\mcL_{\nu} \boxtimes \mcO_{\CC}(\varpi_{d}) \ra M_{\nu} \times \CC$, we pair the same tensor power in each factor by forming the \textbf{Segre product}
\[
	R~\#_{\CC}~S := \bigoplus\nolimits_{m \geq 0}R_{m} \otimes_{\CC} S_{m}.
\]
For its invariant-ring description, let an auxiliary $\CC^{\ast}$ act on $R \otimes_{\CC} S$ with weight $a - b$ on $R_{a} \otimes_{\CC} S_{b}$.
Clearly, its invariants are the terms with $a = b$, so $R~\#_{\CC}~S = ( R \otimes_{\CC} S )^{\CC^{\ast}}$ and it is graded by this common degree.
Moreover, since $R \otimes_{\CC} S$ is finitely generated and the auxiliary torus is reductive, the Segre product is finitely generated as well.
The construction is studied under the name Segre product in \cite[Chapter~4]{GotoWatanabe1978}; see also \cite[\S2.2]{BrionProcesi1990} and, as the \emph{Cartesian product}, \cite[Exercise~II.5.11]{Hartshorne1977}.
Note, however, that those references impose degree-zero hypotheses in the statements cited, whereas here $R_{0}$ need not be equal to $\CC$; only the definition above and the invariant-ring description are used below.

For the section rings $R$ and $S$ above, the K\"unneth formula \cite[Proposition~9.2.4]{Kempf1993}, applied on $\mcN^{\sstable}(\nu) \times \CC$ before taking $K_{\CC}$-invariants, shows that the common degree-$m$ piece is
\[
	( R~\#_{\CC}~S )_{m} \cong H^{0}( M_{\nu} \times \CC ; \mcL_{\nu}^{\otimes m} \boxtimes \mcO_{\CC}(\varpi_{d})^{\otimes m} ).
\]
We now take the invariants for the cutting $\CC^{\ast}$-action:
\[
	( R~\#_{\CC}~S )^{\CC^{\ast}} = \bigoplus_{m \geq 0}
	( R~\#_{\CC}~S )_{m}^{\CC^{\ast}} = \bigoplus_{m \geq 0}
	H^{0}
	\left( M_{\nu} \times \CC ; \big( \mcL_{\nu} \boxtimes \mcO_{\CC}(\varpi_{d}) \big)^{\otimes m} \right)^{\CC^{\ast}},
\]
which is the invariant section ring used to define the algebraic cut.

\begin{definition}[Algebraic Cut]
	\label{def:algebraic-cut}
	For an integer $d$, the \textbf{algebraic cut} of $M_{\nu}$ at level $d$ is the GIT quotient of $M_{\nu} \times \CC$ for the cutting $\CC^{\ast}$-action linearised on $\mcL_{\nu} \boxtimes \mcO_{\CC}(\varpi_{d})$,
	\begin{equation}
		\label{eqn:algebraic-cut}
		M_{\nu}^{\alg , \leq d} := \Proj \big( ( R~\#_{\CC}~S )^{\CC^{\ast}} \big).
	\end{equation}
\end{definition}
Since the cutting $\CC^{\ast}$ is reductive, the invariant ring $( R~\#_{\CC}~S )^{\CC^{\ast}}$ is finitely generated as well.
For a positive regular integral level $d$, the algebraic cut is the coarse space of the symplectic cut of Subsection~\ref{subsection:subsection-2.3-arrangement-core-cuts}.
Indeed, the cutting $\CC^{\ast}$-action is the complexification of the $U_{1}$-action on $M_{\nu} \times \CC$, and the character $\varpi_{d}$ specifies the moment-map normalisation $\rho - d$ used in Subsection~\ref{subsection:subsection-3.2-pre-quantum-line-orbibundles}.

The positive reals $\RR_{> 0} \subset \CC^{\ast}$ act by the gradient flow of $\rho$, so $\rho$ is strictly increasing along every non-trivial $\RR_{> 0}$-orbit.
Since $\rho$ decreases as $\tau \ra 0$ in $\RR_{>0}$, each half-orbit $\{ \tau \cdot (p, \zeta) ~|~ 0 < \tau \leq 1 \}$ lies in $\rho^{-1}( [ 0 , \rho(p, \zeta) ] )$, which is compact by Lemma~\ref{lemma:proper-moment-map}.
Since the orbit maps are algebraic, the limit of $\tau \cdot ( p , \zeta )$ as $\tau \to 0$ exists for every point and is a fixed point, at which $\rho$ coincides with a critical value of $\Phi$; as $\tau \to \infty$, the limit exists if and only if $\zeta = 0$ and $p \in \mcC$ \cite[Corollary~3.6]{HaradaProudfoot2004}, and otherwise $\rho \to \infty$ along the orbit.

By Proposition~\ref{prop:positive-ray-semi-invariants}, $( R~\#_{\CC}~S )^{\CC^{\ast}}$ is the ring of semi-invariants on $\mcN \times \CC$ for the torus $K_{\CC} \times \CC^{\ast}$ and the characters $m(\nu , d)$, so $M_{\nu}^{\alg , \leq d}$ is the corresponding GIT quotient of the affine variety $\mcN \times \CC$.
By the form of the Kempf--Ness theorem for characters due to King \cite[Proposition~3.21 and Corollary~3.22]{Nakajima1999}, applied to the closed invariant subvariety $\mcN \times \CC$ of $T^{\ast}\CC^{N} \times \CC$ as in \cite[Theorem~3.24]{Nakajima1999}, a point of $\mcN \times \CC$ is semistable if and only if the closure of its $K_{\CC} \times \CC^{\ast}$-orbit meets $( \phi_{\RR}^{-1}(\nu) \times \CC ) \cap \rho^{-1}(d)$.
Reducing first by $K$ through \eqref{eqn:hypertoric-kempf-ness}, this says that the closure of the $\CC^{\ast}$-orbit of $( p , \zeta )$ in $M_{\nu} \times \CC$ meets $\rho^{-1}(d)$.
This closure is the orbit together with its limits, at which $\rho \neq d$, so $( p , \zeta )$ is semistable if and only if its orbit meets $\rho^{-1}(d)$, in which case it meets $\rho^{-1}(d)$ in a single $U_{1}$-orbit.
The limits of a semistable orbit are fixed points with $\rho \neq d$, hence unstable, so every semistable orbit is closed in the semistable locus, and the GIT quotient is geometric.
The inclusion of $\rho^{-1}(d)$ into the semistable locus induces an isomorphism of complex orbifolds $\rho^{-1}(d) / U_{1} \cong [ (M_{\nu} \times \CC)^{\sstable}(\varpi_{d}) / \CC^{\ast} ]$, with the complex structure on $\rho^{-1}(d)/U_{1}$ induced from this quotient, as in \cite[Theorems~3.2 and~3.4]{Hausel1998}.
Under this identification, the line orbibundle obtained by descending $\mcL_{\nu} \boxtimes \mcO_{\CC}(\varpi_{d})$ to this quotient is the cut pre-quantum line orbibundle $\mcL_{\nu}^{\leq d}$ constructed in Subsection~\ref{subsection:subsection-3.2-pre-quantum-line-orbibundles}.
We henceforth write $M_{\nu}^{\leq d}$ for this common complex orbifold, whose coarse space is $M_{\nu}^{\alg, \leq d}$, at a positive regular integral level $d$.
At an integral critical level, \eqref{eqn:algebraic-cut} still defines the corresponding algebraic GIT quotient, although it need not be an orbifold.

\begin{remark}
	\label{remark:algebraic-cut-generality}
	The construction above applies more generally to a polarised semiprojective variety equipped with an algebraic $\CC^{\ast}$-action and a compatible linearisation, provided that the associated section ring is finitely generated.
	Our section-space and combinatorial results that follow, however, use the additional structure of hypertoric varieties.
	In particular, Subsection~\ref{subsection:subsection-5.3-cut-chambers-vanishing} describes the unstable locus of every cut in terms of the fibre-rotation cells of $M_{\nu}$.
	The hypertoric coordinate ring then allows us to compute the residual $T^{n}$-weight multiplicities explicitly in Theorem~\ref{theorem:combinatorial-multiplicity}.
\end{remark}

The cutting linearisation imposes an explicit upper bound on fibre-rotation weights.
For $s_{i} \in H^{0}(M_{\nu} ; \mcL_{\nu})_{i}$ and $j \in \ZZ_{\geq 0}$, the product $s_{i}y^{j} \otimes \one_{d}$ is a section of $\mcL_{\nu} \boxtimes \mcO_{\CC}(\varpi_{d})$ on $M_{\nu} \times \CC$.
Since $y$ and $\one_{d}$ have pullback weights $+1$ and $-d$, its cutting weight is $i + j - d$, and hence $s_{i} y^{j} \otimes \one_{d}$ has cutting weight zero if and only if $i + j - d = 0$.
A non-negative exponent $j$ exists only when $i \leq d$, in which case $j = d - i$.
For every integer $d \geq 0$, we therefore define the \textbf{cut-filtration piece}
\[
	F_{d} H^{0}(M_{\nu} ; \mcL_{\nu}) := \bigoplus\nolimits_{i \leq d} H^{0}(M_{\nu} ; \mcL_{\nu})_{i}.
\]
These spaces form an increasing filtration that is exhaustive because the $\CC^{\ast}$-action is rational, so every section is a finite sum of weight vectors.
Theorem~\ref{theorem:combinatorial-multiplicity} will show that the grading is non-negative and that each $F_{d}$ is finite-dimensional.
For $d \geq 1$, the successive quotients recover the individual fibre-rotation weight spaces:
\[
	\frac{F_{d}H^{0}(M_{\nu} ; \mcL_{\nu})}{F_{d - 1}H^{0}(M_{\nu} ; \mcL_{\nu})} \cong H^{0}(M_{\nu} ; \mcL_{\nu})_{d}.
\]

For a positive regular integral level $d$, write $U_{d} := ( M_{\nu} \times \CC )^{\sstable}(\varpi_{d})$ and $\mcG_{d} := \mcL_{\nu} \boxtimes \mcO_{\CC}(\varpi_{d}) \lra M_{\nu} \times \CC$.
Here, $U_{d}$ is the semistable locus for the cutting $\CC^{\ast}$-action; the variation of this locus with the cut level is discussed further in Subsection~\ref{subsection:subsection-5.3-cut-chambers-vanishing}.
As shown after Definition~\ref{def:algebraic-cut}, a point is semistable if and only if its $\CC^{\ast}$-orbit meets $\rho^{-1}(d)$.

The cut-filtration pieces are the invariant sections on $M_{\nu} \times \CC$.

\begin{lemma}[Invariant Sections]
	\label{lemma:filtration-invariant-sections}
	For every integer $d \geq 0$, there is a natural $T^{n}$-equivariant isomorphism
	\begin{equation}
		\label{eqn:filtration-invariant-sections}
		F_{d}H^{0}(M_{\nu} ; \mcL_{\nu}) \overset{\sim}{\longrightarrow} H^{0}\left( M_{\nu} \times \CC ; \mcL_{\nu} \boxtimes \mcO_{\CC}(\varpi_{d}) \right)^{\CC^{\ast}},
	\end{equation}
	given on weight vectors by $s_{i} \mapsto s_{i}y^{d - i} \otimes \one_{d}$.
\end{lemma}

\begin{proof}
	Since the auxiliary factor $\CC$ is affine, the K\"unneth formula \cite[Proposition~9.2.4]{Kempf1993}, applied on $\mcN^{\sstable}(\nu) \times \CC$ before taking $K_{\CC}$-invariants, gives
	\[
		H^{0}( M_{\nu} \times \CC ; \mcL_{\nu} \boxtimes \mcO_{\CC}(\varpi_{d}) )
		\cong H^{0}(M_{\nu} ; \mcL_{\nu}) \otimes_{\CC} ( \CC[y] \cdot \one_{d} ).
	\]
	By the weight calculation before the definition of $F_{d}$, for $s_{i} \in H^{0}(M_{\nu} ; \mcL_{\nu})_{i}$ the section $s_{i}y^{j} \otimes \one_{d}$ is invariant if and only if $j = d - i \geq 0$, and every invariant section is uniquely a finite sum of terms of this form, so
	\[
		H^{0}( M_{\nu} \times \CC ; \mcL_{\nu} \boxtimes \mcO_{\CC}(\varpi_{d}) )^{\CC^{\ast}} \cong
		\bigoplus\nolimits_{i \leq d} H^{0}(M_{\nu} ; \mcL_{\nu})_{i} y^{d - i} \otimes \one_{d}.
	\]
	Thus, the map defined on weight vectors by $s_{i} \mapsto s_{i}y^{d - i} \otimes \one_{d}$ gives the isomorphism \eqref{eqn:filtration-invariant-sections}.
\end{proof}

\subsection{Cut Chambers and Vanishing of Higher Cohomology}
\label{subsection:subsection-5.3-cut-chambers-vanishing}

Theorem~\ref{thm:equivariant-KRR-index-theorem} computes the equivariant holomorphic Euler characteristic defined in \eqref{eqn:holomorphic-Euler-characteristic}.
We seek the identity
\[
	\chi_{T^{n}}(M_{\nu}^{\leq d} ; \mcL_{\nu}^{\leq d})
	= [ H^{0}(M_{\nu}^{\leq d} ; \mcL_{\nu}^{\leq d}) ] \in R(T^{n}).
\]
For this, it suffices to prove $H^{q}(M_{\nu}^{\leq d} ; \mcL_{\nu}^{\leq d}) = 0$ for every $q > 0$.
Theorem~\ref{theorem:cut-cohomology-every-level} below proves this at every positive regular integral level, and identifies $H^{0}(M_{\nu}^{\leq d} ; \mcL_{\nu}^{\leq d})$ with $F_{d}H^{0}(M_{\nu} ; \mcL_{\nu})$.

Throughout this subsection, let $d$ be a positive regular integral level.
Write
\[
	Y_{d} := M_{\nu}^{\leq d}, \qquad
	D_{d} := Z_{\nu}^{d}, \qquad
	\mcL_{d} := \mcL_{\nu}^{\leq d}, \qquad
	\w_{d} := \w_{\RR}^{\leq d}
\]
Recall that $Y_{d}$ is a compact K\"ahler orbifold of complex dimension $\dim_{\CC}M_{\nu} = 2n$, with K\"ahler form $\w_{d}$, cut divisor $D_{d} \subset Y_{d}$, and pre-quantum line orbibundle $\mcL_{d} \ra Y_{d}$.
Cohomology on $M_{\nu}$, on $M_{\nu} \times \CC$, and on the quotient stacks below is algebraic coherent cohomology.
Points of $M_{\nu}$ and of $M_{\nu} \times \CC$ mean points of their coarse spaces, and limits and closures are taken there.
Weights of sections are taken for the pullback action of Notation~\ref{notation:representation-convention}, and $(\, \cdot\, )_{w}$ denotes the weight-$w$ part of a rational $\CC^{\ast}$-module.
For a cocharacter $\theta$ of a torus fixing a point $x$, the \textbf{$\theta$-weight} of an equivariant line bundle at $x$ is the geometric weight of $\theta$ on its fibre at $x$; this is the convention of \cite{HalpernLeistner2015}, in which the conormal bundle of a blade, restricted to its fixed component, has positive weights \cite[Definition~2.2, (S3)]{HalpernLeistner2015}.

We deduce Theorem~\ref{theorem:cut-cohomology-every-level} below from Teleman's quantisation theorem \cite{Teleman2000}, in the form proved by Halpern-Leistner \cite[Theorem~3.29]{HalpernLeistner2015}: removing an unstable stratum does not change the invariant cohomology of a line bundle whose weight along the stratum is negative.
Hypertoric cuts enter only through the verification of the theorem's hypotheses: a global quotient presentation of the cut, the list of unstable strata (Lemma~\ref{lemma:cut-unstable-locus}), their identification with blades (Proposition~\ref{prop:cells-are-blades}), and one weight computation for each stratum (Lemma~\ref{lemma:stratum-weights}).

We first present $M_{\nu}$, the cut $Y_{d}$, and the auxiliary quotient $[\widetilde{X} / \widetilde{K}]$ as global quotient stacks.
Write $\mcN^{\sstable} := \mcN^{\sstable}(\nu)$ and $Q : \mcN^{\sstable} \ra M_{\nu}^{\alg}$ for the quotient map.
Since $( \nu , 0 )$ is a regular value of $\phi_{\HK}$, the infinitesimal $K$-action is injective along $\phi_{\HK}^{-1}(\nu , 0)$, and the stabiliser in $K_{\CC}$ of a point of $T^{\ast}\CC^{N}$ is the intersection of the kernels of the characters of its non-zero coordinates, so the stabiliser of each point of $\phi_{\HK}^{-1}(\nu, 0)$ is finite.
Every $K_{\CC}$-orbit in $\mcN^{\sstable}$ meets this level set by \eqref{eqn:hypertoric-kempf-ness}, so every point of $\mcN^{\sstable}$ has finite stabiliser, $\phi_{\CC}$ is submersive along $\mcN^{\sstable}$, and $\mcN^{\sstable}$ is smooth.
Since stability and semistability coincide for $\nu$, as noted in Subsection~\ref{subsection:subsection-2.2-algebraic-quotients-prequantum}, every point of $\mcN^{\sstable}$ is properly stable \cite[Definition~1.8]{MFK1994}, and $Q$ is a uniform geometric quotient \cite[Theorem~1.10(iii) and Definition~0.7]{MFK1994}.
The fibre-rotation action preserves $\mcN^{\sstable}$ and commutes with $K_{\CC}$.
Hence, the torus $\widetilde{K} := K_{\CC} \times \CC^{\ast}$, whose second factor is the cutting $\CC^{\ast}$, acts on the smooth quasi-projective variety $\widetilde{X} := \mcN^{\sstable} \times \CC$ by $( k , \tau ) \cdot ( x , \zeta ) := ( k \cdot ( \tau \cdot x ) , \tau\zeta )$, where $\tau \cdot x$ denotes the fibre-rotation action.
This action is linearisable, as $\widetilde{X}$ is an invariant locally closed subvariety of the representation $T^{\ast}\CC^{N} \times \CC$.
The base change $\widetilde{Q} := Q \times \Identity : \widetilde{X} \ra M_{\nu}^{\alg} \times \CC$ is a geometric quotient by $K_{\CC}$.
Put $\widetilde{X}^{\sstable} := \widetilde{Q}^{-1}(U_{d})$; we regard $M_{\nu} = [ \mcN^{\sstable} / K_{\CC} ]$ and
\[
	[ \widetilde{X} / \widetilde{K} ] = [ ( M_{\nu} \times \CC ) / \CC^{\ast} ],
	\qquad
	Y_{d} = [ \widetilde{X}^{\sstable} / \widetilde{K} ]
\]
as quotient stacks.
The quotient stacks displayed above are smooth.
The stabilisers of $M_{\nu}$ and $Y_{d}$ are finite, and only finitely many stabiliser groups occur, so these two stacks have finite inertia and are tame in characteristic zero \cite[Proposition~2.7 and Theorem~3.2]{AbramovichOlssonVistoli2008}.
Since these stacks and their coarse spaces have affine diagonals, their coarse-space maps are good moduli spaces without higher direct images of quasi-coherent sheaves \cite[Remarks~3.5 and~7.3]{Alper2013}, and their coherent cohomology agrees with the coherent orbifold cohomology used in this paper, as in the proof of Proposition~\ref{prop:positive-ray-semi-invariants}.
The line bundle $L_{\nu} \boxtimes \mcO_{\CC}(\varpi_{d})$ on $\widetilde{X}$, linearised by \eqref{eqn:positive-character-linearisation} for $K_{\CC}$ and by the lifts of Subsection~\ref{subsection:subsection-5.2-algebraic-cutting} for the cutting factor, gives $\mcG_{d}$ on $[ \widetilde{X} / \widetilde{K} ]$ and $\mcL_{d}$ on $Y_{d}$; we also write $\mcG_{d}$ for the bundle on $\widetilde{X}$ and for its descent to $M_{\nu} \times \CC$.
The larger quotient $[ \widetilde{X}/\widetilde{K} ]$ is a smooth Artin stack with positive-dimensional stabilisers along the loci $\widetilde{Z}_{\beta}$ of Lemma~\ref{lemma:stabilisers-fixed-components}; for this stack, we use only $H^{q}( [ \widetilde{X} / \widetilde{K} ] ; \mcG_{d} ) = H^{q}( \widetilde{X} ; \mcG_{d} )^{\widetilde{K}}$.

Let $\CC^{\ast}$ act on $\widetilde{X}$ by $( x , \zeta ) \mapsto ( \tau \cdot x , \zeta )$, lifted to $\mcG_{d}$ by the lift of Subsection~\ref{subsection:subsection-5.2-algebraic-cutting} on the factor $L_{\nu}$ and trivially on $\mcO_{\CC}(\varpi_{d})$, so that $y$ and $\one_{d}$ are fixed.
This action commutes with $\widetilde{K}$ and with $T^{n}$, and its lift preserves the invariant sections of each power of $\mcG_{d}$, so it preserves $U_{d}$ and descends to an action on $Y_{d}$ with a lift to $\mcL_{d}$, which we call the \textbf{residual fibre-rotation action}.

For every positive regular integral level $d$, let $\Psi_{d}$ be the composite of \eqref{eqn:filtration-invariant-sections} with restriction to $U_{d}$ and descent, so that $\Psi_{d}(s)$ is the descent of $s y^{d - w} \otimes \one_{d}$ for $s$ of fibre-rotation weight $w \leq d$.
For the residual fibre-rotation action, this section has weight $w$, so $\Psi_{d}$ is $T^{n} \times \CC^{\ast}$-equivariant.

The following theorem shows that, in every degree, cutting at a positive regular integral level $d$ keeps the fibre-rotation weights $w \leq d$ of the cohomology of $\mcL_{\nu}$ and removes the weights $w > d$.

\begin{theorem}
	\label{theorem:cut-cohomology-every-level}
	Let $d$ be a positive regular integral level.
	Then, for every $q \geq 0$, restriction and descent induce natural $T^{n}$-equivariant isomorphisms
	\[
		H^{q}(Y_{d} ; \mcL_{d})_{w} \cong H^{q}(M_{\nu} ; \mcL_{\nu})_{w} \quad ( w \leq d ),
		\qquad
		H^{q}(Y_{d} ; \mcL_{d})_{w} = 0 \quad ( w > d ),
	\]
	where the weights are taken for the residual fibre-rotation action on $Y_{d}$ and for the fibre-rotation action on $M_{\nu}$.
	In particular:
	\begin{enumerate}[(i)]
		\item the map $\Psi_{d} : F_{d}H^{0}(M_{\nu} ; \mcL_{\nu}) \ra H^{0}(Y_{d} ; \mcL_{d})$ is a $T^{n} \times \CC^{\ast}$-equivariant isomorphism;
		\item $H^{q}(Y_{d} ; \mcL_{d}) = 0$ for every $q > 0$;
		\item $\chi_{T^{n}}(Y_{d} ; \mcL_{d}) = [ F_{d}H^{0}(M_{\nu} ; \mcL_{\nu}) ]$ in $R(T^{n})$;
		\item if $d - 1$ is also a positive regular integral level, then $j_{d} := \Psi_{d} \circ ( F_{d - 1} \hookrightarrow F_{d} ) \circ \Psi_{d - 1}^{-1}$ is a $T^{n}$-equivariant injection $H^{0}(Y_{d - 1} ; \mcL_{d - 1}) \hookrightarrow H^{0}(Y_{d} ; \mcL_{d})$, and
		      \[
			      H^{0}(M_{\nu} ; \mcL_{\nu})_{d}
			      \cong \frac{H^{0}(Y_{d} ; \mcL_{d})}{j_{d} H^{0}(Y_{d - 1} ; \mcL_{d - 1})}
		      \]
		      as residual $T^{n}$-representations.
	\end{enumerate}
\end{theorem}

We only need the following case of Halpern-Leistner's theorem \cite[Theorem~3.29]{HalpernLeistner2015}.
Teleman \cite[Proposition~2.6(a), Proposition~2.11(a), and Remark~2.12(i)]{Teleman2000} proves the stratum-removal principle stated before Theorem~\ref{theorem:cut-cohomology-every-level}.

\begin{theorem}[{Teleman, Halpern-Leistner}]
	\label{theorem:quantisation-theorem}
	Let a torus $G$ act linearisably on a smooth quasi-projective variety $X$.
	Let $\theta$ be a cocharacter of $G$ and $Z$ a connected component of the fixed locus $X^{\theta}$, and suppose that the \textbf{blade}
	\[
		S := \{ x \in X ~|~ \lim\nolimits_{t \ra 0} \theta(t) \cdot x \text{ exists in } X \text{ and lies in } Z \}
	\]
	is closed in $X$.
	Let $E$ be a $G$-equivariant line bundle on $X$ whose $\theta$-weight along $Z$ is negative.
	Then, restriction is an isomorphism
	\[
		R\Gamma( [ X / G ] ; E ) \overset{\sim}{\lra} R\Gamma( [ ( X \setminus S ) / G ] ; E ),
	\]
	and $H^{q}( [ X / G ] ; E ) = H^{q}(X ; E)^{G}$, and likewise for $X \setminus S$.
\end{theorem}

\begin{proof}
	Since $G$ is a torus, condition (S2) of \cite[Definition~2.2]{HalpernLeistner2015} is automatic, and (S1) and (S3) hold for every blade \cite[Remark~2.1]{HalpernLeistner2015}, so $S$ is a closed KN stratum; Properties (A) and (L+) hold because $X$ is smooth \cite[Lemma~2.7]{HalpernLeistner2015}.
	Write $j : S \ra X$ and $\sigma : Z \ra S$ for the inclusions.
	The structure sheaf lies in $D^{-}( [ X / G ] )_{\geq 0}$, since $L\sigma^{\ast}Lj^{\ast}\mcO = \mcO_{Z}$ has $\theta$-weight $0$ \cite[Definition~3.24 and the definition preceding Proposition~3.9]{HalpernLeistner2015}.
	By Property~(A) and \cite[Proposition~3.31]{HalpernLeistner2015}, $E$ lies in $D^{+}( [ X / G ] )_{< 0}$ if and only if $\sigma^{\ast}j^{!}E \in D^{+}(Z)_{< 0}$.
	The proof of \cite[Lemma~2.7]{HalpernLeistner2015} shows that $j$ is a regular embedding; if $c$ is its codimension, then  $\sigma^{\ast}j^{!}E \cong E|_{Z} \otimes \det(N_{S/X})|_{Z}[ -c ]$, as in the proof of \cite[Lemma~2.9]{HalpernLeistner2015}.
	If $\eta \geq 0$ is the $\theta$-weight of $\det(N_{S/X}^{\vee})|_{Z}$ \cite[(4)]{HalpernLeistner2015}, then this complex has the single $\theta$-weight $\wt_{\theta}( E|_{Z} ) - \eta < 0$.
	Halpern-Leistner's theorem \cite[Theorem~3.29]{HalpernLeistner2015}, applied to the pair $( \mcO , E )$ with $w = v = 0$, therefore shows that the restriction map $R \mathrm{Hom}(\mcO , E) = R\Gamma(E)$ is an isomorphism.
	Since a torus is reductive, cohomology on $[ X / G ]$ is the $G$-invariant part of cohomology on $X$ \cite[(A.3) and~(A.5)]{Teleman2000}.
\end{proof}

By the comparison in Theorem~\ref{theorem:cut-cohomology-every-level}, the higher cohomology groups on the cut vanish once those on $M_{\nu}$ do; the following lemma shows that they do, using the Koszul resolution in the proof of Proposition~\ref{prop:positive-ray-semi-invariants}.

\begin{lemma}
	\label{lemma:hypertoric-higher-vanishing}
	For every $m \geq 0$ and $q > 0$, we have $H^{q}(M_{\nu} ; \mcL_{\nu}^{\otimes m}) = 0$.
\end{lemma}

\begin{proof}
	In the proof of Proposition~\ref{prop:positive-ray-semi-invariants}, the pushforward of $\mcL_{\nu}^{\otimes m}$ to the Lawrence quotient $W_{\nu}$ is resolved by a Koszul complex with terms $\mcA_{\nu}^{\otimes m} \otimes \bigwedge^{k} \mfk_{\CC}$ for $0 \leq k \leq r$, and these terms have no higher cohomology by \cite[Theorem~5.3]{Fujino2007}.
	A finite resolution by acyclic sheaves computes cohomology, so $H^{q}(M_{\nu} ; \mcL_{\nu}^{\otimes m})$ is the $q$-th cohomology of the complex of global sections of the resolution.
	That complex sits in degrees $-r , \ldots , 0$, so its cohomology vanishes in positive degrees.
\end{proof}

To locate the unstable points of $M_{\nu} \times \CC$, we follow each point of $M_{\nu}$ along its fibre-rotation orbit in both directions and record the limits of $\Phi$ at the two ends.
The following lemma records their basic properties, and identifies the complement of the cut divisor in $Y_{e}$ with an open subset of $M_{\nu}$.

\begin{lemma}
	\label{lemma:limit-levels}
	For $p \in M_{\nu}$ and real $\tau$, put
	\[
		\Phi_{0}(p) := \inf_{0 < \tau \leq 1} \Phi( \tau \cdot p ), \qquad
		\Phi_{\infty}(p) := \sup_{\tau \geq 1} \Phi( \tau \cdot p ) \in [ 0 , \infty ].
	\]
	\begin{enumerate}[(i)]
		\item the function $\tau \mapsto \Phi( \tau \cdot p )$ on $( 0 , \infty )$ is non-decreasing, and it is strictly increasing unless $p$ is fixed by the fibre-rotation action;
		\item the limit $\lim_{\tau \ra 0} \tau \cdot p$ exists, is a fixed point, and has $\Phi$-value $\Phi_{0}(p)$; the limit $\lim_{\tau \ra \infty} \tau \cdot p$ exists in $M_{\nu}$ if and only if $p \in \mcC$, if and only if $\Phi_{\infty}(p) < \infty$, and it then has $\Phi$-value $\Phi_{\infty}(p)$;
		\item we have $\Phi_{0} \leq \Phi \leq \Phi_{\infty}$, $\Phi_{0}$ and $\Phi_{\infty}$ are invariant under $T^{n} \times \CC^{\ast}$, their finite values lie in the finite set $\CritVal(\Phi)$, $\Phi_{0}$ is upper semicontinuous, and $\Phi_{\infty}$ is lower semicontinuous;
		\item for every positive regular integral level $e$, the map $\iota_{e} : \{ \Phi_{0} < e \} \ra Y_{e}$, $p \mapsto [ ( p , 1 ) ]$, is a $T^{n}$-equivariant open immersion onto $Y_{e} \setminus D_{e}$ that intertwines the fibre-rotation action on $M_{\nu}$ with the residual fibre-rotation action on $Y_{e}$, and the residual fibre-rotation action fixes $D_{e}$ pointwise.
		      If $e \in \Cout$, then the domain of $\iota_{e}$ is $M_{\nu}$.
	\end{enumerate}
\end{lemma}

\begin{proof}
	The facts recorded after Definition~\ref{def:algebraic-cut}, applied to the points $( p , 0 )$, for which $\rho( \tau \cdot ( p , 0 ) ) = \Phi( \tau \cdot p )$, give~(i) and~(ii): along $M_{\nu} \times \{ 0 \}$, the positive reals act by the gradient flow of $\Phi$, whose stationary points are the fixed points of the fibre-rotation action, and the orbit of a point that is not fixed contains no fixed point.
	For~(iii), the inequalities are immediate from the definitions, and by~(ii), the finite values of $\Phi_{0}$ and $\Phi_{\infty}$ are values of $\Phi$ at fixed points, which are critical; $\CritVal(\Phi)$ is finite by \cite[Proposition~4.7]{BrownI}.
	Invariance follows from $\lim_{\tau} \tau \cdot ( s \cdot p ) = \lim_{\tau} ( \tau s ) \cdot p$ and from $T^{n}$ commuting with the fibre-rotation action and preserving $\Phi$.

	Finally, $\Phi_{0}$ is an infimum and $\Phi_{\infty}$ a supremum of continuous functions, so they are upper and lower semicontinuous respectively.
	For~(iv), let $\zeta \neq 0$.
	The function $\tau \mapsto \rho( \tau \cdot ( p , \zeta ) ) = \Phi( \tau \cdot p ) + \tfrac{1}{2}\tau^{2}|\zeta|^{2}$ on $( 0 , \infty )$ is continuous and strictly increasing, with limits $\Phi_{0}(p)$ and $\infty$, so $( p , \zeta ) \in U_{e}$ if and only if $\Phi_{0}(p) < e$.
	The invariant map $( p , \zeta ) \mapsto \zeta^{-1} \cdot p$ makes $M_{\nu} \times \CC^{\ast}$ a trivial $\CC^{\ast}$-torsor over $M_{\nu}$, with section $p \mapsto ( p , 1 )$.
	Since $\Phi_{0}( \zeta^{-1} \cdot p ) = \Phi_{0}(p)$, this map identifies $Y_{e} \setminus D_{e}$ with $\{ \Phi_{0} < e \}$, with inverse $\iota_{e}$.
	This identification is equivariant, since $T^{n}$ and the residual fibre-rotation action act on the first factor and preserve $\{ \zeta = 1 \}$.
	The residual fibre-rotation action of $\tau$ agrees with the cutting action of $\tau$ on the points $( p , 0 ) \in U_{e}$, so it fixes $D_{e}$ pointwise.
	If $e \in \Cout$, then $\Phi_{0} \leq c_{\max} < e$ by~(iii), so the domain of $\iota_{e}$ is $M_{\nu}$.
\end{proof}

The fixed components of the fibre-rotation action decompose $M_{\nu}$ into attracting cells, and the core $\mcC$ into repelling cells.
Since $M_{\nu}$ is not proper, we obtain these cells from the Bia\l{}ynicki-Birula decomposition of the compact cut $Y_{\delta}$ at a level $\delta$ in the outer cut chamber, which contains $M_{\nu}$ as the complement of its cut divisor by Lemma~\ref{lemma:limit-levels}(iv).

\begin{proposition}
	\label{prop:BB-cells}
	Let $\{ \mcZ_{\beta} \}_{\beta}$ be the connected components of the fixed locus of the fibre-rotation action on $M_{\nu}$, and fix an integral level $\delta \in \Cout$.
	\begin{enumerate}[(i)]
		\item there are finitely many $\mcZ_{\beta}$, each compact, and $\Phi$ takes a constant value $\Phi(\mcZ_{\beta}) \in \CritVal(\Phi)$ on each;
		\item the attracting sets $A_{\beta}^{+} := \{ p ~|~ \lim_{\tau \ra 0} \tau \cdot p \in \mcZ_{\beta} \}$ partition $M_{\nu}$, and the indices can be ordered as a sequence $\beta_{1} , \beta_{2} , \ldots$ so that the closure of $A_{\beta_{k}}^{+}$ lies in $A_{\beta_{1}}^{+} \cup \cdots \cup A_{\beta_{k}}^{+}$;
		\item the repelling sets of the fixed components of $Y_{\delta}$, defined by $\tau \ra \infty$ for the residual fibre-rotation action, partition $Y_{\delta}$, and they can be ordered as a sequence $R_{1} , R_{2} , \ldots$ so that the closure of $R_{k}$ lies in $R_{1} \cup \cdots \cup R_{k}$.
		      For each $\beta$, the repelling set of $\iota_{\delta}(\mcZ_{\beta})$ is $\iota_{\delta}(A_{\beta}^{-})$, where $A_{\beta}^{-} := \{ p ~|~ \lim_{\tau \ra \infty} \tau \cdot p \text{ exists in } M_{\nu} \text{ and lies in } \mcZ_{\beta} \}$;
		\item $A_{\beta}^{+} \cap A_{\beta}^{-} = \mcZ_{\beta}$;
		\item if $\Sigma \subseteq M_{\nu}$ is a closed union of sets $A_{\beta}^{+}$, or a compact union of sets $A_{\beta}^{-}$, then it has a finite filtration by subsets that are closed in $M_{\nu}$; the non-empty successive differences of this filtration are single sets $A_{\beta}^{\pm}$, each closed in the complement of the previous step.
	\end{enumerate}
\end{proposition}

\begin{proof}
	Since $\delta$ is regular, the cutting action has finite stabilisers on $U_{\delta}$, so the stack $Y_{\delta}$ is smooth and tame.
	It is also proper, and its coarse space is the projective variety $\Proj$ in \eqref{eqn:algebraic-cut}.
	By \cite[Theorem~3.5]{Skowera2013}, the residual fibre-rotation action decomposes $Y_{\delta}$ into disjoint, locally closed, invariant substacks that are equivariant affine fibrations over the fixed components, in the form of \cite[Proposition~5]{Oprea2006}.
	The equivariant affine \'{e}tale atlas that Oprea requires is supplied by \cite[Proposition~3.2]{Skowera2013}, after the reparametrisation of \cite[Remark~3.4]{Skowera2013}, since $Y_{\delta}$ is irreducible, being smooth and connected, and has trivial generic stabiliser.
	The proof of \cite[Proposition~5]{Oprea2006} shows that these substacks descend from the plus cells of \cite{BialynickiBirula1973} on the atlas, so they are the attracting sets of the fixed components as $\tau \ra 0$, and by \cite[Remark~3.6]{Skowera2013} they form a filtration in the sense of \cite[\S2.3, (a)]{Oprea2006}.
	The same results, applied to the inverse action, give the first assertion of~(iii).

	By Lemma~\ref{lemma:limit-levels}(iv), $\iota_{\delta}$ identifies $M_{\nu}$ equivariantly with $Y_{\delta} \setminus D_{\delta}$, and $D_{\delta}$ is fixed pointwise.
	Under the identification $Y_{\delta} = \rho^{-1}(\delta) / U_{1}$ made after Definition~\ref{def:algebraic-cut}, the function $[(q, \zeta)] \mapsto \Phi(q)$ is continuous on $Y_{\delta}$ and equal to $\delta$ on $D_{\delta}$.
	Each $p \in \mcZ_{\beta}$ is fixed, so $\iota_{\delta}(p)$ is represented by $(p , (2\delta - 2\Phi(p))^{1/2}) \in \rho^{-1}(\delta)$, where the function takes the value $\Phi(p)$.
	As $\mcZ_{\beta}$ is connected and consists of critical points, $\Phi$ is constant on it, so the function is constant on $\iota_{\delta}(\mcZ_{\beta})$, equal to $\Phi(\mcZ_{\beta}) \leq c_{\max} < \delta$.
	Each of the finitely many connected components of the fixed locus of the proper stack $Y_{\delta}$ therefore lies in $D_{\delta}$ or is some $\iota_{\delta}(\mcZ_{\beta})$.
	Each $\mcZ_{\beta}$ is closed in $\Phi^{-1}( \Phi(\mcZ_{\beta}) )$, which is compact because $\rho$ is proper by Lemma~\ref{lemma:proper-moment-map}; this gives~(i).
	For~(ii), the limit of $\tau \cdot \iota_{\delta}(p)$ as $\tau \ra 0$ is $\iota_{\delta}( \lim_{\tau \ra 0} \tau \cdot p ) \not\in D_{\delta}$ by Lemma~\ref{lemma:limit-levels}(ii), so the attracting sets of $Y_{\delta}$ are the components of $D_{\delta}$ and the sets $\iota_{\delta}(A_{\beta}^{+})$.
	The filtration therefore restricts to $M_{\nu}$, since the closure in $M_{\nu}$ of a subset of $M_{\nu}$ is the preimage under $\iota_{\delta}$ of its closure in $Y_{\delta}$.
	For the last assertion of~(iii), points of $D_{\delta}$ are fixed, so the repelling set of $\iota_{\delta}(\mcZ_{\beta})$ lies in $\iota_{\delta}(M_{\nu})$, where it is $\iota_{\delta}(A_{\beta}^{-})$ by equivariance.
	For~(iv), if $p \in A_{\beta}^{+} \cap A_{\beta}^{-}$, then $\Phi( \tau \cdot p )$ tends to $\Phi(\mcZ_{\beta})$ both as $\tau \ra 0$ and as $\tau \ra \infty$, so $p$ is fixed by Lemma~\ref{lemma:limit-levels}(i).

	For~(v), in the first case put $\Sigma_{k} := \Sigma \cap ( A_{\beta_{1}}^{+} \cup \cdots \cup A_{\beta_{k}}^{+} )$.
	By~(ii), the union of the first $k$ attracting sets contains the closure of each of its members, so $\Sigma_{k}$ is closed.
	The difference $\Sigma_{k} \setminus \Sigma_{k - 1} = \Sigma \cap A_{\beta_{k}}^{+}$ is either empty or $A_{\beta_{k}}^{+}$, whose closure lies in $\Sigma_{k}$ and so meets $M_{\nu} \setminus \Sigma_{k - 1}$ only in $A_{\beta_{k}}^{+}$.
	In the second case, $\iota_{\delta}(\Sigma)$ is compact, hence closed in $Y_{\delta}$, and it is a union of repelling sets by~(iii); the same argument in $Y_{\delta}$, with the ordering $R_{1} , R_{2} , \ldots$, followed by intersection with the open set $\iota_{\delta}(M_{\nu})$, gives the claim.
\end{proof}

Since a point of $M_{\nu} \times \CC$ is semistable if and only if its cutting orbit meets $\rho^{-1}(d)$, it is unstable if and only if $\rho$ stays above $d$, or stays below $d$, along that orbit.
The following lemma describes both kinds of unstable point in terms of the limits of Lemma~\ref{lemma:limit-levels} and the cells of Proposition~\ref{prop:BB-cells}.

\begin{lemma}
	\label{lemma:cut-unstable-locus}
	Put $\Sigma_{d}^{+} := \{ \Phi_{0} > d \}$ and $\Sigma_{d}^{-} := \{ \Phi_{\infty} < d \}$.
	Then
	\[
		( M_{\nu} \times \CC ) \setminus U_{d} = ( \Sigma_{d}^{+} \times \CC ) \sqcup ( \Sigma_{d}^{-} \times \{ 0 \} ), \qquad \Sigma_{d}^{+} = \bigcup_{\Phi(\mcZ_{\beta}) > d} A_{\beta}^{+}, \qquad \Sigma_{d}^{-} = \bigcup_{\Phi(\mcZ_{\beta}) < d} A_{\beta}^{-}.
	\]
	Both $\Sigma_{d}^{\pm}$ are closed in $M_{\nu}$ and invariant under $T^{n} \times \CC^{\ast}$, and $\Sigma_{d}^{-}$ is a compact subset of the core $\mcC$.
\end{lemma}

\begin{proof}
	For $\zeta \neq 0$, the proof of Lemma~\ref{lemma:limit-levels}(iv) shows that $( p , \zeta ) \in U_{d}$ if and only if $\Phi_{0}(p) < d$.
	For $\zeta = 0$, the function $\tau \mapsto \rho( \tau \cdot ( p , 0 ) ) = \Phi( \tau \cdot p )$ is constant with value $\Phi(p)$ if $p$ is fixed, and otherwise continuous and strictly increasing, with image $( \Phi_{0}(p) , \Phi_{\infty}(p) )$; since $d$ is not a critical value, and so differs from $\Phi(p)$ when $p$ is fixed, $( p , 0 ) \in U_{d}$ if and only if $\Phi_{0}(p) < d < \Phi_{\infty}(p)$.
	This gives the displayed decomposition, which is disjoint because $\Phi_{0} \leq \Phi_{\infty}$, and Lemma~\ref{lemma:limit-levels}(ii) gives the descriptions by cells.
	Since $d \not\in \CritVal(\Phi)$, Lemma~\ref{lemma:limit-levels}(iii) gives $\Sigma_{d}^{+} = \{ \Phi_{0} \geq d \}$ and $\Sigma_{d}^{-} = \{ \Phi_{\infty} \leq d \}$; by Lemma~\ref{lemma:limit-levels}(ii) and~(iii), these sets are closed and invariant, and $\Sigma_{d}^{-}$ lies in the core, which is compact (Subsection~\ref{subsection:subsection-2.3-arrangement-core-cuts}).
\end{proof}

We now identify the unstable strata as blades in $\widetilde{X}$.
The torus $\widetilde{K}$ acts linearly and diagonally on $T^{\ast}\CC^{N} \times \CC$ in the coordinates $( z , w , \zeta )$, so the following lemma applies to its orbits: every point in the closure of an orbit is, up to translation, a limit along a one-parameter subgroup.

\begin{lemma}
	\label{lemma:torus-orbit-closures}
	Let a torus $H$ act linearly on a finite-dimensional complex vector space $V$, diagonally in a basis, with characters $\chi_{1} , \ldots , \chi_{\ell}$.
	Let $v , v_{0} \in V$ with $v_{0}$ in the closure of $H \cdot v$.
	Then, there are a cocharacter $\theta$ of $H$ and an element $h_{0} \in H$ such that $\lim_{t \ra 0} \theta(t) \cdot v = h_{0} \cdot v_{0}$.
\end{lemma}

\begin{proof}
	Let $I$ and $J$ be the sets of indices of the non-zero coordinates of $v$ and $v_{0}$, respectively, so that $J \subseteq I$.
	For $a \in \ZZ_{\geq 0}^{I}$, write $u^{a} := \prod_{i} u_{i}^{a_{i}}$ for $u \in V$, and $\chi^{a} := \sum_{i} a_{i}\chi_{i}$.
	Since $( h \cdot v )^{a} = \chi^{a}(h) v^{a}$, the function $u \mapsto u^{a}v^{b} - u^{b}v^{a}$ vanishes on $H \cdot v$ whenever $\chi^{a} = \chi^{b}$, hence on its closure, so
	\begin{equation}
		\label{eqn:binomial-relation}
		v_{0}^{a}v^{b} = v_{0}^{b}v^{a}
		\qquad
		\text{whenever } a , b \in \ZZ_{\geq 0}^{I} \text{ and } \chi^{a} = \chi^{b}.
	\end{equation}
	There is no relation $\sum_{i \in I \setminus J} a_{i}\chi_{i} = \sum_{j \in J} c_{j}\chi_{j}$ with $a_{i} \in \QQ_{\geq 0}$ not all zero and $c_{j} \in \QQ$.
	Otherwise, after clearing denominators and writing $c = c^{+} - c^{-}$ with $c^{\pm} \geq 0$, applying \eqref{eqn:binomial-relation} to $a' := ( a , c^{-} )$ and $b' := ( 0 , c^{+} )$ in $\ZZ_{\geq 0}^{I \setminus J} \times \ZZ_{\geq 0}^{J}$ gives $v_{0}^{a'}v^{b'} = v_{0}^{b'}v^{a'}$, whose left side vanishes and whose right side does not.
	By Motzkin's transposition theorem over $\QQ$ \cite[Corollary~7.1k and the convention at the start of Chapter~7]{Schrijver86}, applied to the strict inequalities $\langle \theta , \chi_{i} \rangle > 0$ for $i \in I \setminus J$ and the inequalities $\pm\langle \theta , \chi_{j} \rangle \geq 0$ for $j \in J$, there is a rational cocharacter $\theta$ with $\langle \theta , \chi_{i} \rangle > 0$ for $i \in I \setminus J$ and $\langle \theta , \chi_{j} \rangle = 0$ for $j \in J$.
	After scaling, we may take $\theta$  to be integral.
	For this $\theta$, $\lim_{t \ra 0} \theta(t) \cdot v = v_{J}$, the vector obtained from $v$ by setting the coordinates outside $J$ to zero.
	Finally, put $\gamma := ( v_{0 , j} / v_{j} )_{j \in J} \in ( \CC^{\ast} )^{J}$.
	By \eqref{eqn:binomial-relation}, $\gamma^{c} = 1$ whenever $c \in \ZZ^{J}$ satisfies $\sum_{j} c_{j}\chi_{j} = 0$, so $\sum_{j} c_{j}\chi_{j} \mapsto \gamma^{c}$ is a well-defined homomorphism from a subgroup of the character lattice of $H$ to $\CC^{\ast}$.
	Since $\CC^{\ast}$ is divisible, this homomorphism extends to the whole lattice.
	Hence, $\gamma = ( \chi_{j}(h_{1}) )_{j \in J}$ for some $h_{1} \in H$, that is, $v_{0} = h_{1} \cdot v_{J}$, and we take $h_{0} := h_{1}^{-1}$.
\end{proof}

For each fixed component $\mcZ_{\beta}$, put $\widetilde{Z}_{\beta} := Q^{-1}(\mcZ_{\beta}) \times \{ 0 \} \subseteq \widetilde{X}$.
A blade is determined by a cocharacter and a component of its fixed locus, and the following lemma supplies the cocharacter of $\widetilde{K}$ that fixes $\widetilde{Z}_{\beta}$.

\begin{lemma}
	\label{lemma:stabilisers-fixed-components}
	There is a unique primitive cocharacter $\theta_{\beta}$ of $\widetilde{K}$, of the form $\theta_{\beta}(t) = ( \kappa_{\beta}(t) , t^{m_{\beta}} )$ with $m_{\beta} > 0$, such that $\theta_{\beta}(\CC^{\ast})$ is the identity component of the stabiliser in $\widetilde{K}$ of every point of $\widetilde{Z}_{\beta}$.
	For every $\widetilde{K}$-invariant open subset $X' \subseteq \widetilde{X}$ containing $\widetilde{Z}_{\beta}$, the set $\widetilde{Z}_{\beta}$ is a connected component of the fixed locus $( X' )^{\theta_{\beta}}$.
\end{lemma}

\begin{proof}
	A geometric quotient is submersive, with fibres the orbits \cite[Definition~0.6]{MFK1994}, here the connected $K_{\CC}$-orbits, so a decomposition of $Q^{-1}(\mcZ_{\beta})$ into two disjoint non-empty closed subsets would consist of saturated sets and would disconnect $\mcZ_{\beta}$.
	Therefore, $\widetilde{Z}_{\beta}$ is connected.
	For $y \in \widetilde{Z}_{\beta}$, the point $\widetilde{Q}(y)$ is fixed by the cutting action, so the stabiliser $\widetilde{K}_{y}$ projects onto the cutting factor, with finite kernel $( K_{\CC} )_{y}$.
	Its identity component $\widetilde{K}_{y}^{\circ}$ is therefore a one-dimensional torus with a unique primitive cocharacter of positive degree over the cutting factor.
	The stabiliser of a point of $T^{\ast}\CC^{N} \times \CC$ is the intersection of the kernels of the characters of its non-zero coordinates, so only finitely many such tori occur, and if $y \in \widetilde{Z}_{\beta}$ is fixed by one of them, $S$, then $S = \widetilde{K}_{y}^{\circ}$.
	Thus, $\widetilde{Z}_{\beta}$ is the disjoint union of the finitely many closed sets $\widetilde{Z}_{\beta} \cap \widetilde{X}^{S}$, and by connectedness only one of them is non-empty; $\theta_{\beta}$ is the primitive cocharacter of positive degree of the corresponding torus $S$.
	For the second claim, a point $( x , \zeta )$ fixed by $\theta_{\beta}$ has $\zeta = 0$, and $Q(x)$ is fixed by the fibre-rotation action, so $( X' )^{\theta_{\beta}}$ lies in $\cup_{\beta^{\prime}} \widetilde{Z}_{\beta^{\prime}}$, in which $\widetilde{Z}_{\beta}$ is open and closed by Proposition~\ref{prop:BB-cells}(i).
	Being connected, $\widetilde{Z}_{\beta}$ is a connected component of $( X' )^{\theta_{\beta}}$.
\end{proof}

With these cocharacters, the pulled-back attracting cells are blades of $\theta_{\beta}$, and the pulled-back repelling cells, at $\zeta = 0$, are blades of $\theta_{\beta}^{-1}$.
As the strata are removed one at a time, we state this inside an arbitrary $\widetilde{K}$-invariant open subset of $\widetilde{X}$.

\begin{proposition}
	\label{prop:cells-are-blades}
	Let $X' \subseteq \widetilde{X}$ be a $\widetilde{K}$-invariant open subset.
	\begin{enumerate}[(i)]
		\item if $Q^{-1}(A_{\beta}^{+}) \times \CC \subseteq X'$, then $Q^{-1}(A_{\beta}^{+}) \times \CC$ is the blade in $X'$ of $\theta_{\beta}$ and $\widetilde{Z}_{\beta}$;
		\item if $Q^{-1}(A_{\beta}^{-}) \times \{ 0 \} \subseteq X'$, then $Q^{-1}(A_{\beta}^{-}) \times \{ 0 \}$ is the blade in $X'$ of $\theta_{\beta}^{-1}$ and $\widetilde{Z}_{\beta}$.
	\end{enumerate}
\end{proposition}

\begin{proof}
	Let $\widetilde{x} \in X'$ lie over $( p , \zeta )$.
	If $\lim_{t \ra 0} \theta_{\beta}(t) \cdot \widetilde{x}$ exists in $X'$ and lies in $\widetilde{Z}_{\beta}$, then applying $\widetilde{Q}$ shows that the cutting orbit $\tau \cdot (p, \zeta)$ converges into $\mcZ_{\beta} \times \{0\}$ as $\tau \ra 0$, so $p \in A_{\beta}^{+}$.
	If instead $\lim_{t \ra 0} \theta_{\beta}(t)^{-1} \cdot \widetilde{x}$ exists in $X^{\prime}$ and lies in $\widetilde{Z}_{\beta}$, then this orbit converges into $\mcZ_{\beta} \times \{0\}$ as $\tau \ra \infty$, so $\zeta = 0$ and $p \in A_{\beta}^{-}$.

	Conversely, in case~(i), let $p \in A_{\beta}^{+}$ and $p_{0} := \lim_{\tau \ra 0} \tau \cdot p \in \mcZ_{\beta}$.
	The closure $W$ of $\widetilde{K} \cdot \widetilde{x}$ in $\widetilde{X}$ is $K_{\CC}$-invariant, so $\widetilde{Q}(W)$ is closed \cite[Definition~0.6]{MFK1994}.
	It contains the cutting orbit of $( p , \zeta )$, and hence its limit $( p_{0} , 0 )$; choose $\widetilde{x}_{0} \in W$ over $( p_{0} , 0 )$.
	Lemma~\ref{lemma:torus-orbit-closures}, applied in $T^{\ast}\CC^{N} \times \CC$, gives a cocharacter $\theta$ of $\widetilde{K}$ and $g_{0} \in \widetilde{K}$ with $\lim_{t \ra 0} \theta(t) \cdot \widetilde{x} = g_{0} \cdot \widetilde{x}_{0}$, which lies over $( p_{0} , 0 )$ and hence in $\widetilde{Z}_{\beta} \subseteq X'$.
	This limit is fixed by $\theta$, so $\theta = k\theta_{\beta}$ for some $k \in \ZZ$ by Lemma~\ref{lemma:stabilisers-fixed-components}.
	If $k > 0$, then $\widetilde{x}$ lies in the blade, and if $k = 0$, then $\widetilde{x} \in \widetilde{Z}_{\beta}$.
	If $k < 0$, then $\lim_{t \ra \infty} \theta_{\beta}(t) \cdot \widetilde{x}$ lies in $\widetilde{Z}_{\beta}$, so $\zeta = 0$ and $p \in A_{\beta}^{-}$ as above.
	By Proposition~\ref{prop:BB-cells}(iv), $p \in A_{\beta}^{+} \cap A_{\beta}^{-} = \mcZ_{\beta}$, so again $\widetilde{x} \in \widetilde{Z}_{\beta}$.
	Case~(ii) is the same, with $\zeta = 0$, $p \in A_{\beta}^{-}$, $p_{0} := \lim_{\tau \ra \infty} \tau \cdot p$, and the roles of $k > 0$ and $k < 0$ exchanged.
\end{proof}

It remains to check the weight hypothesis of Theorem~\ref{theorem:quantisation-theorem} along each stratum.
The sign of the weight depends only on whether the critical value $\Phi(\mcZ_{\beta})$ lies above or below $d$.

\begin{lemma}
	\label{lemma:stratum-weights}
	The $\theta_{\beta}$-weight of $\mcG_{d}$ along $\widetilde{Z}_{\beta}$ is $m_{\beta}( d - \Phi(\mcZ_{\beta}) )$.
	In particular, it is negative if $\Phi(\mcZ_{\beta}) > d$, and the $\theta_{\beta}^{-1}$-weight is negative if $\Phi(\mcZ_{\beta}) < d$.
\end{lemma}

\begin{proof}
	Write $\kappa := \kappa_{\beta}$ and $m := m_{\beta}$.
	By \eqref{eqn:hypertoric-kempf-ness}, every $K_{\CC}$-orbit in $\mcN^{\sstable}$ meets $\phi_{\HK}^{-1}(\nu , 0)$, so we may take $y = ( x , 0 ) \in \widetilde{Z}_{\beta}$ with $x = ( z , w )$ in that level set.
	The group $K_{\CC}$ acts with weight $+\alpha_{i}$ on $z_{i}$ and $-\alpha_{i}$ on $w_{i}$, so
	\[
		\theta_{\beta}(t) \cdot x = \big( t^{\langle \alpha_{i} , \kappa \rangle}z_{i} ,\, t^{m - \langle \alpha_{i} , \kappa \rangle}w_{i} \big)_{i}.
	\]
	Since $\theta_{\beta}$ fixes $y$, we have $\langle \alpha_{i} , \kappa \rangle = 0$ whenever $z_{i} \neq 0$ and $\langle \alpha_{i} , \kappa \rangle = m$ whenever $w_{i} \neq 0$.
	Pairing $\phi_{\RR}(x) = \nu$ with $\kappa$ therefore gives
	\[
		\langle \nu , \kappa \rangle = \tfrac{1}{2} \sum\nolimits_{i}\big( |z_{i}|^{2} - |w_{i}|^{2} \big) \langle \alpha_{i} , \kappa \rangle = -\tfrac{m}{2}\| w \|^{2} = -m\, \Phi(\mcZ_{\beta}).
	\]
	On the fibre of $L_{\nu}$ at $y$, the component $\kappa(t)$ acts by $\chi_{\nu}( \kappa(t) ) = t^{\langle \nu , \kappa \rangle}$ by \eqref{eqn:positive-character-linearisation}, and the cutting factor acts trivially, since its lift fixes the fibre coordinate.
	On the fibre of $\mcO_{\CC}(\varpi_{d})$, the cutting factor has geometric weight $+d$ and $K_{\CC}$ acts trivially.
	The $\theta_{\beta}$-weight of $\mcG_{d}$ at $y$ is therefore $-m\Phi(\mcZ_{\beta}) + md$, and it is the same at every point of the connected set $\widetilde{Z}_{\beta}$.
	Equivalently, the fibre-rotation action has geometric weight $\langle \nu , \kappa \rangle / m = -\Phi(\mcZ_{\beta})$ on the fibre of $\mcL_{\nu}$ over $\mcZ_{\beta}$, as the Kostant formula predicts with the conventions of Notation~\ref{notation:weight-conventions}.
\end{proof}

\begin{proof}[Proof of Theorem~\ref{theorem:cut-cohomology-every-level}]
	By Lemma~\ref{lemma:cut-unstable-locus}, $\widetilde{X} \setminus \widetilde{X}^{\sstable}$ is the disjoint union of the closed sets $Q^{-1}(\Sigma_{d}^{+}) \times \CC$ and $Q^{-1}(\Sigma_{d}^{-}) \times \{ 0 \}$.
	Filtering $\Sigma_{d}^{+}$ and then $\Sigma_{d}^{-}$ as in Proposition~\ref{prop:BB-cells}(v) lists strata $S_{1} , \ldots , S_{b}$, each of the form $Q^{-1}(A_{\beta}^{+}) \times \CC$ with $\Phi(\mcZ_{\beta}) > d$ or $Q^{-1}(A_{\beta}^{-}) \times \{ 0 \}$ with $\Phi(\mcZ_{\beta}) < d$, such that every union $S_{1} \cup \cdots \cup S_{k}$ is closed in $\widetilde{X}$.
	Put $X^{( k )} := \widetilde{X} \setminus ( S_{1} \cup \cdots \cup S_{k} )$, so that $X^{( 0 )} = \widetilde{X}$, $X^{( b )} = \widetilde{X}^{\sstable}$, and $S_{k + 1}$ is closed in $X^{( k )}$ and invariant under $\widetilde{K}$, under $T_{\CC}^{N}$ acting on the first factor, and under the residual fibre-rotation action.
	By Proposition~\ref{prop:cells-are-blades} and Lemma~\ref{lemma:stabilisers-fixed-components}, applied with $X' = X^{( k )}$, the stratum $S_{k + 1}$ is the blade in $X^{( k )}$ of $\theta_{\beta}$ or $\theta_{\beta}^{-1}$ and the connected component $\widetilde{Z}_{\beta}$ of the fixed locus $( X^{(k)} )^{\theta_{\beta}}$.
	By Lemma~\ref{lemma:stratum-weights}, the corresponding weight of $\mcG_{d}$ along $\widetilde{Z}_{\beta}$ is negative.
	Theorem~\ref{theorem:quantisation-theorem}, applied to $X^{(k)}$ and $S_{k + 1}$ for $k = 0 , \ldots , b - 1$, therefore gives a restriction isomorphism
	\[
		R\Gamma( [ \widetilde{X} / \widetilde{K} ] ; \mcG_{d} ) \overset{\sim}{\lra} R\Gamma( [ \widetilde{X}^{\sstable} / \widetilde{K} ] ; \mcG_{d} ) = R\Gamma(Y_{d} ; \mcL_{d}),
	\]
	which is equivariant for the residual fibre-rotation action and for $T_{\CC}^{N}$, linearised by $\chi_{\lambda}$, since both commute with $\widetilde{K}$ and preserve every $X^{( k )}$.
	The action of $T_{\CC}^{N}$ induces the residual $T^{n}$-action on $\widetilde{K}$-invariants.

	For the left-hand side,
	\[
		H^{q}( [ \widetilde{X} / \widetilde{K} ] ; \mcG_{d} ) = H^{q}(\widetilde{X} ; \mcG_{d})^{\widetilde{K}}
	\]
	by Theorem~\ref{theorem:quantisation-theorem}.
	The projection $\widetilde{X} \ra \mcN^{\sstable}$ is affine, with direct image $\mcO_{\mcN^{\sstable}}[y] = \bigoplus_{j \geq 0} \mcO_{\mcN^{\sstable}} y^{j}$, and cohomology on the Noetherian scheme $\mcN^{\sstable}$ commutes with direct sums.
	Hence, $H^{q}(\widetilde{X} ; \mcG_{d}) = H^{q}(\mcN^{\sstable} ; L_{\nu}) \otimes \CC[y] \cdot \one_{d}$, whose $K_{\CC}$-invariants are $H^{q}(M_{\nu} ; \mcL_{\nu}) \otimes \CC[y] \cdot \one_{d}$, again by \cite[(A.5)]{Teleman2000}.
	As in the proof of Lemma~\ref{lemma:filtration-invariant-sections}, a term $s y^{j} \otimes \one_{d}$, with $s$ of fibre-rotation weight $w$, has cutting weight $w + j - d$, so it is invariant if and only if $j = d - w \geq 0$, and it then has residual weight $w$.
	This gives the displayed isomorphisms of the theorem.
	For $q = 0$, the composite sends $s$ to the descent of $s y^{d - w} \otimes \one_{d}$, which is $\Psi_{d}(s)$; this proves~(i).
	Part~(ii) follows from Lemma~\ref{lemma:hypertoric-higher-vanishing} with $m = 1$, part~(iii) from~(i), (ii), and the definition \eqref{eqn:holomorphic-Euler-characteristic}, and part~(iv) from~(i) and $F_{d}H^{0}(M_{\nu} ; \mcL_{\nu}) / F_{d - 1}H^{0}(M_{\nu} ; \mcL_{\nu}) \cong H^{0}(M_{\nu} ; \mcL_{\nu})_{d}$.
\end{proof}

\begin{remark}[Attribution and Chamber Structure]
	\label{remark:quantisation-attribution}
	For compact pre-quantised Hamiltonian circle actions, the truncation of circle-weight multiplicities under symplectic cutting was established by Duistermaat, Guillemin, Meinrenken, and Wu \cite[Theorem~1]{DGMW1995}; see also Meinrenken \cite[Corollary~4.6]{Meinrenken1998}.
	Theorem~\ref{theorem:cut-cohomology-every-level}(i) gives a corresponding identification of section spaces in our non-compact setting.
	For non-compact pre-quantised Hamiltonian spaces with proper moment map, Weitsman defines a quantisation through the Riemann--Roch numbers of the reduced spaces, and proves its existence by symplectic cutting \cite[Definition~4.1, Theorems~1 and~4.3]{Weitsman2001}.
	Paradan shows that the index of a cut space agrees with this formal geometric quantisation on the weights below a bound that grows with the cut \cite[Definition~1.2 and Proposition~3.26]{Paradan2009}, and Ma and Zhang identify it with a transversal index \cite[Theorems~0.1 and~0.2]{MaZhang2014}.
	Note, however, that these are statements about indices of Dirac operators, whereas Theorem~\ref{theorem:cut-cohomology-every-level} concerns holomorphic sections and higher cohomology.

	For K\"ahler manifolds with a holomorphic $\CC^{\ast}$-action, a proper moment map that is bounded on one side, finitely many fixed components, and a free circle action on the level set, Wu proves the analogous isomorphisms \cite[Theorem~3.1(2)]{Wu2000}.
	Since $\Phi_{0}$ and $\Phi_{\infty}$ take their finite values in $\CritVal(\Phi)$, Lemma~\ref{lemma:cut-unstable-locus} shows that $U_{d}$~--~and hence the complex orbifold $Y_{d}$~--~depends only on the cut chamber containing $d$.
	Compare \cite[Theorems~3.2--3.3]{Hausel1998} and \cite[\S\S1.2--1.3]{BrionProcesi1990}; for variation of geometric invariant theory more generally, see \cite{Thaddeus1996, DolgachevHu1998}.

	If $d - 1$ and $d$ lie in the same chamber, then $\varpi_{d} = \varpi_{d - 1}\varpi_{1}$, and the cutting-invariant section $y \otimes \one_{1}$ descends on the common semistable locus to a section $s_{\zeta}$ of a line orbibundle on $Y_{d - 1} = Y_{d}$.
	This section vanishes to first order along $D_{d}$ and nowhere else, since on each local uniformising chart meeting $D_{d}$, $\zeta$ is a coordinate that cuts out $D_{d}$.
	Hence, $\mcL_{d} \cong \mcL_{d - 1} \otimes \mcO_{Y_{d}}(D_{d})$, the injection $j_{d}$ is multiplication by $s_{\zeta}$, and taking first Chern classes gives the Duistermaat--Heckman variation $[ \w_{d} ] = [ \w_{d - 1} ] + [ D_{d} ]$; compare \cite[Theorem~7.10 with Corollary~7.7]{Hausel1998}.
\end{remark}

%% file: sections/6-combinatorial-quantisation.tex
% !TeX root = ../main.tex
\section{Combinatorial Quantisation}
\label{section:section-6-combinatorial-quantisation}

At every positive regular integral level $d$, Theorem~\ref{theorem:cut-cohomology-every-level}(i) identifies $F_{d}H^{0}(M_{\nu} ; \mcL_{\nu})$ with $H^{0}(M_{\nu}^{\leq d} ; \mcL_{\nu}^{\leq d})$.
In this section, we compute the residual $T^{n}$-weight multiplicities in this filtration directly from the coordinate ring of the hypertoric quotient and obtain a weighted lattice-point description on the moment polyptych.
We then compare this description with the localisation formula of Theorem~\ref{theorem:character-hypertoric-orbifold} in the examples of \cite{BrownI}.

\subsection{Combinatorial Interpretation}
\label{subsection:subsection-6.1-combinatorial-interpretation}

We now compute $F_{d}H^{0}(M_{\nu} ; \mcL_{\nu})$ using the polynomial semi-invariants in Proposition~\ref{prop:positive-ray-semi-invariants}.
Recall the fixed integral lift $\lambda$ of $\nu$ from Section~\ref{section:section-2-conventions-geometric-input} and the affine-linear functionals $\ell_{i}$ of Subsection~\ref{subsection:subsection-2.3-arrangement-core-cuts}, and consider a monomial $z^{a}w^{b} \in P$, with $a, b \in \ZZ_{\geq 0}^{N}$.
A monomial $z^{a}w^{b}$ has $K_{\CC}$-character $\nu$ if and only if its exponent difference satisfies $\imath_{\ZZ}^{\ast}( a - b ) = \nu$.
Since $\nu = \imath_{\ZZ}^{\ast}(\lambda)$, this is equivalent to $a - b - \lambda \in \ker \imath_{\ZZ}^{\ast} = \Image \pi_{\ZZ}^{\ast}$.
Thus, there is a unique residual weight $x \in ( \mft_{\ZZ}^{n} )^{\ast}$ such that
\begin{equation}
	\label{eqn:exponent-differences}
	a_{i} - b_{i} = \langle x ,\, u_{i} \rangle + \lambda_{i} = \ell_{i}(x), \qquad
	\text{for each } i = 1 , \ldots , N.
\end{equation}
By \eqref{eqn:positive-ray-semi-invariants}, the corresponding class in $P / J_{\CC}P$ is then a section of $\mcL_{\nu}$ of residual $T^{n}$-weight $x$.

Fix a residual weight $x$.
For each $i$, the coordinatewise minimal non-negative exponents satisfying $a_{i} - b_{i} = \ell_{i}(x)$ are $a_{i}^{\min} = \max\{ 0 , \ell_{i}(x) \}$ and $b_{i}^{\min} = \max\{ 0 , -\ell_{i}(x) \}$.
Substituting these choices into $z^{a}w^{b}$ gives the monomial
\begin{equation}
	\label{eqn:minimal-monomial}
	p_{x} := \prod_{i = 1}^{N} z_{i}^{\max\{ 0 ,\, \ell_{i}(x) \}} w_{i}^{\max\{ 0 ,\, - \ell_{i}(x) \}},
\end{equation}
which we call the \textbf{minimal monomial} associated with the residual weight $x$.

Any other pair of exponents $( a_{i} , b_{i} ) \in \ZZ_{\geq 0}^{2}$ satisfying $a_{i} - b_{i} = \ell_{i}(x)$ must differ from the minimal choice by adding the same non-negative integer to both entries.
That is to say, for each $i = 1 , \ldots , N$, there is a unique integer $s_{i} \in \ZZ_{\geq 0}$ such that $a_{i} = \max\{ 0 ,\, \ell_{i}(x) \} + s_{i}$ and $b_{i} = \max\{ 0 ,\, - \ell_{i}(x) \} + s_{i}$.
Every monomial of $K_{\CC}$-character $\nu$ and residual weight $x$ therefore factors uniquely as $p_{x}$ times a monomial in the residual-invariant products $z_{i}w_{i}$:
\[
	z^{a}w^{b} = \prod_{i = 1}^{N}
	z_{i}^{\max\{ 0 ,\, \ell_{i}(x) \} + s_{i}} w_{i}^{\max\{ 0 ,\, - \ell_{i}(x) \} + s_{i}} = p_{x} \cdot \prod_{i = 1}^{N} ( z_{i}w_{i} )^{s_{i}},
\]
so $p_{x}$ records the part of the monomial forced by the residual weight $x$, while the exponents $s_{i}$ count the additional invariant factors.

The same factorisation describes the entire residual-weight-$x$ subspace before imposing the moment-map equations.
Set $h_{i} := z_{i}w_{i}$ and $B_{0} := \CC[h_{1} , \ldots , h_{N}]$, and let $P_{x} \subseteq P$ be the span of the monomials satisfying \eqref{eqn:exponent-differences}.
Since multiplication by $h_{i}$ preserves the exponent differences, $P_{x}$ is a $B_{0}$-module.
Furthermore, the unique factorisation identifies its monomial basis with $p_{x}$ times the monomial basis of $B_{0}$, giving the $B_{0}$-module isomorphism $B_{0} \overset{\sim}{\ra} P_{x}$, $f \mapsto p_{x}f$.
Thus, $P_{x} = p_{x}B_{0}$ is free of rank one, a description we will use when quotienting by the complex moment-map equations in the variables $h_{i}$.

To impose the cutting condition, adjoin the auxiliary coordinate $y$.
Then a monomial of residual weight $x$ on $T^{\ast}\CC^{N} \times \CC$ has the unique normal form
\begin{equation}
	\label{eqn:monomial-normal-form}
	\psi = z^{a}w^{b} y^{c} = p_{x} \cdot \left( \prod_{i = 1}^{N} ( z_{i} w_{i} )^{s_{i}} \right) y^{c}, \qquad c \in \ZZ_{\geq 0}.
\end{equation}
Under the pullback convention of Notation~\ref{notation:representation-convention}, the cutting weights are $\wt_{\CC^{\ast}}(z_{i}) = 0$ and $\wt_{\CC^{\ast}}(w_{i}) = \wt_{\CC^{\ast}}(y) = 1$.
On $P[y]$, we call this weight the \textbf{cut degree}, $\deg_{\cut} := \wt_{\CC^{\ast}}$; it is distinct from the tensor-power degree in the section rings of Subsection~\ref{subsection:subsection-5.2-algebraic-cutting}.
In particular, each $h_{i} = z_{i}w_{i}$ has residual weight zero but cut degree one.

The exponent of $w_{i}$ in $p_{x}$ is $\max\{ 0 ,\, - \ell_{i}(x) \}$, and therefore $\deg_{\cut}(p_{x}) = \sum_{i = 1}^{N} \max\{ 0 ,\, - \ell_{i}(x) \} = \mcS(x)$, where $\mcS$ is the slack function of \eqref{eqn:slack-function}.
Thus, $\mcS(x)$ is the minimum cut degree required to realise the residual weight $x$; every additional $h_{i}$- or $y$-factor increases the degree by one.
For the monomial in \eqref{eqn:monomial-normal-form}, this gives $\deg_{\cut}(\psi) = \wt_{\CC^{\ast}}(\psi) = \mcS(x) + \sum\nolimits_{i = 1}^{N} s_{i} + c$.
Recall that the constant section $\one_{d}$ of $\mcO_{\CC}(\varpi_{d})$ has pullback weight $-d$.
The section $\psi \otimes \one_{d}$ is therefore cutting-invariant if and only if $s_{1} + \ldots + s_{N} + c = d - \mcS(x)$.
Clearly, there are no solutions when $\mcS(x) > d$; otherwise, counting monomials of degree $d - \mcS(x)$ in the $N + 1$ variables $h_{1} , \ldots , h_{N} , y$ gives $\tbinom{d - \mcS(x) + N}{N}$.
This count is taken before imposing the complex moment-map equations; to obtain the section multiplicity, we must determine the remaining independent classes.

The complex moment-map equations are linear in the $h_{i}$ variables.
Let $J_{\CC} \subseteq B_{0}$ be the ideal generated by the components of $\phi_{\CC}(z , w) = \imath^{\ast}(h_{1} , \ldots , h_{N}) = \sum_{i = 1}^{N} h_{i}\alpha_{i}$.
That is, after choosing a basis $\kappa_{1} , \ldots , \kappa_{r}$ of $\mfk_{\CC}$, where $r = \dim_{\CC} \mfk_{\CC}$, the ideal $J_{\CC}$ is generated by the homogeneous linear forms $g_{a} := \langle \phi_{\CC} ,\, \kappa_{a} \rangle = \sum\nolimits_{i = 1}^{N} \langle \alpha_{i} ,\, \kappa_{a} \rangle h_{i}$ for $a = 1 , \ldots , r$.
Its extension $J_{\CC}P$ is the complex moment-map ideal of Proposition~\ref{prop:positive-ray-semi-invariants}.

Equip $P$ with the exponent-difference grading given by the character lattice $( \mft_{\ZZ}^{N} )^{\ast}$, determined by $\deg_{\diff}(z_{i}) = \epsilon_{i}$ and $\deg_{\diff}(w_{i}) = -\epsilon_{i}$, where $\epsilon_{1} , \ldots , \epsilon_{N}$ is the dual basis introduced in Subsection~\ref{subsection:subsection-2.1-hypertoric-quotients}.
In this grading, we have $\deg_{\diff}(z^{a}w^{b}) = \sum_{i = 1}^{N} ( a_{i} - b_{i} ) \epsilon_{i}$, and $P_{x}$ is the homogeneous component of degree $\lambda + \pi_{\ZZ}^{\ast}(x)$.
Since the generators of $J_{\CC}$ have exponent-difference degree zero, the component of $J_{\CC}P$ of degree $\lambda + \pi_{\ZZ}^{\ast}(x)$ is $( J_{\CC}P ) \cap P_{x} = J_{\CC}P_{x} = p_{x}J_{\CC}$, where the last equality uses $P_{x} = p_{x}B_{0}$.
Multiplication by $p_{x}$ therefore induces the isomorphism $B_{0}/J_{\CC} \overset{\sim}{\ra} ( P / J_{\CC}P )_{x}$ given by $\widebar{f} \mapsto \widebar{p_{x}f}$.

Since $\imath^{\ast} : ( \mft^{N} )^{\ast} \ra \mfk^{\ast}$ is surjective and $\dim_{\CC} \mfk_{\CC} = r$, the generators $g_{1} , \ldots , g_{r}$ are $r$ independent linear forms in the $h_{1} , \ldots , h_{N}$ variables.
Choose classes $q_{1} , \ldots , q_{n}$ forming a basis of the quotient $\Span\{ h_{1} , \ldots , h_{N} \} / \Span\{ g_{1} , \ldots , g_{r} \}$, whose dimension is $N - r = n$.
Then $B_{0}/J_{\CC} \cong \CC[q_{1} , \ldots , q_{n}]$, where each $q_{j}$ is represented by a homogeneous linear combination of the variables $h_{i} = z_{i}w_{i}$ and therefore has cut degree one.
After adjoining the auxiliary coordinate function $y$, we obtain $( ( P / J_{\CC} P )[y] )_{x} \cong p_{x} \CC[q_{1} , \ldots , q_{n} , y]$.
Hence, the residual-weight-$x$, cut-degree-$d$ subspace has the basis
\begin{equation}
	\label{eqn:weight-x-basis}
	\{ p_{x} q_{1}^{k_{1}} \cdots q_{n}^{k_{n}} y^{c} ~|~ k_{1} , \ldots , k_{n} , c \in \ZZ_{\geq 0}, ~k_{1} + \ldots + k_{n} + c = d - \mcS(x) \},
\end{equation}
and this set is empty when $\mcS(x) > d$.

The basis \eqref{eqn:weight-x-basis} determines which residual weights occur as well as their multiplicities.
The minimum cut degree is $\mcS(x)$, and the remaining degree is distributed among the $n + 1$ variables $q_{1} , \ldots , q_{n} , y$.
The following theorem also establishes the non-negativity and finite-dimensionality used in Subsection~\ref{subsection:subsection-5.2-algebraic-cutting}.

\begin{theorem}
	\label{theorem:combinatorial-multiplicity}
	For each $d \in \ZZ_{\geq 0}$, the residual $T^{n}$-weight $x \in ( \mft_{\ZZ}^{n} )^{\ast}$ occurs in the filtration piece $F_{d}H^{0}(M_{\nu} ; \mcL_{\nu})$ with multiplicity
	\begin{equation}
		\label{eqn:multiplicity}
		\mult_{d}(x) = \begin{cases}
			\binom{d - \mcS(x) + n}{n}, \qquad & \text{if } \mcS(x) \leq d, \\
			0, \qquad                          & \text{if } \mcS(x) > d,
		\end{cases}
	\end{equation}
	and no negative fibre-rotation weights occur in $H^{0}(M_{\nu} ; \mcL_{\nu})$.
	Furthermore, the support of $\mult_{d}$ is the finite set $\supp( \mult_{d} ) = \{ x \in ( \mft_{\ZZ}^{n} )^{\ast} \mid \mcS(x) \leq d \}$, and
	\[
		F_{d} H^{0}(M_{\nu} ; \mcL_{\nu}) \cong \bigoplus_{x \in ( \mft_{\ZZ}^{n} )^{\ast}} \CC_{x}^{\oplus\mult_{d}(x)}, \qquad \dim_{\CC} F_{d} H^{0}(M_{\nu} ; \mcL_{\nu}) = \sum_{x \in ( \mft_{\ZZ}^{n} )^{\ast}} \mult_{d}(x),
	\]
	where $\CC_{x}$ denotes the one-dimensional residual $T^{n}$-representation of weight $x$.

	When $d$ is a positive regular integral level, this support is $|\Delta_{\nu}^{\leq d}| \cap ( \mft_{\ZZ}^{n} )^{\ast}$ by \cite[Proposition~4.6]{BrownI}.
	In that case, Theorem~\ref{theorem:cut-cohomology-every-level}(i) identifies $F_{d} H^{0}(M_{\nu} ; \mcL_{\nu})$ with $H^{0}(M_{\nu}^{\leq d} ; \mcL_{\nu}^{\leq d})$, and the same multiplicity and dimension formulae hold for sections of $\mcL_{\nu}^{\leq d} \ra M_{\nu}^{\leq d}$.
\end{theorem}

\begin{proof}
	Proposition~\ref{prop:positive-ray-semi-invariants} with $m = 1$ identifies $( P/J_{\CC}P )_{\nu}$ with $H^{0}(M_{\nu} ; \mcL_{\nu})$.
	By \eqref{eqn:filtration-invariant-sections}, adjoining $y$ and tensoring the cut-degree-$d$ component with $\one_{d}$ gives the invariant sections that correspond to $F_{d}H^{0}(M_{\nu} ; \mcL_{\nu})$.

	Fix a residual weight $x$.
	The basis \eqref{eqn:weight-x-basis} shows that its multiplicity is the number of non-negative integer solutions of $k_{1} + \ldots + k_{n} + c = d - \mcS(x)$.
	Clearly, if $\mcS(x) > d$, there are no solutions; otherwise, counting monomials of degree $d - \mcS(x)$ in $n + 1$ variables gives $\mult_{d}(x) = \binom{d - \mcS(x) + n}{n}$, proving \eqref{eqn:multiplicity}.
	The same description shows that no negative fibre-rotation weights occur: $\mcS(x) \geq 0$ always holds, and the additional generators have positive cut degree.

	Choose a regular level $\delta^{\prime} > d$, which is possible because $\Phi$ has only finitely many critical values by \cite[Proposition~4.7]{BrownI}.
	By \cite[Proposition~4.6]{BrownI}, the set $\{ \mcS \leq \delta^{\prime} \}$ is compact.
	Since $\mcS$ is continuous, $\{ \mcS \leq d \}$ is a closed subset of this compact set and therefore contains only finitely many lattice points.
	The filtration piece is consequently finite-dimensional, and summing its weight spaces gives the stated decomposition and dimension formula.
	When $d$ is a positive regular integral level, the weight support is the set of lattice points in $|\Delta_{\nu}^{\leq d}|$ by \cite[Proposition~4.6]{BrownI}.
	Theorem~\ref{theorem:cut-cohomology-every-level}(i) then identifies the filtration piece with $H^{0}(M_{\nu}^{\leq d} ; \mcL_{\nu}^{\leq d})$.
	The same multiplicity and dimension formulae therefore apply there.
\end{proof}

\begin{remark}[Relation to the Monopole Formula]
	\label{remark:monopole-formula}
	By Proposition~\ref{prop:positive-ray-semi-invariants} and the description of $( P / J_{\CC}P )_{x}$ before \eqref{eqn:weight-x-basis}, the residual-weight-$x$ part of $H^{0}(M_{\nu} ; \mcL_{\nu})$ is $p_{x}\CC[q_{1} , \ldots , q_{n}]$, where $p_{x}$ has fibre-rotation weight $\mcS(x)$ and each $q_{j}$ has weight one.
	Hence, the character of $H^{0}(M_{\nu} ; \mcL_{\nu})$ for $T^{n}$ and the fibre-rotation action is
	\[
		\sum_{x \in ( \mft_{\ZZ}^{n} )^{\ast}} \frac{\tau^{\mcS(x)}}{( 1 - \tau )^{n}}\, e^{\langle x ,\, \xi \rangle}.
	\]
	Since $|\ell| = 2 \max\{ 0 , -\ell \} + \ell$, we have
	\[
		\sum_{i} |\ell_{i}(x)|
		= 2\mcS(x) + \left\langle x , \sum_{i} u_{i} \right\rangle + \sum_{i} \lambda_{i}.
	\]
	Thus, the substitution $	\tau = t^{2}$ and $e^{\langle x,\xi \rangle} \longmapsto z^{x} t^{\langle x,\sum_{i} u_{i} \rangle}$ turns this series into
	\[
		t^{-\sum_{i} \lambda_{i}} \sum_{x}
		\frac{z^{x}t^{\sum_{i} |\ell_{i}(x)|}}{( 1 - t^{2} )^{n}}.
	\]
	Up to the overall factor $t^{-\sum_{i} \lambda_{i}}$, this is Nakajima's formula for the sections of the line bundle determined by a background flux, in the grading by dilation of $T^{\ast}\CC^{N}$ \cite[\S5(ii)]{Nakajima2016}, and so the abelian case of the monopole formula with background flux \cite[\S2.2]{CremonesiHananyMekareeyaZaffaroni2014}.
	His proof uses the same minimal monomials.
	For $\nu = 0$, the corresponding Hilbert series of the affine hypertoric variety is computed in \cite[\S6]{CremonesiHananyZaffaroni2014}.
	Braverman, Finkelberg, and Nakajima realise the partial resolution of the Coulomb branch as the $\Proj$ of a ring of semi-invariants whose graded pieces are counted by the same formula \cite[Proposition~3.25 and Remark~3.26]{BravermanFinkelbergNakajima2018}.
	The argument above adds Proposition~\ref{prop:positive-ray-semi-invariants}, which identifies the sections of $\mcL_{\nu}^{\otimes m}$ with the semi-invariants of weight $m\nu$ for every $m \geq 0$; the minimal-monomial argument depends on this step.
	It also adds the fibre-rotation grading and  the description of the support by the slack function.
\end{remark}

The following corollary expresses the character of $H^{0}(M_{\nu}^{\leq d} ; \mcL_{\nu}^{\leq d})$ as a weighted sum over the lattice points of $|\Delta_{\nu}^{\leq d}|$.

\begin{corollary}[Weighted Lattice-Point Character]
	\label{cor:weighted-lattice-point-character}
	Let $d$ be a  positive regular integral level.
	Then
	\begin{equation}
		\label{eqn:weighted-lattice-point-character}
		\operatorname{char}_{T^{n}}
		H^{0}(M_{\nu}^{\leq d} ; \mcL_{\nu}^{\leq d})(e^{\xi}) =
		\sum_{x \in |\Delta_{\nu}^{\leq d}| \cap ( \mft_{\ZZ}^{n} )^{\ast}}
		\binom{d - \mcS(x) + n}{n} e^{\langle x,\, \xi \rangle}.
	\end{equation}
	Moreover, setting $\xi = 0$ gives
	\[
		\dim_{\CC} H^{0}(M_{\nu}^{\leq d} ; \mcL_{\nu}^{\leq d}) =
		\sum_{x \in |\Delta_{\nu}^{\leq d}| \cap ( \mft_{\ZZ}^{n} )^{\ast}}
		\binom{d - \mcS(x) + n}{n}.
	\]
	By Theorem~\ref{theorem:cut-cohomology-every-level}(iii), the character in \eqref{eqn:weighted-lattice-point-character} agrees with the Kawasaki localisation formula in Theorem~\ref{theorem:character-hypertoric-orbifold}.
\end{corollary}

\begin{proof}
	By Theorem~\ref{theorem:cut-cohomology-every-level}(i), the representation $H^{0}(M_{\nu}^{\leq d} ; \mcL_{\nu}^{\leq d})$ has the weight decomposition given by Theorem~\ref{theorem:combinatorial-multiplicity}.
	Summing the corresponding characters gives \eqref{eqn:weighted-lattice-point-character}, and setting $\xi = 0$ gives the dimension formula.
	Theorem~\ref{theorem:cut-cohomology-every-level}(iii) identifies this section-space character with the equivariant index, which Theorem~\ref{theorem:character-hypertoric-orbifold} computes by orbifold localisation.
\end{proof}

\begin{remark}[Enlarged Moment Image]
	\label{remark:enlarged-cut-moment-image}
	Suppose that $n = 1$ and that $d$ is a positive regular integral level in the outer cut chamber.
	In this case, the multiplicity $\mult_{d}(x)$ counts the lattice points above $x$ in a two-dimensional moment polygon, obtained by using the auxiliary-circle moment map as a second coordinate.
	Rotation of the auxiliary coordinate $\zeta$ descends to a Hamiltonian circle action on $M_{\nu}^{\leq d}$, with moment map $h$ induced by $\tfrac{1}{2}|\zeta|^{2} = d - \Phi$ on $\rho^{-1}(d)$.
	Together with the residual $T^{1}$-action, this gives an effective Hamiltonian $T^{1} \times U_{1}$-action on the four-real-dimensional cut, making it a symplectic toric orbifold.

	With the residual weight lattice $(\mft_{\ZZ}^{1})^{\ast}$ identified with $\ZZ$, the joint moment map $(\mu_{\RR} , h)$ has image
	\begin{equation}
		\label{eqn:enlarged-cut-moment-image}
		\mcP_{\nu}^{\leq d} = \{ ( x , h ) \in \RR^{2} ~|~ 0 \leq h \leq d - \mcS(x) \}.
	\end{equation}
	Indeed, $h \leq d - \mcS(x)$ by \cite[Lemma~4.3]{BrownI}, and the graph $h = d - \mcS(x)$ is attained on the cut extended core by \cite[Proposition~4.4]{BrownI}.
	Its endpoints lie on $h = 0$, so orbifold convexity \cite[Theorem~5.2]{LermanTolman1997} fills the region between this graph and the horizontal axis.
	For each integral $x$ with $\mcS(x) \leq d$, the lattice points in the vertical slice are $\{ ( x , j ) ~|~ j \in \ZZ , ~ 0 \leq j \leq d - \mcS(x) \}$, whose number is $d - \mcS(x) + 1 = \mult_{d}(x)$.
	Thus, the multiplicity at $x$ counts the lattice points that project to $x$, and the weighted count becomes
	\[
		\dim_{\CC}H^{0}(M_{\nu}^{\leq d} ; \mcL_{\nu}^{\leq d}) = \#( \mcP_{\nu}^{\leq d} \cap \ZZ^{2} ).
	\]
	We illustrate this interpretation in Figure~\ref{figure:rank-one-enlarged-moment-polytope}, using the $M_{\nu} \cong T^{\ast}\CC\PP^{1}$ case as an example.

	For $n \geq 2$, the same moment-image description holds in $( \mft^{n} )^{\ast} \times \RR$, but the enlarged torus has dimension $n + 1 < 2n$, so this action is not toric.
	Hence, the $d - \mcS(x) + 1$ lattice heights do not in general account for the full multiplicity $\binom{d - \mcS(x) + n}{n}$ when $n \geq 2$.

	\begin{figure}[htbp]
		\centering
		\resizebox{\ifdim\linewidth<250pt\linewidth\else250pt\fi}{!}{%
			\begin{tikzpicture}[
					x=0.50cm,
					y=0.50cm,
					font=\small
				]

				\filldraw[
					fill=gray!70,
					draw=gray!100,
					line width=1.2pt
				]
				(-4 , 0) -- (0 , 4) -- (6 , 4) -- (10,0) -- cycle;

				\foreach \X/\H in {
						-2/2,
						2/4,
						4/4,
						6/4,
						8/2
					}{
						\draw[
							gray!75!black,
							dashed,
							line width=0.85pt
						] (\X , 0) -- (\X,\H);
					}

				\draw[->, line width=0.85pt]
				(-5.7 , 0) -- (11.1 , 0);
				\draw[->, line width=0.85pt]
				(0 , -0.10) -- (0 , 4.95);

				\node[left] at (-0.24 , 5.15) {$h$};
				\node[below] at (11.35 , -0.12) {$x$};

				\foreach \X/\Y in {
						-4/0,
						-2/0, -2/2,
						0/0,  0/2,  0/4,
						2/0,  2/2,  2/4,
						4/0,  4/2,  4/4,
						6/0,  6/2,  6/4,
						8/0,  8/2,
						10/0
					}{
						\filldraw[
							fill=black!78,
							draw=black!90,
							line width=0.3pt
						] (\X,\Y) circle[radius=1.65pt];
					}

				\node
				at (-2 , -0.80) {$\Delta_{1}^{\leq d}$};
				\node
				at (3 , -0.80) {$\Delta_{\emptyset}$};
				\node
				at (8 , -0.80) {$\Delta_{2}^{\leq d}$};

			\end{tikzpicture}
		}
		\caption{
			The enlarged moment polygon $\mcP_{\nu}^{\leq d}$ for $M_{\nu} \cong T^{\ast}\CC\PP^{1}$, whose upper boundary is the graph $h = d - \mcS(x)$.
			For each integral residual weight $x \in |\Delta_{\nu}^{\leq d}|$, the vertical slice contains $\mult_{d}(x) = d - \mcS(x) + 1$ lattice points, and summing these multiplicities gives the dimension of the cut section space.
		}
		\figalttext{
			A shaded trapezoid resting on a horizontal axis labelled x, with a vertical axis labelled h through its upper-left corner.
			Its left side rises at slope one to a flat top, and its right side falls at slope one back to the axis.
			Black dots mark the lattice points of the trapezoid, and dashed vertical segments run from the axis to the upper boundary through the dots in each column.
			Below the axis, the intervals under the left side, the top, and the right side are labelled Delta sub 1 with superscript at most d, Delta sub empty set, and Delta sub 2 with superscript at most d.
		}
		\label{figure:rank-one-enlarged-moment-polytope}
	\end{figure}
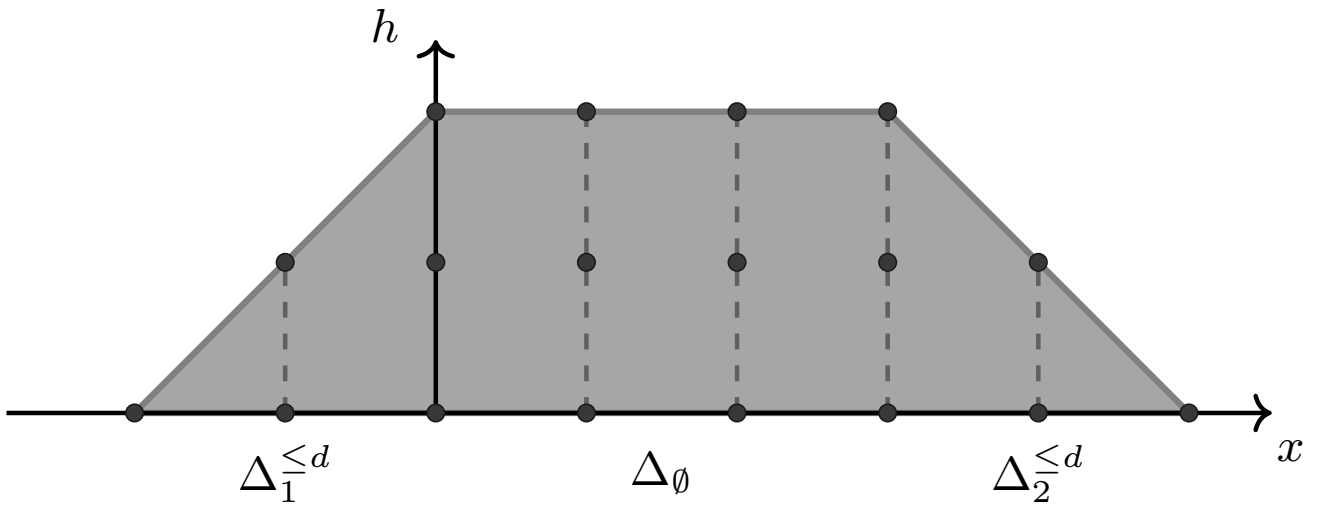
\end{remark}

\begin{remark}[Reduced Spaces]
	\label{remark:reduced-spaces}
	The multiplicities of Theorem~\ref{theorem:combinatorial-multiplicity} have a reduced-space description in every rank, which for $n = 1$ recovers the vertical slices of Remark~\ref{remark:enlarged-cut-moment-image}.
	Let $d$ be a positive regular integral level, let $x \in ( \mft_{\ZZ}^{n} )^{\ast}$ satisfy $\mcS(x) < d$, and set $r := d - \mcS(x)$.
	Let $A_{d} := ( R~\#_{\CC}~S )^{\CC^{\ast}}$ be the invariant ring of Definition~\ref{def:algebraic-cut}, graded by the tensor degree $m$, and let $( A_{d} )_{m , mx}$ denote the residual-weight-$mx$ part of its degree-$m$ piece.
	In degree $m$, the lift is $m\lambda$ and the cut degree is $md$.
	Replacing $( \lambda , x )$ by $( m\lambda , mx )$ multiplies every value $\ell_{i}(x)$ by $m$, so the minimal monomial at $mx$ is $p_{x}^{m}$, of cut degree $m\mcS(x)$, and the argument leading to \eqref{eqn:weight-x-basis} gives $( A_{d} )_{m , mx} = p_{x}^{m}\, \CC[q_{1} , \ldots , q_{n} , y]_{mr}$.
	Multiplication by $p_{x}^{m}$ is injective on these polynomials, again by that argument, and $p_{x}^{m}p_{x}^{l} = p_{x}^{m + l}$, so the pieces $( A_{d} )_{m , mx}$ for $m \geq 0$ form a graded ring isomorphic to the $r$-th Veronese subring of $\CC[q_{1} , \ldots , q_{n} , y]$.
	The GIT quotient $M_{\nu}^{\alg , \leq d} \sslash_{x} T_{\CC}^{n}$ by the complexified residual torus at $x$ is the $\Proj$ of this ring, and is therefore $\CC\PP^{n}$, with homogeneous coordinates $q_{1} , \ldots , q_{n} , y$ and polarisation $\mcO(r)$.
	The auxiliary circle acts on it by rotating $y$ and fixing $q_{1} , \ldots , q_{n}$.
	Since the $q_{j}$ may be taken to be the components of $\mu_{\CC}$, the fixed-point set of the auxiliary circle on the quotient consists of the image $q = 0$ of the extended core and the quotient $\{ y = 0 \}$ of the cut divisor.
	In degree one, Theorem~\ref{theorem:combinatorial-multiplicity} gives $\mult_{d}(x) = \dim_{\CC} H^{0}(\CC\PP^{n} ; \mcO_{\CC\PP^{n}}(r))$, as quantisation commuting with reduction predicts when $x$ is a regular value of the residual moment map on the cut \cite[Corollary~1.2]{Meinrenken1998}.
	At a critical value, the quotient stack can have positive-dimensional stabilisers, and the statement concerns the coarse quotient only.

	For $n = 1$ and an outer-chamber level, the auxiliary circle acts on the quotient $\CC\PP^{1}$ as its torus, and the $r + 1$ lattice points of its moment image are those above $x$ in \eqref{eqn:enlarged-cut-moment-image}.
	For $n \geq 2$, the auxiliary circle acts on $\CC\PP^{n}$ only through a single circle subgroup of its $n$-torus, so the $r + 1$ lattice points of its moment image do not account for the multiplicity $\binom{r + n}{n}$, which counts the lattice points of the standard $n$-simplex dilated by $r$.
	This description is algebraic, and the reduced symplectic form on the quotient need not be the Fubini--Study form.
\end{remark}

\subsection{Worked Examples}
\label{subsection:subsection-6.2-worked-examples}

We now continue the examples of \cite{BrownI}, comparing the localisation formula of Theorem~\ref{theorem:character-hypertoric-orbifold} with the weighted lattice-point character of Corollary~\ref{cor:weighted-lattice-point-character}.
In rank one, the cuts of $T^{\ast}\CC\PP^{1}$ are smooth, while the outer-chamber cuts of the $A_{2}$-resolution have a fixed point with local group $\ZZ / 2\ZZ$; in rank two, the cuts of $T^{\ast}\CC\PP^{2}$ and $T^{\ast}( \CC\PP^{1} \times \CC\PP^{1} )$ are smooth.
With $\chi_{v}$ as in Section~\ref{section:section-4-equivariant-quantisation-cut-spaces}, we write
\begin{equation}
	\label{eqn:cut-character-shorthand}
	\chi_{d}(e^{\xi}) := \chi_{T^{n}}(M_{\nu}^{\leq d} ; \mcL_{\nu}^{\leq d})(e^{\xi}) = \sum\nolimits_{v \in \Verts(\Delta_{\nu}^{\leq d})} \chi_{v}(e^{\xi}).
\end{equation}
By Theorem~\ref{theorem:cut-cohomology-every-level}, $\chi_{d}$ is the character of $H^{0}( M_{\nu}^{\leq d} ; \mcL_{\nu}^{\leq d} )$ at every positive regular integral level, so its value at the identity is the dimension.
We write $t = e^{\xi}$ in rank one and $(t_{1}, t_{2}) = (e^{\xi_{1}}, e^{\xi_{2}})$ in rank two.

\begin{example}[$T^{\ast}\CC\PP^{1}$, Continued]
	\label{example:TCP1-example-2}
	We continue \cite[Example~4.10]{BrownI}, with affine-linear functionals $\ell_{1}(x) = x$ and $\ell_{2}(x) = m - x$, where $m , d \in \ZZ_{> 0}$.
	The cut polyptych has components $[ -d , 0 ]$, $[ 0 , m ]$, and $[ m , m + d ]$.
	The four fixed points of the cut have trivial local groups, and their tangent weights are shown in Figure~\ref{fig:isotropy-TCP1}; the cutting weight occurs twice at each boundary vertex.

	\begin{figure}[H]
		\centering
		\begin{tikzpicture}[x=1cm , y=1cm]

			\tikzset{
				polyboundary/.style={line width=0.9pt},
				weight/.style={->,line width=0.8pt},
				cutweight/.style={->,line width=0.8pt},
				fixedpoint/.style={font=\scriptsize},
				weightlabel/.style={font=\scriptsize},
				regionlabel/.style={font=\scriptsize},
			}

			\coordinate (b1) at (-6.75 , 0);
			\coordinate (v1) at (-2.55 , 0);
			\coordinate (v2) at ( 2.55 , 0);
			\coordinate (b2) at ( 6.75 , 0);

			\foreach \p in {v1, v2}
			\draw[dashed, dash pattern=on 1.5pt off 1.5pt, line width=0.35pt, black!55]
			($(\p)+(0 , 0.25)$) -- ($(\p)+(0 , -0.52)$);

			\foreach \p in {b1, b2}
			\draw[dashed, dash pattern=on 1.5pt off 1.5pt, line width=0.35pt, black!55]
			($(\p)+(0 , 0.25)$) -- ($(\p)+(0 , -0.52)$);

			\draw[polyboundary] (b1) -- (b2);
			\foreach \p in {b1,v1,v2,b2}
			\fill (\p) circle (1.35pt);

			\draw[weight]
			($(v1)+(0 , 0.25)$) -- ++(0.85 , 0)
			node[midway , above right=2pt, weightlabel] {$\varrho_{0, 1} = +1$};
			\draw[weight]
			($(v1)+(0 , -0.25)$) -- ++(-0.85 , 0)
			node[midway , below left=2pt, weightlabel] {$\varsigma_{0, 1} = -1$};

			\draw[weight]
			($(v2)+(0 , 0.25)$) -- ++(-0.85 , 0)
			node[midway , above left=2pt, weightlabel] {$\varrho_{m, 2} = -1$};
			\draw[weight]
			($(v2)+(0 , -0.25)$) -- ++(0.85 , 0)
			node[midway , below right=2pt, weightlabel] {$\varsigma_{m, 2} = +1$};

			\draw[cutweight]
			($(b1)+(0 , 0.22)$) -- ++(0.85 , 0);
			\draw[cutweight]
			($(b1)+(0 , -0.22)$) -- ++(0.85 , 0);

			\node[weightlabel , above=3pt] at ($(b1)+(0.43 , 0.22)$)
			{$\vartheta_{-d} = +1$};

			\draw[cutweight]
			($(b2)+(0 , 0.22)$) -- ++(-0.85 , 0);
			\draw[cutweight]
			($(b2)+(0 , -0.22)$) -- ++(-0.85 , 0);

			\node[weightlabel , above=3pt] at ($(b2)+(-0.43 , 0.22)$)
			{$\vartheta_{m + d} = -1$};

			\node[fixedpoint , below=17pt] (lv1) at (v1) {$0$};
			\node[fixedpoint , anchor=base] at (b1 |- lv1.base) {$-d$};
			\node[fixedpoint , anchor=base] at (v2 |- lv1.base) {$m$};
			\node[fixedpoint , anchor=base] at (b2 |- lv1.base) {$m + d$};

			\coordinate (r1) at ($(b1)!0.5!(v1)$);
			\coordinate (r2) at ($(v1)!0.5!(v2)$);
			\coordinate (r3) at ($(v2)!0.5!(b2)$);

			\node[regionlabel , anchor=base] at (r1 |- lv1.base) {$\Delta_{1}^{\leq d}$};
			\node[regionlabel , anchor=base] at (r2 |- lv1.base) {$\Delta_{\emptyset}$};
			\node[regionlabel , anchor=base] at (r3 |- lv1.base) {$\Delta_{2}^{\leq d}$};

		\end{tikzpicture}
		\caption{
			Local tangent weights at the four residual fixed points of the cut of $T^{\ast}\CC\PP^{1}$.
			At each boundary vertex, the cutting weight occurs with multiplicity two and is represented by a pair of parallel arrows.
			}
		\figalttext{
			A horizontal segment with four marked points, labelled minus d, 0, m, and m plus d.
			At the same height, the three intervals between them are labelled Delta sub 1 with superscript at most d, Delta sub empty set, and Delta sub 2 with superscript at most d.
			At 0 and at m, one arrow points each way.
			At each end, a pair of parallel arrows points into the segment, labelled theta sub minus d equals plus 1 on the left and theta sub m plus d equals minus 1 on the right.
		}
		\label{fig:isotropy-TCP1}
	\end{figure}
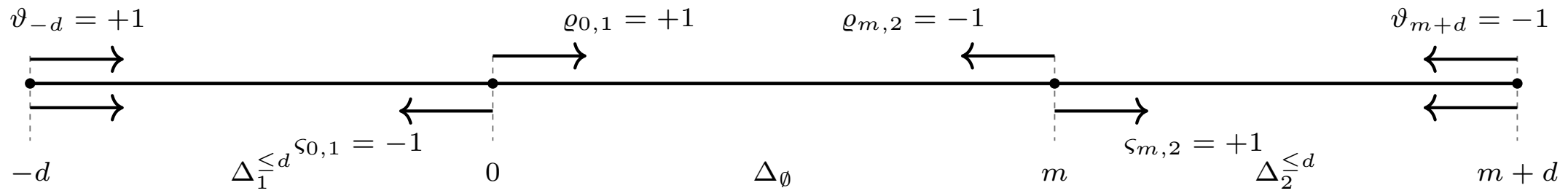

	Corollary~\ref{corollary:character-smooth-hypertoric-cut} gives
	\[
		\chi_{-d}(t) = \frac{t^{-d}}{( 1 - t )^{2}}, \quad
		\chi_{0}(t) =\frac{1}{( 1 - t )( 1 - t^{-1} )}, \quad
		\chi_{m}(t) = \frac{t^{m}}{( 1 - t )( 1 - t^{-1} )}, \quad
		\chi_{m + d}(t) = \frac{t^{m + d}}{( 1 - t^{-1} )^{2}}.
	\]
	Their sum factors as the product of the interval characters for $[ 0 , m + d ]$ and $[ -d , 0 ]$:
	\begin{equation}
		\label{eqn:TCP1-example-2-character}
		\chi_{d}(t) = \left( \sum_{r = 0}^{m + d} t^{r} \right)
		\left( \sum_{q = 0}^{d} t^{-q} \right).
	\end{equation}
	On the other hand, $\mcS(x) = \max\{ 0 , -x  \} + \max\{ 0 , x - m \}$, so Theorem~\ref{theorem:combinatorial-multiplicity} gives the multiplicity $d - \mcS(x) + 1$ at each integral $x \in [-d, m + d]$.
	This is also the coefficient of $t^{x}$ in \eqref{eqn:TCP1-example-2-character}, so both descriptions give
	\[
		\chi_{d}(t) = \sum_{x = -d}^{m + d} ( d - \mcS(x) + 1 ) t^{x}, \qquad \dim_{\CC}H^{0}(M_{\nu}^{\leq d} ; \mcL_{\nu}^{\leq d}) = ( m + d + 1 )( d + 1 ).
	\]

	For this example, the enlarged polygon of Remark~\ref{remark:enlarged-cut-moment-image} is $\Conv\{ ( -d , 0 ) , ( 0 , d ) , ( m , d ) , ( m + d , 0 ) \}$.
	Its vertical slice above an integral $x$ contains $d - \mcS(x) + 1$ lattice points.
	Thus, its ordinary lattice-point count equals the quantisation dimension, as it must for the smooth toric four-manifold $M_{\nu}^{\leq d}$ \cite{Hamilton2008}.
\end{example}

\begin{example}[$A_{2}$-Resolution, Continued]
	\label{example:A2-resolution-example-2}
	We continue \cite[Example~4.11]{BrownI} with $r = 2$, so that $\ell_{1}(x) = x$, $\ell_{2}(x) = m_{1} - x$, and $\ell_{3}(x) = m_{2} - x$, where $0 < m_{1} < m_{2}$.
	Set $e := d - ( m_{2} - m_{1})$.
	Since $\mcS(0) = \mcS(m_{1}) = 0$ and $\mcS(m_{2}) = m_{2} - m_{1}$, a positive regular integral level $d$ lies in the outer chamber if and only if $e > 0$.
	For $0 < d < m_{2} - m_{1}$, the cut does not reach the vertex $m_{2}$, and the calculation in Example~\ref{example:TCP1-example-2} applies verbatim with $m = m_{1}$; we therefore take $e > 0$.
	The polyptych then has components $[ -d , 0 ]$, $[ 0 , m_{1} ]$, $[ m_{1} , m_{2} ]$, and $[ m_{2} , v ]$, where $v := ( m_{1} + m_{2} + d ) / 2 = m_{2} + e/2$.
	At the right boundary vertex, the full labelled normal $u_{23} := u_{\{ 2 , 3 \}} = 2$ is non-primitive, so $\Gamma_{v} \cong \ZZ / 2\ZZ$ and $\vartheta_{v} = -1/2$; the other four fixed points have trivial local groups (Figure~\ref{fig:isotropy-A2}).

	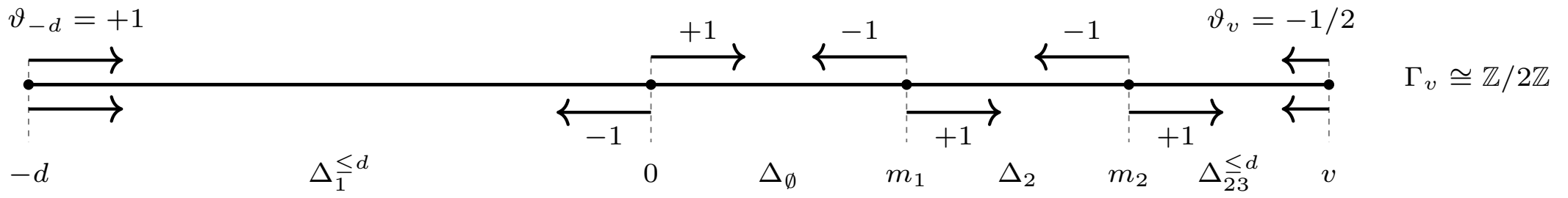
\begin{figure}[H]
		\centering
		\begin{tikzpicture}[x=1cm , y=1cm]

			\tikzset{
				polyboundary/.style={line width=0.9pt},
				weight/.style={->,line width=0.8pt},
				cutweight/.style={->,line width=0.8pt},
				fixedpoint/.style={font=\scriptsize},
				weightlabel/.style={font=\scriptsize},
				grouplabel/.style={font=\scriptsize},
				regionlabel/.style={font=\scriptsize},
			}

			\coordinate (b1)  at (-5.6 , 0);
			\coordinate (v1)  at ( 0.0 , 0);
			\coordinate (v2)  at ( 2.3 , 0);
			\coordinate (v3)  at ( 4.3 , 0);
			\coordinate (b23) at ( 6.1 , 0);

			\foreach \p in {v1,v2,v3}
			\draw[dashed, dash pattern=on 1.5pt off 1.5pt, line width=0.35pt, black!55]
			($(\p)+(0 , 0.25)$) -- ($(\p)+(0 , -0.52)$);

			\foreach \p in {b1,b23}
			\draw[dashed, dash pattern=on 1.5pt off 1.5pt, line width=0.35pt, black!55]
			($(\p)+(0 , 0.22)$) -- ($(\p)+(0 , -0.52)$);

			\draw[polyboundary] (b1) -- (b23);

			\foreach \p in {b1,v1,v2,v3,b23}
			\fill (\p) circle (1.35pt);

			\draw[cutweight]
			($(b1)+(0 , 0.22)$) -- ++(0.85 , 0);
			\draw[cutweight]
			($(b1)+(0 , -0.22)$) -- ++(0.85 , 0);

			\node[weightlabel , above=3pt]
			at ($(b1)+(0.43 , 0.22)$)
			{$\vartheta_{-d} = +1$};

			\draw[weight]
			($(v1)+(0 , 0.25)$) -- ++(0.85 , 0)
			node[midway , above , weightlabel] {$+1$};

			\draw[weight]
			($(v1)+(0 , -0.25)$) -- ++(-0.85 , 0)
			node[midway , below , weightlabel] {$-1$};

			\draw[weight]
			($(v2)+(0 , 0.25)$) -- ++(-0.85 , 0)
			node[midway , above , weightlabel] {$-1$};

			\draw[weight]
			($(v2)+(0 , -0.25)$) -- ++(0.85 , 0)
			node[midway , below , weightlabel] {$+1$};

			\draw[weight]
			($(v3)+(0 , 0.25)$) -- ++(-0.85 , 0)
			node[midway , above , weightlabel] {$-1$};

			\draw[weight]
			($(v3)+(0 , -0.25)$) -- ++(0.85 , 0)
			node[midway , below , weightlabel] {$+1$};

			\draw[cutweight]
			($(b23)+(0 , 0.22)$) -- ++(-0.425 , 0);
			\draw[cutweight]
			($(b23)+(0 , -0.22)$) -- ++(-0.425 , 0);

			\node[weightlabel , above=3pt]
			at ($(b23)+(-0.43 , 0.22)$)
			{$\vartheta_{v} = -1/2$};

			\node[grouplabel, anchor=west] at ($(b23)+(0.55 , 0.05)$)
			{$\Gamma_{v} \cong \ZZ/ 2\ZZ$};

			\node[fixedpoint , below=17pt] (lv1) at (v1)  {$0$};
			\node[fixedpoint , anchor=base] at (b1 |- lv1.base)  {$-d$};
			\node[fixedpoint , anchor=base] at (v2 |- lv1.base)  {$m_{1}$};
			\node[fixedpoint , anchor=base] at (v3 |- lv1.base)  {$m_{2}$};
			\node[fixedpoint , anchor=base] at (b23 |- lv1.base) {$v$};

			\coordinate (r1) at ($(b1)!0.5!(v1)$);
			\coordinate (r2) at ($(v1)!0.5!(v2)$);
			\coordinate (r3) at ($(v2)!0.5!(v3)$);
			\coordinate (r4) at ($(v3)!0.5!(b23)$);

			\node[regionlabel , anchor=base] at (r1 |- lv1.base) {$\Delta_{1}^{\leq d}$};
			\node[regionlabel , anchor=base] at (r2 |- lv1.base) {$\Delta_{\emptyset}$};
			\node[regionlabel , anchor=base] at (r3 |- lv1.base) {$\Delta_{2}$};
			\node[regionlabel , anchor=base] at (r4 |- lv1.base) {$\Delta_{23}^{\leq d}$};

		\end{tikzpicture}

		\caption{
				Local tangent weights at the five residual fixed points of the cut of the $A_{2}$-resolution.
				At the fixed point above the right boundary vertex, the full labelled normal $u_{23} = 2$ gives the tangent weight $\vartheta_{v} = -1/2$ with multiplicity two and the local orbifold group $\Gamma_{v} \cong \ZZ / 2\ZZ$, where $v = ( m_{1} + m_{2} + d ) / 2$.
			}
		\figalttext{
			A horizontal segment with five marked points, labelled minus d, 0, m sub 1, m sub 2, and v.
			At the same height, the four intervals between them are labelled Delta sub 1 with superscript at most d, Delta sub empty set, Delta sub 2, and Delta sub 2 3 with superscript at most d.
			At each of the three middle points, one arrow points each way.
			At minus d, a pair of parallel arrows points right, labelled theta sub minus d equals plus 1.
			At v, a pair of shorter parallel arrows points left, labelled theta sub v equals minus one half, beside the label Gamma sub v isomorphic to Z mod 2Z.
		}
		\label{fig:isotropy-A2}
	\end{figure}

	Write $d + m_{1} + m_{2} = 2a + b$ with $b \in \{ 0 , 1 \}$, so that $a = m_{2} + \lfloor e / 2 \rfloor$ and $b \equiv e \pmod{2}$.
	Theorem~\ref{theorem:character-hypertoric-orbifold} at the four smooth fixed points, together with Lemma~\ref{lemma:boundary-vertex-closed-form} at $v$ with $k = 2$, gives
	\[
		\chi_{d}(t) = \frac{t^{-d}}{( 1 - t )^{2}} + \frac{1 + t^{m_{1}} + t^{m_{2}}}{( 1 - t )( 1 - t^{-1} )} + t^{a}\, \frac{( b + 1 ) + ( 1 - b )\, t^{-1}}{( 1 - t^{-1} )^{2}}.
	\]
	The contributions at $-d$, $0$, and $m_{1}$ are those of Example~\ref{example:TCP1-example-2} with $m = m_{1}$, and comparing the remaining terms with \eqref{eqn:TCP1-example-2-character} gives
	\begin{equation}
		\label{eqn:A2-outer-character}
		\chi_{d}(t) = \left( \sum_{r = 0}^{m_{1} + d} t^{r} \right)
		\left( \sum_{q = 0}^{d} t^{-q} \right)
		- t^{m_{2}} \left( \sum_{k = 1}^{\lceil e / 2 \rceil} t^{k} \right)
		\left( \sum_{l = 0}^{\lfloor e / 2 \rfloor} t^{l} \right).
	\end{equation}
	The correction is supported on $( m_{2} , m_{1} + d ]$: beyond $m_{2}$, the slack function has slope $u_{23} = 2$, which also moves the right end of the support from $m_{1} + d$ to $\lfloor v \rfloor = a$.
	On $[-d, a]$, the coefficient of $t^{x}$ in \eqref{eqn:A2-outer-character} is $d - \mcS(x) + 1$, as Corollary~\ref{cor:weighted-lattice-point-character} requires.
	Evaluating at $t = 1$ gives
	\[
		\dim_{\CC} H^{0}( M_{\nu}^{\leq d} ; \mcL_{\nu}^{\leq d} ) =
		( m_{1} + d + 1 )( d + 1 ) - \lceil e/2 \rceil \big( \lfloor e/2 \rfloor + 1 \big).
	\]
	When $e$ is odd, $v$ is a half-integer and $\chi_{v}^{\fib}$ is non-trivial, but the averaged contribution still involves only integral powers of $t$.
\end{example}

\begin{example}[$T^{\ast}\CC\PP^{2}$, Continued]
	\label{example:TCP2-example-2}
	We continue \cite[Example~4.12]{BrownI} with $m , d \in \ZZ_{> 0}$.
	As noted there, $c_{\max} = 0$ and every cut is smooth, so its nine fixed points have trivial local groups; Figure~\ref{fig:isotropy-TCP2} shows their tangent weights.

	\begin{figure}[H]
		\centering
		\begin{tikzpicture}[scale=0.66]

			\tikzset{
				polyboundary/.style={draw=gray!100, line width=0.72pt},
				polyedge/.style={draw=gray!100, line width=0.42pt},
				weight/.style={->,draw=black , line width=0.82pt , shorten <=1.2pt},
				fixedpoint/.style={font=\scriptsize},
				weightlabel/.style={font=\scriptsize},
				regionlabel/.style={font=\scriptsize},
			}

			\def\mval{4.3}
			\def\dval{2.13}

			\def\alen{0.80}
			\def\diag{0.56}

			\def\alenTwo{0.68}
			\def\diagTwo{0.48}

			\coordinate (v12) at (0 , 0);
			\coordinate (v23) at (\mval , 0);
			\coordinate (v13) at (0,\mval);

			\coordinate (bL)  at (-\dval , 0);
			\coordinate (bB)  at (0 , -\dval);
			\coordinate (bBR) at ({\mval+\dval} , -\dval);
			\coordinate (bR)  at ({\mval+\dval},0);
			\coordinate (bT)  at (0,{\mval+\dval});
			\coordinate (bTL) at (-\dval,{\mval+\dval});

			\filldraw[
				fill=gray!70,
				draw=gray!100
			]
			(v12) -- (v23) -- (v13) -- cycle;

			\draw[polyboundary]
			(bL) -- (bB) -- (bBR) -- (bR)
			-- (bT) -- (bTL) -- cycle;

			\draw[polyedge] (bB)  -- (bT);
			\draw[polyedge] (bL)  -- (bR);
			\draw[polyedge] (bTL) -- (bBR);

			\node[regionlabel] at ({\mval/3} , {\mval/3}) {$\Delta_{\emptyset}$};

			\node[fixedpoint , above right=5pt]       at (v12)  {$( 0 , 0 )$};
			\node[fixedpoint , below=5pt]       at ($(v23)+(-0.50 , 0)$)  {$( m , 0 )$};
			\node[fixedpoint , left=5pt]       at ($(v13)+(0 , -0.35)$)  {$( 0 , m )$};
			\node[fixedpoint , left=5pt]       at (bL)  {$( -d , 0 )$};
			\node[fixedpoint , left=5pt]      at ($(bB)+(0 , -0.12)$)  {$( 0 , -d )$};
			\node[fixedpoint , right=5pt] at (bBR) {$( m + d , -d )$};
			\node[fixedpoint , right=5pt]      at (bR)  {$( m + d , 0 )$};
			\node[fixedpoint , right=5pt] at ($(bT)+(0 , 0.10)$)  {$( 0 , m + d )$};
			\node[fixedpoint , left=5pt] at (bTL) {$( -d , m + d )$};

			% Interior tangent weights
			\draw[weight] (v12) -- ++(\alen , 0);
			\draw[weight] (v12) -- ++(0,\alen);
			\draw[weight] (v12) -- ++(-\alen , 0);
			\draw[weight] (v12) -- ++(0 , -\alen);

			\draw[weight] (v23) -- ++(\alen , 0);
			\draw[weight] (v23) -- ++(-\alen , 0);
			\draw[weight] (v23) -- ++(-\diag,\diag);
			\draw[weight] (v23) -- ++(\diag , -\diag);

			\draw[weight] (v13) -- ++(0,\alen);
			\draw[weight] (v13) -- ++(0 , -\alen);
			\draw[weight] (v13) -- ++(-\diag,\diag);
			\draw[weight] (v13) -- ++(\diag , -\diag);

			% Boundary tangent weights
			\draw[weight] (bL) -- ++(\alen , 0);
			\draw[weight] (bL) -- ++(\alenTwo , 0);
			\draw[weight] (bL) -- ++(0,\alen);
			\draw[weight] (bL) -- ++(\diag , -\diag);

			\draw[weight] (bB) -- ++(0,\alen);
			\draw[weight] (bB) -- ++(0,\alenTwo);
			\draw[weight] (bB) -- ++(\alen , 0);
			\draw[weight] (bB) -- ++(-\diag,\diag);

			\draw[weight] (bR) -- ++(-\alen , 0);
			\draw[weight] (bR) -- ++(-\alenTwo , 0);
			\draw[weight] (bR) -- ++(-\diag,\diag);
			\draw[weight] (bR) -- ++(0 , -\alen);

			\draw[weight] (bBR) -- ++(-\diag,\diag);
			\draw[weight] (bBR) -- ++(-\diagTwo,\diagTwo);
			\draw[weight] (bBR) -- ++(0,\alen);
			\draw[weight] (bBR) -- ++(-\alen , 0);

			\draw[weight] (bT) -- ++(0 , -\alen);
			\draw[weight] (bT) -- ++(0 , -\alenTwo);
			\draw[weight] (bT) -- ++(\diag , -\diag);
			\draw[weight] (bT) -- ++(-\alen , 0);

			\draw[weight] (bTL) -- ++(\diag , -\diag);
			\draw[weight] (bTL) -- ++(\diagTwo , -\diagTwo);
			\draw[weight] (bTL) -- ++(\alen , 0);
			\draw[weight] (bTL) -- ++(0 , -\alen);

			% Legend
			\begin{scope}[shift={(0.13 , 0.13)}]
				\coordinate (L) at (7.55 , 4.35);

				\node[font=\small , anchor=west] at (7.02 , 5.68) {weight directions};

				\draw[weight] (L) -- ++(0.91 , 0);
				\node[weightlabel , anchor=west] at ($(L)+(1.07 , 0)$) {$t_1$};

				\draw[weight] (L) -- ++(0 , 0.91);
				\node[weightlabel , anchor=west] at ($(L)+(0.21 , 0.80)$) {$t_2$};

				\coordinate (D1) at (7.58 , 3.91);
				\coordinate (D2) at (8.11 , 1.95);

				\draw[weight] (D1) -- ++(0.60 , -0.60);
				\node[weightlabel , anchor=west]
				at (8.59 , 3.47) {$t_1t_2^{-1}$};

				\draw[weight] (D2) -- ++(-0.60 , 0.60);
				\node[weightlabel , anchor=west]
				at (8.59 , 2.39) {$t_1^{-1}t_2$};
			\end{scope}

			\foreach \p in {v12,v23,v13,bL,bB,bBR,bR,bT,bTL}
			\fill[black] (\p) circle (1.35pt);
		\end{tikzpicture}

		\caption{
			Local tangent weights at the nine residual fixed points of the cut of $T^{\ast}\CC\PP^{2}$.
			The arrows indicate the geometric tangent-weight directions.
			At each boundary vertex, the cutting weight occurs with multiplicity two and is represented by two closely spaced arrowheads on the same ray.
		}
		\figalttext{
			A hexagon divided by a horizontal, a vertical, and a diagonal line into seven polygons, with the central triangle shaded, labelled Delta sub empty set, and its vertices labelled (0, 0), (m, 0), and (0, m).
			Short arrows at each of the nine marked vertices show the tangent-weight directions; at the six outer vertices, a doubled arrowhead marks the cutting weight.
			A legend shows the directions t sub 1, t sub 2, t sub 1 t sub 2 inverse, and t sub 1 inverse t sub 2.
		}
		\label{fig:isotropy-TCP2}
	\end{figure}
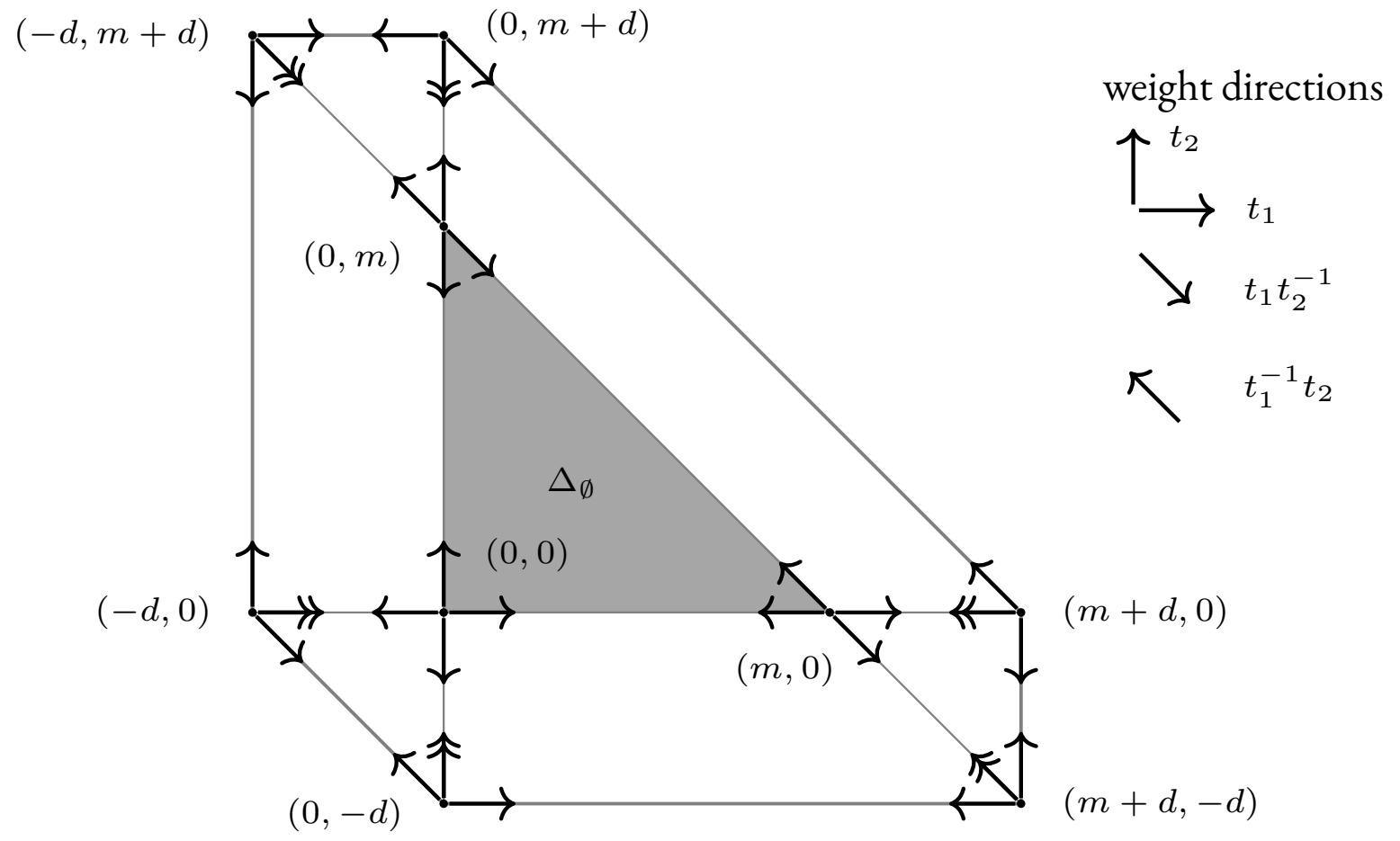
	Summing the nine contributions of Corollary~\ref{corollary:character-smooth-hypertoric-cut} gives the product of the lattice-point characters of the base triangle $\Conv\{ ( 0 , 0 ) , ( m + d , 0 ) , ( 0 , m + d ) \}$ and the fibre triangle $\Conv\{ ( 0 , 0 ) , ( -d , 0 ) , ( 0 , -d ) \}$:
	\[
		\chi_{d}(t_{1} , t_{2}) =
		\left( \sum\nolimits_{\substack{r , q \geq 0 \\ r + q \leq m + d}} t_{1}^{r}t_{2}^{q} \right)
		\left( \sum\nolimits_{\substack{k , l \geq 0 \\ k + l \leq d}} t_{1}^{-k}t_{2}^{-l} \right),
		\qquad
		\dim_{\CC} H^{0}( M_{\nu}^{\leq d} ; \mcL_{\nu}^{\leq d} ) =
		\binom{m + d + 2}{2} \binom{d + 2}{2}.
	\]
	The coefficient of $\chi_{d}$ at $x$ is $\binom{d - \mcS(x) + 2}{2}$, the multiplicity of Theorem~\ref{theorem:combinatorial-multiplicity} on the hexagonal support $\{ \mcS \leq d \}$.
\end{example}

\begin{example}[$T^{\ast}( \CC\PP^{1} \times \CC\PP^{1} )$, Continued]
	\label{example:TCP1xCP1-example-2}
	We continue \cite[Example~4.13]{BrownI} with $m_{1} , m_{2} , d \in \ZZ_{> 0}$.
	Again, every cut is smooth; Figure~\ref{fig:isotropy-TCP1xCP1} shows the tangent weights at its twelve fixed points.

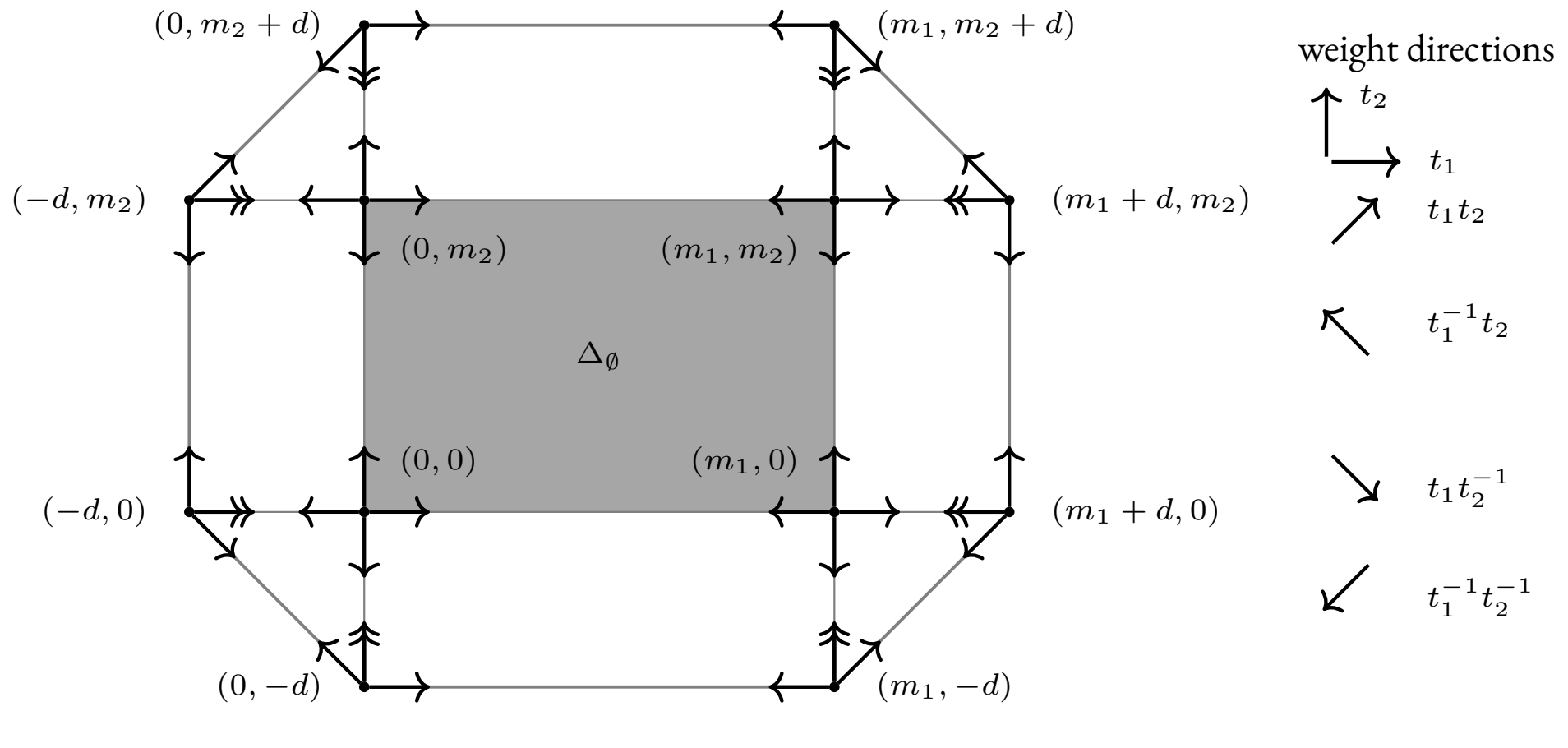
\begin{figure}[H]
    \centering
    \begin{tikzpicture}[scale=0.88]
        \tikzset{
            polyboundary/.style={draw=gray!100 , line width=0.72pt},
            polyedge/.style={draw=gray!100 , line width=0.42pt},
            weight/.style={->,draw=black , line width=0.82pt , shorten <=1.2pt},
            fixedpoint/.style={font=\scriptsize},
            weightlabel/.style={font=\scriptsize},
            regionlabel/.style={font=\scriptsize},
        }

        \def\moneval{4.3}
        \def\mtwoval{2.85}
        \def\dval{1.6}

        \def\alen{0.60}
        \def\alenTwo{0.51}
        \def\diag{0.42}
        \def\diagTwo{0.36}

        \coordinate (v00) at (0 , 0);
        \coordinate (v10) at (\moneval , 0);
        \coordinate (v11) at (\moneval , \mtwoval);
        \coordinate (v01) at (0 , \mtwoval);

        \coordinate (bL1) at (-\dval , 0);
        \coordinate (bB1) at (0 , -\dval);
        \coordinate (bB2) at (\moneval , -\dval);
        \coordinate (bR1) at ({\moneval+\dval} , 0);
        \coordinate (bR2) at ({\moneval+\dval} , \mtwoval);
        \coordinate (bT2) at (\moneval , {\mtwoval+\dval});
        \coordinate (bT1) at (0 , {\mtwoval+\dval});
        \coordinate (bL2) at (-\dval , \mtwoval);

        \filldraw[
        	fill=gray!70,
        	draw=gray!100
        ]
        (v00) -- (v10) -- (v11) -- (v01) -- cycle;
        \draw[polyboundary]
            (bL1) -- (bB1) -- (bB2) -- (bR1) -- (bR2)
            -- (bT2) -- (bT1) -- (bL2) -- cycle;
        \draw[polyedge] (bB1) -- (bT1);
        \draw[polyedge] (bB2) -- (bT2);
        \draw[polyedge] (bL1) -- (bR1);
        \draw[polyedge] (bL2) -- (bR2);

        \node[regionlabel] at ({\moneval/2} , {\mtwoval/2}) {$\Delta_{\emptyset}$};

        \node[fixedpoint , anchor=south west] at ($(v00)+(0.18 , 0.18)$) {$(0 , 0)$};
        \node[fixedpoint , anchor=south east] at ($(v10)+(-0.18 , 0.18)$) {$(m_{1} , 0)$};
        \node[fixedpoint , anchor=north east] at ($(v11)+(-0.18 , -0.18)$) {$(m_{1} , m_{2})$};
        \node[fixedpoint , anchor=north west] at ($(v01)+(0.18 , -0.18)$) {$(0 , m_{2})$};

        \node[fixedpoint , left=6pt]  at (bL1) {$(-d , 0)$};
        \node[fixedpoint , left=6pt] at (bB1) {$(0 , -d)$};
        \node[fixedpoint , right=6pt] at (bB2) {$(m_{1} , -d)$};
        \node[fixedpoint , right=6pt] at (bR1) {$(m_{1} + d , 0)$};
        \node[fixedpoint , right=6pt] at (bR2) {$(m_{1} + d , m_{2})$};
        \node[fixedpoint , right=6pt] at (bT2) {$(m_{1} , m_{2} + d)$};
        \node[fixedpoint , left=6pt] at (bT1) {$(0 , m_{2} + d)$};
        \node[fixedpoint , left=6pt]  at (bL2) {$(-d , m_{2})$};

        % Interior fixed points
        \foreach \p in {v00,v10,v11,v01}{
            \draw[weight] (\p) -- ++(\alen , 0);
            \draw[weight] (\p) -- ++(-\alen , 0);
            \draw[weight] (\p) -- ++(0 , \alen);
            \draw[weight] (\p) -- ++(0 , -\alen);
        }

        % Boundary fixed points
        \draw[weight] (bL1) -- ++(\alen , 0);
        \draw[weight] (bL1) -- ++(\alenTwo , 0);
        \draw[weight] (bL1) -- ++(0 , \alen);
        \draw[weight] (bL1) -- ++(\diag , -\diag);

        \draw[weight] (bB1) -- ++(0 , \alen);
        \draw[weight] (bB1) -- ++(0 , \alenTwo);
        \draw[weight] (bB1) -- ++(\alen , 0);
        \draw[weight] (bB1) -- ++(-\diag , \diag);

        \draw[weight] (bB2) -- ++(0 , \alen);
        \draw[weight] (bB2) -- ++(0 , \alenTwo);
        \draw[weight] (bB2) -- ++(-\alen , 0);
        \draw[weight] (bB2) -- ++(\diag , \diag);

        \draw[weight] (bR1) -- ++(-\alen , 0);
        \draw[weight] (bR1) -- ++(-\alenTwo , 0);
        \draw[weight] (bR1) -- ++(0 , \alen);
        \draw[weight] (bR1) -- ++(-\diag , -\diag);

        \draw[weight] (bR2) -- ++(-\alen , 0);
        \draw[weight] (bR2) -- ++(-\alenTwo , 0);
        \draw[weight] (bR2) -- ++(0 , -\alen);
        \draw[weight] (bR2) -- ++(-\diag , \diag);

        \draw[weight] (bT2) -- ++(0 , -\alen);
        \draw[weight] (bT2) -- ++(0 , -\alenTwo);
        \draw[weight] (bT2) -- ++(-\alen , 0);
        \draw[weight] (bT2) -- ++(\diag , -\diag);

        \draw[weight] (bT1) -- ++(0 , -\alen);
        \draw[weight] (bT1) -- ++(0 , -\alenTwo);
        \draw[weight] (bT1) -- ++(\alen , 0);
        \draw[weight] (bT1) -- ++(-\diag , -\diag);

        \draw[weight] (bL2) -- ++(\alen , 0);
        \draw[weight] (bL2) -- ++(\alenTwo , 0);
        \draw[weight] (bL2) -- ++(0 , -\alen);
        \draw[weight] (bL2) -- ++(\diag , \diag);

        % Legend
        \begin{scope}[shift={(-0.25 , -0.25)}]
            \node[font=\small , anchor=west] at (8.65 , 4.45) {weight directions};
            \coordinate (L) at (9.05 , 3.45);
            \draw[weight] (L) -- ++(0.68 , 0);
            \node[weightlabel , anchor=west] at ($(L)+(0.80 , 0)$) {$t_{1}$};
            \draw[weight] (L) -- ++(0 , 0.68);
            \node[weightlabel , anchor=west] at ($(L)+(0.16 , 0.60)$) {$t_{2}$};

            \coordinate (D1) at (9.07 , 2.67);
            \coordinate (D2) at (9.47 , 1.65);
            \coordinate (D3) at (9.07 , 0.80);
            \coordinate (D4) at (9.47 , -0.20);

            \draw[weight] (D1) -- ++(0.45 , 0.45);
            \node[weightlabel , anchor=west] at (9.83 , 3.00) {$t_{1}t_{2}$};
            \draw[weight] (D2) -- ++(-0.45 , 0.45);
            \node[weightlabel , anchor=west] at (9.83 , 1.98) {$t_1^{-1}t_{2}$};
            \draw[weight] (D3) -- ++(0.45 , -0.45);
            \node[weightlabel , anchor=west] at (9.83 , 0.47) {$t_{1}t_2^{-1}$};
            \draw[weight] (D4) -- ++(-0.45 , -0.45);
            \node[weightlabel , anchor=west] at (9.83 , -0.53) {$t_1^{-1}t_2^{-1}$};
        \end{scope}

        \foreach \p in {v00,v10,v11,v01,bL1,bB1,bB2,bR1,bR2,bT2,bT1,bL2}
            \fill[black] (\p) circle (1.35pt);
    \end{tikzpicture}

    \caption{
        Local tangent weights at the twelve residual fixed points of the cut of $T^{\ast}(\CC\PP^{1} \times \CC\PP^{1})$.
        The arrows indicate the geometric tangent-weight directions.
        At each boundary vertex, the cutting weight occurs with multiplicity two and is represented by two closely spaced arrowheads on the same ray.
    }
    \figalttext{
    	An octagon divided by two horizontal and two vertical lines into nine polygons, with the central rectangle shaded, labelled Delta sub empty set, and its vertices labelled (0, 0), (m sub 1, 0), (m sub 1, m sub 2), and (0, m sub 2).
    	Short arrows at each of the twelve marked vertices show the tangent-weight directions; at the eight outer vertices, a doubled arrowhead marks the cutting weight.
    	A legend shows the weight directions.
    }
    \label{fig:isotropy-TCP1xCP1}
\end{figure}
	Grouping the twelve contributions of Corollary~\ref{corollary:character-smooth-hypertoric-cut} by the four vertices of the base rectangle and summing gives
	\[
		\chi_{d}(t_{1} , t_{2}) =
		\sum\nolimits_{\substack{k , l \geq 0 \\ k + l \leq d}} t_{1}^{-k}t_{2}^{-l}
		\left( \sum_{r = 0}^{m_{1} + 2k} t_{1}^{r} \right)
		\left( \sum_{q = 0}^{m_{2} + 2l} t_{2}^{q} \right),
	\]
	whose coefficient at $x$ is again $\binom{d - \mcS(x) + 2}{2}$, now on the octagonal support.
	Evaluating at $( 1 , 1 )$ gives
	\[
		\begin{aligned}
			\dim_{\CC} H^{0}( M_{\nu}^{\leq d} ; \mcL_{\nu}^{\leq d} )
			 & = \sum_{k = 0}^{d} \sum_{l = 0}^{d - k} ( m_{1} + 2k + 1 )( m_{2} + 2l + 1 ) \\
			 & = \frac{ ( d + 1 )( d + 2 ) \big( d^{2} + 2d( m_{1} + m_{2} ) + 3d + 3m_{1}m_{2} + 3m_{1} + 3m_{2} + 3 \big) }{6}.
		\end{aligned}
	\]
\end{example}

%% file: sections/7-conclusion.tex
% !TeX root = ../main.tex
\section{Conclusion}
\label{section:section-7-conclusion}

At every positive regular integral level $d$, Theorem~\ref{theorem:cut-cohomology-every-level} identifies the sections of $\mcL_{\nu}^{\leq d}$ with the filtration piece $F_{d}H^{0}(M_{\nu} ; \mcL_{\nu})$ and shows that the higher cohomology groups of $\mcL_{\nu}^{\leq d}$ vanish.
The Kawasaki formula of Theorem~\ref{theorem:character-hypertoric-orbifold} therefore computes the character of the space of sections, which Corollary~\ref{cor:weighted-lattice-point-character} also gives as a weighted lattice-point sum over the moment polyptych.
Since Theorem~\ref{theorem:cut-cohomology-every-level} holds in every cut chamber, its identifications embed the section space at a positive regular integral level $d^{\prime} < d$ in that at $d$ as a $T^{n} \times \CC^{\ast}$-representation, with quotient the sum of the weight spaces $H^{0}(M_{\nu} ; \mcL_{\nu})_{w}$ for $d^{\prime} < w \leq d$.
This holds also when a critical value lies between $d^{\prime}$ and $d$, although the moment polyptych changes there, and with it the fixed points entering Theorem~\ref{theorem:character-hypertoric-orbifold}.

Two questions follow from this comparison.
First, can the multiplicities be realised geometrically in higher rank?
By Remark~\ref{remark:reduced-spaces}, the multiplicity of a residual weight $x$ with $\mcS(x) < d$ is the dimension of the space of sections on the coarse quotient $M_{\nu}^{\alg , \leq d} \sslash_{x} T_{\CC}^{n} \cong \CC\PP^{n}$.
For $n \geq 2$, the auxiliary circle acts on this quotient only through a single circle subgroup of its $n$-torus.
Do the algebraic toric structures of these quotients come from a completely integrable system on the cut whose Hamiltonians include the components of the moment map for $T^{n} \times U_{1}$, or from a toric degeneration of the cut?
Such a system would turn the weighted lattice-point sum of Corollary~\ref{cor:weighted-lattice-point-character} into an ordinary lattice-point count on the cut itself, as Remark~\ref{remark:enlarged-cut-moment-image} does for $n = 1$ in the outer chamber.
For the moduli space of flat $\SU(2)$-connections discussed in Section~\ref{section:section-1-introduction}, the corresponding lattice-point count is that of a toric compactification with a different complex structure, whereas on the cut, the count and the section spaces would refer to the same complex structure.

Second, which polarised semiprojective varieties admit analogous cohomology comparisons?
As in Remark~\ref{remark:algebraic-cut-generality}, the algebraic cut and its filtration can be defined more generally for a compatible linearisation whose section ring is finitely generated.
In the hypertoric case, Theorem~\ref{theorem:cut-cohomology-every-level} uses the description of the unstable strata by fibre-rotation cells and the fact that the higher cohomology groups on $M_{\nu}$ vanish; the explicit multiplicities additionally use the hypertoric coordinate ring.

A motivating case is the moduli space of rank-two Higgs bundles with fixed determinant of degree one and trace-free Higgs field on a curve of genus $g \geq 2$.
Halpern-Leistner's proof of the equivariant Verlinde formula uses a parallel comparison: for line bundles of positive level, a quantisation theorem identifies cohomology on the stack of Higgs bundles with cohomology on its semistable locus, and a vanishing theorem then identifies the graded index with the graded dimension of the space of sections \cite[Proposition~3.4 and Theorem~3.11]{HalpernLeistner2016}.
Hausel's compactification of this moduli space is the cut above the largest critical value of the Higgs-field norm \cite[Definition~7.1]{Hausel1998}, and its singular points are $\ZZ / 2\ZZ$-points on the divisor $Z$ at infinity \cite[Theorems~6.2 and~7.3]{Hausel1998}, so the theorem of Wu recalled in Remark~\ref{remark:quantisation-attribution} does not cover it.
For this outer cut, Hausel's description of the compactification gives both parts of the comparison.
Hartogs's theorem reduces the comparison of sections to the cotangent bundle of the moduli space $N$ of stable bundles, since the compactification is isomorphic, outside a subset of codimension at least two, to the first-chamber cut $\PP( T^{\ast}N \oplus \mcO_{N} )$ \cite[Theorem~7.4]{Hausel1998}.
The anticanonical line orbibundle of the compactification is $\mcO((3g - 2)Z)$, and $\mcO(2Z)$ is the pullback of the hyperplane bundle under the extended Hitchin map \cite[Corollaries~7.7 and~7.9]{Hausel1998}.
The anticanonical line orbibundle is therefore semipositive, and the higher cohomology groups of a positive pre-quantum line orbibundle vanish by Baily's vanishing theorem for V-manifolds \cite[\S4, (1)]{Baily1957}.
Once the line bundles and moment maps have been matched with those in the fixed-determinant version of the equivariant Verlinde formula, the equivariant index of the compactification at level $d$ is then the part of fibre-rotation weight at most $d$ of the equivariant Verlinde series, which an orbifold localisation formula with positive-dimensional fixed components computes.
Does the comparison of Theorem~\ref{theorem:cut-cohomology-every-level}, including the vanishing of higher cohomology, also hold for Hausel's cuts in the inner chambers?

\section*{Funding}

This work was supported by a School of Mathematics Studentship from The University of Edinburgh.

\section*{Acknowledgements}
\label{section:acknowledgements}

This article grew out of the author's doctoral research at The University of Edinburgh.
The author is especially grateful to Johan Martens for suggesting the problem that became the starting point for this work, and for his continued guidance during the doctoral research in which the localisation formula of Section~\ref{section:section-4-equivariant-quantisation-cut-spaces} was obtained.
The author also thanks Andrew Dancer for many helpful conversations and for hosting a visit to the Mathematical Institute at The University of Oxford, whose hospitality and support are gratefully acknowledged.

%% file: bibliography.bib
@article{AbramovichOlssonVistoli2008,
		author  = {Abramovich, Dan and Olsson, Martin and Vistoli, Angelo},
		title   = {Tame stacks in positive characteristic},
		journal = {Ann. Inst. Fourier (Grenoble)},
		year    = {2008},
		volume  = {58},
		number  = {4},
		pages   = {1057--1091},
		doi     = {10.5802/aif.2378},
	}

@article{Alper2013,
		author  = {Alper, Jarod},
		title   = {Good moduli spaces for {A}rtin stacks},
		journal = {Ann. Inst. Fourier (Grenoble)},
		year    = {2013},
		volume  = {63},
		number  = {6},
		pages   = {2349--2402},
		doi     = {10.5802/aif.2833},
	}

@misc{AndersenGukovPei2016,
		author        = {Andersen, J{\o}rgen Ellegaard and Gukov, Sergei and Pei, Du},
		title         = {The {V}erlinde formula for {H}iggs bundles},
		year          = {2016},
		eprint        = {1608.01761},
		archivePrefix = {arXiv},
		primaryClass  = {math.AG},
	}

@article{Baily1957,
		author  = {Baily, Walter L.},
		title   = {On the imbedding of {V}-manifolds in projective space},
		journal = {Amer. J. Math.},
		year    = {1957},
		volume  = {79},
		number  = {2},
		pages   = {403--430},
		doi     = {10.2307/2372689},
	}

@misc{BellamyCrawSchedler2024,
    author        = {Bellamy, Gwyn and Craw, Alastair and Schedler, Travis},
    title         = {The semi-invariant ring as the {C}ox ring of a {GIT} quotient},
    year          = {2024},
    eprint        = {2404.12225},
    archivePrefix = {arXiv},
    primaryClass  = {math.AG},
}

@article{BellamyKuwabara2012,
    author  = {Bellamy, Gwyn and Kuwabara, Toshiro},
    title   = {On deformation quantizations of hypertoric varieties},
    journal = {Pacific J. Math.},
    year    = {2012},
    volume  = {260},
    number  = {1},
    pages   = {89--127},
    doi     = {10.2140/pjm.2012.260.89},
}

@article{BialynickiBirula1973,
		author  = {Bia{\l}ynicki-Birula, Andrzej},
		title   = {Some theorems on actions of algebraic groups},
		journal = {Ann. of Math. (2)},
		year    = {1973},
		volume  = {98},
		pages   = {480--497},
		doi     = {10.2307/1970915},
	}

@article{BielawskiDancer2000,
		author  = {Bielawski, Roger and Dancer, Andrew S.},
		title   = {The geometry and topology of toric hyperk{\"a}hler manifolds},
		journal = {Comm. Anal. Geom.},
		year    = {2000},
		volume  = {8},
		number  = {4},
		pages   = {727--760},
		doi     = {10.4310/CAG.2000.v8.n4.a2},
	}

@article{BravermanFinkelbergNakajima2018,
		author  = {Braverman, Alexander and Finkelberg, Michael and Nakajima, Hiraku},
		title   = {Towards a mathematical definition of {C}oulomb branches of $3$-dimensional $\mathcal{N}=4$ gauge theories, {II}},
		journal = {Adv. Theor. Math. Phys.},
		year    = {2018},
		volume  = {22},
		number  = {5},
		pages   = {1071--1147},
		doi     = {10.4310/ATMP.2018.v22.n5.a1},
	}

@incollection{BrionProcesi1990,
		author    = {Brion, Michel and Procesi, Claudio},
		title     = {Action d'un tore dans une vari{\'e}t{\'e} projective},
		booktitle = {Operator Algebras, Unitary Representations, Enveloping Algebras, and Invariant Theory},
		editor    = {Connes, Alain and Duflo, Michel and Joseph, Anthony and Rentschler, Rudolf},
		series    = {Progress in Mathematics},
		volume    = {92},
		pages     = {509--539},
		publisher = {Birkh{\"a}user Boston},
		address   = {Boston, MA},
		year      = {1990},
	}

@unpublished{BrownI,
	author = {Brown, Benjamin C. W.},
	title  = {Symplectic Cuts and Moment Polyptychs of Hypertoric Varieties},
	year   = {2026},
	eprint        = {2610.37724},
	archivePrefix = {arXiv},
	primaryClass  = {math.SG}
}

@phdthesis{CannasDaSilva1996,
    author = {Cannas da Silva, Ana},
    title  = {Multiplicity formulas for orbifolds},
    school = {Massachusetts Institute of Technology},
    year   = {1996},
    url    = {http://hdl.handle.net/1721.1/38409},
}

@article{CremonesiHananyMekareeyaZaffaroni2014,
		author  = {Cremonesi, Stefano and Hanany, Amihay and Mekareeya, Noppadol and Zaffaroni, Alberto},
		title   = {Coulomb branch {H}ilbert series and {H}all--{L}ittlewood polynomials},
		journal = {J. High Energy Phys.},
		year    = {2014},
		volume  = {2014},
		number  = {9},
		pages   = {178},
		doi     = {10.1007/JHEP09(2014)178},
	}

@article{CremonesiHananyZaffaroni2014,
		author  = {Cremonesi, Stefano and Hanany, Amihay and Zaffaroni, Alberto},
		title   = {Monopole operators and {H}ilbert series of {C}oulomb branches of $3d$ $\mathcal{N}=4$ gauge theories},
		journal = {J. High Energy Phys.},
		year    = {2014},
		volume  = {2014},
		number  = {1},
		pages   = {005},
		doi     = {10.1007/JHEP01(2014)005},
	}

@article{DGMW1995,
		shorthand = {DGMW},
		author     = {Duistermaat, Johannes J. and Guillemin, Victor W. and Meinrenken, Eckhard and Wu, Siye},
		title      = {Symplectic reduction and {R}iemann--{R}och for circle actions},
		journal    = {Math. Res. Lett.},
		year       = {1995},
		volume     = {2},
		number     = {3},
		pages      = {259--266},
		doi        = {10.4310/MRL.1995.v2.n3.a3},
	}

@article{DolgachevHu1998,
		author  = {Dolgachev, Igor V. and Hu, Yi},
		title   = {Variation of geometric invariant theory quotients},
		journal = {Inst. Hautes {\'E}tudes Sci. Publ. Math.},
		year    = {1998},
		volume  = {87},
		pages   = {5--51},
		doi     = {10.1007/BF02698859},
	}

@article{EdidinGraham1998,
		author  = {Edidin, Dan and Graham, William},
		title   = {Algebraic cuts},
		journal = {Proc. Amer. Math. Soc.},
		year    = {1998},
		volume  = {126},
		number  = {3},
		pages   = {677--685},
		doi     = {10.1090/S0002-9939-98-04439-6},
	}

@article{Fujino2007,
    author  = {Fujino, Osamu},
    title   = {Multiplication maps and vanishing theorems for toric varieties},
    journal = {Math. Z.},
    year    = {2007},
    volume  = {257},
    number  = {3},
    pages   = {631--641},
    doi     = {10.1007/s00209-007-0140-5},
}

@article{GotoWatanabe1978,
    author  = {Goto, Shiro and Watanabe, Kei-ichi},
    title   = {On graded rings, I},
    journal = {J. Math. Soc. Japan},
    year    = {1978},
    volume  = {30},
    number  = {2},
    pages   = {179--213},
    doi     = {10.2969/jmsj/03020179},
}

@article{Grothendieck1961,
    author  = {Grothendieck, Alexander},
    title   = {{\'E}l{\'e}ments de g{\'e}om{\'e}trie alg{\'e}brique. {II}. {\'E}tude globale {\'e}l{\'e}mentaire de quelques classes de morphismes},
    journal = {Inst. Hautes {\'E}tudes Sci. Publ. Math.},
    year    = {1961},
    volume  = {8},
    pages   = {5--222},
    doi     = {10.1007/BF02699291},
}

@book{GuilleminGinzburgKarshon2002,
		author    = {Ginzburg, Viktor L. and Guillemin, Victor W. and Karshon, Yael},
		title     = {Moment maps, cobordisms, and {H}amiltonian group actions},
		series    = {Mathematical Surveys and Monographs},
		volume    = {98},
		note      = {Appendix J by Maxim Braverman},
		publisher = {American Mathematical Society},
		address   = {Providence, RI},
		year      = {2002},
		isbn      = {0-8218-0502-9},
		doi       = {10.1090/surv/098},
	}

@article{GukovPei2017,
		author  = {Gukov, Sergei and Pei, Du},
		title   = {Equivariant {V}erlinde formula from fivebranes and vortices},
		journal = {Comm. Math. Phys.},
		year    = {2017},
		volume  = {355},
		number  = {1},
		pages   = {1--50},
		doi     = {10.1007/s00220-017-2931-9},
	}

@article{HalpernLeistner2015,
		author  = {Halpern-Leistner, Daniel},
		title   = {The derived category of a {GIT} quotient},
		journal = {J. Amer. Math. Soc.},
		year    = {2015},
		volume  = {28},
		number  = {3},
		pages   = {871--912},
		doi     = {10.1090/S0894-0347-2014-00815-8},
	}

@misc{HalpernLeistner2016,
		author        = {Halpern-Leistner, Daniel},
		title         = {The equivariant {V}erlinde formula on the moduli of {H}iggs bundles},
		year          = {2016},
		eprint        = {1608.01754},
		archivePrefix = {arXiv},
		primaryClass  = {math.AG},
		note = {Appendix by Constantin Teleman},
	}

@incollection{Hamilton2008,
		author    = {Hamilton, Mark D.},
		title     = {The quantization of a toric manifold is given by the integer lattice points in the moment polytope},
		booktitle = {Toric Topology},
        editor    = {Harada, Megumi and Karshon, Yael and Masuda, Mikiya and Panov, Taras},
		series    = {Contemporary Mathematics},
		volume    = {460},
		pages     = {131--140},
		publisher = {American Mathematical Society},
		address   = {Providence, RI},
		year      = {2008},
		doi       = {10.1090/conm/460/09015},
	}

@article{HaradaProudfoot2004,
		author  = {Harada, Megumi and Proudfoot, Nicholas J.},
		title   = {Properties of the residual circle action on a hypertoric variety},
		journal = {Pacific J. Math.},
		year    = {2004},
		volume  = {214},
		number  = {2},
		pages   = {263--284},
		doi     = {10.2140/pjm.2004.214.263},
	}

@book{Hartshorne1977,
		author    = {Hartshorne, Robin},
		title     = {Algebraic geometry},
		series    = {Graduate Texts in Mathematics},
		volume    = {52},
		publisher = {Springer-Verlag},
		address   = {New York},
		year      = {1977},
		isbn      = {0-387-90244-9},
		doi       = {10.1007/978-1-4757-3849-0},
	}

@article{Hausel1998,
		author  = {Hausel, Tam{\'a}s},
		title   = {Compactification of moduli of {H}iggs bundles},
		journal = {J. Reine Angew. Math.},
		year    = {1998},
		volume  = {503},
		pages   = {169--192},
		doi     = {10.1515/crll.1998.096},
	}

@article{HauselRodriguezVillegas2015,
		author  = {Hausel, Tam{\'a}s and Rodr{\'i}guez Villegas, Fernando},
		title   = {Cohomology of large semiprojective hyperk{\"a}hler varieties},
		journal = {Ast{\'e}risque},
		year    = {2015},
		volume  = {370},
		pages   = {113--156},
		doi     = {10.24033/ast.967},
	}

@article{HauselSturmfels2002,
		author  = {Hausel, Tam{\'a}s and Sturmfels, Bernd},
		title   = {Toric hyperk{\"a}hler varieties},
		journal = {Doc. Math.},
		year    = {2002},
		volume  = {7},
		pages   = {495--534},
		doi     = {10.4171/DM/130},
	}

@article{HitchinKarlhedeLindstromRocek1987,
		author  = {Hitchin, Nigel J. and Karlhede, Anders and Lindstr{\"o}m, Ulf and Ro\v{c}ek, Martin},
		title   = {Hyper-{K}{\"a}hler metrics and supersymmetry},
		journal = {Comm. Math. Phys.},
		year    = {1987},
		volume  = {108},
		number  = {4},
		pages   = {535--589},
		doi = {10.1007/BF01214418},
	}

@article{JeffreyWeitsman1992,
		author  = {Jeffrey, Lisa C. and Weitsman, Jonathan},
		title   = {Bohr-{S}ommerfeld orbits in the moduli space of flat connections and the {V}erlinde dimension formula},
		journal = {Comm. Math. Phys.},
		year    = {1992},
		volume  = {150},
		number  = {3},
		pages   = {593--630},
		doi     = {10.1007/BF02096964},
	}

@article{JeffreyWeitsman1994,
		author  = {Jeffrey, Lisa C. and Weitsman, Jonathan},
		title   = {Toric structures on the moduli space of flat connections on a {R}iemann surface: volumes and the moment map},
		journal = {Adv. Math.},
		year    = {1994},
		volume  = {106},
		number  = {2},
		pages   = {151--168},
		doi     = {10.1006/aima.1994.1054},
	}

@article{Kawasaki1979,
		author  = {Kawasaki, Tetsuro},
		title   = {The {R}iemann--{R}och theorem for complex {$V$}-manifolds},
		journal = {Osaka Math. J.},
		year    = {1979},
		volume  = {16},
		number  = {1},
		pages   = {151--159},
		url     = {http://projecteuclid.org/euclid.ojm/1200771835},
	}

@book{Kempf1993,
		author    = {Kempf, George},
		title     = {Algebraic varieties},
		series    = {London Mathematical Society Lecture Note Series},
		volume    = {172},
		publisher = {Cambridge University Press},
		address   = {Cambridge},
		year      = {1993},
		isbn      = {0-521-42613-8},
		doi       = {10.1017/CBO9781107359956},
	}

@article{Konno2000,
		author  = {Konno, Hiroshi},
		title   = {Cohomology rings of toric hyperk{\"a}hler manifolds},
		journal = {Internat. J. Math.},
		year    = {2000},
		volume  = {11},
		number  = {8},
		pages   = {1001--1026},
		doi     = {10.1142/S0129167X00000490},
	}

@incollection{Konno2008,
		author    = {Konno, Hiroshi},
		title     = {The geometry of toric hyperk{\"a}hler varieties},
		booktitle = {Toric Topology},
        editor = {Harada, Megumi and Karshon, Yael and Masuda, Mikiya and Panov, Taras},
		series    = {Contemporary Mathematics},
		volume    = {460},
		pages     = {241--260},
		publisher = {American Mathematical Society},
		address   = {Providence, RI},
		year      = {2008},
		doi       = {10.1090/conm/460/09022},
	}

@article{Lerman1995,
		author  = {Lerman, Eugene},
		title   = {Symplectic cuts},
		journal = {Math. Res. Lett.},
		year    = {1995},
		volume  = {2},
		number  = {3},
		pages   = {247--258},
		doi     = {10.4310/MRL.1995.v2.n3.a2},
	}

@article{LermanTolman1997,
		author  = {Lerman, Eugene and Tolman, Susan},
		title   = {{H}amiltonian torus actions on symplectic orbifolds and toric varieties},
		journal = {Trans. Amer. Math. Soc.},
		year    = {1997},
		volume  = {349},
		number  = {10},
		pages   = {4201--4230},
		doi     = {10.1090/S0002-9947-97-01821-7},
	}

@article{Martens2008,
		author  = {Martens, Johan},
		title   = {Equivariant volumes of non-compact quotients and instanton counting},
		journal = {Comm. Math. Phys.},
		year    = {2008},
		volume  = {281},
		number  = {3},
		pages   = {827--857},
		doi     = {10.1007/s00220-008-0501-x},
	}

@article{MaZhang2014,
	author  = {Ma, Xiaonan and Zhang, Weiping},
	title   = {Geometric Quantization for Proper Moment Maps: the {V}ergne Conjecture},
	journal = {Acta Math.},
	volume  = {212},
	pages   = {11--57},
	year    = {2014},
	doi     = {10.1007/s11511-014-0108-3}
}

@article{Meinrenken1998,
		author  = {Meinrenken, Eckhard},
        title = {Symplectic surgery and the {${\mathrm Spin}^{c}$}-{D}irac operator},
		journal = {Adv. Math.},
		year    = {1998},
		volume  = {134},
		number  = {2},
		pages   = {240--277},
		doi     = {10.1006/aima.1997.1701},
	}

@book{MFK1994,
		author    = {Mumford, David and Fogarty, John and Kirwan, Frances C.},
		title     = {Geometric invariant theory},
		edition   = {Third},
		series    = {Ergebnisse der Mathematik und ihrer Grenzgebiete. 2. Folge},
		volume    = {34},
		publisher = {Springer-Verlag},
		address   = {Berlin},
		year      = {1994},
		doi = {10.1007/978-3-642-57916-5},
		isbn      = {978-3-540-56963-3},
	}

@book{Nakajima1999,
    author    = {Nakajima, Hiraku},
    title     = {Lectures on {H}ilbert schemes of points on surfaces},
    series    = {University Lecture Series},
    volume    = {18},
    publisher = {American Mathematical Society},
    address   = {Providence, RI},
    year      = {1999},
    isbn      = {978-0-8218-1956-2},
    doi       = {10.1090/ulect/018},
}

@article{Nakajima2016,
		author  = {Nakajima, Hiraku},
		title   = {Towards a mathematical definition of {C}oulomb branches of $3$-dimensional $\mathcal{N}=4$ gauge theories, {I}},
		journal = {Adv. Theor. Math. Phys.},
		year    = {2016},
		volume  = {20},
		number  = {3},
		pages   = {595--669},
		doi     = {10.4310/ATMP.2016.v20.n3.a4},
		eprint  = {1503.03676},
	}

@article{Oprea2006,
		author  = {Oprea, Dragos},
		title   = {Tautological classes on the moduli spaces of stable maps to {$\mathbb{P}^r$} via torus actions},
		journal = {Adv. Math.},
		year    = {2006},
		volume  = {207},
		number  = {2},
		pages   = {661--690},
		doi     = {10.1016/j.aim.2006.01.002},
	}

@article{Paradan2009,
	author  = {Paradan, Paul-{\'E}mile},
	title   = {Formal Geometric Quantization},
	journal = {Ann. Inst. Fourier (Grenoble)},
	volume  = {59},
	number  = {1},
	pages   = {199--238},
	year    = {2009}
}

@phdthesis{Proudfoot2004,
		author = {Proudfoot, Nicholas J.},
		title  = {Hyperk{\"a}hler analogues of K{\"a}hler quotients},
		school = {University of California, Berkeley},
		year   = {2004},
	}

@book{Schrijver86,
	author    = {Schrijver, Alexander},
	title     = {Theory of linear and integer programming},
	series    = {Wiley-Interscience Series in Discrete Mathematics},
	address  = {Chichester},
	publisher = {John Wiley \& Sons Ltd.},
	year      = {1986},
	isbn      = {0-471-90854-1},
}

@article{Skowera2013,
		author  = {Skowera, Jonathan},
		title   = {Bia{\l}ynicki-{B}irula decomposition of {D}eligne--{M}umford stacks},
		journal = {Proc. Amer. Math. Soc.},
		year    = {2013},
		volume  = {141},
		number  = {6},
		pages   = {1933--1937},
		doi     = {10.1090/S0002-9939-2012-11464-9},
	}

@article{Teleman2000,
		author  = {Teleman, Constantin},
		title   = {The quantization conjecture revisited},
		journal = {Ann. of Math. (2)},
		year    = {2000},
		volume  = {152},
		number  = {1},
		pages   = {1--43},
		doi     = {10.2307/2661378},
	}

@article{Thaddeus1996,
	author  = {Thaddeus, Michael},
	title   = {Geometric invariant theory and flips},
	journal = {J. Amer. Math. Soc.},
	year    = {1996},
	volume  = {9},
	number  = {3},
	pages   = {691--723},
	doi     = {10.1090/S0894-0347-96-00204-4},
}

@article{Vergne1996,
		author  = {Vergne, Mich\`ele},
		title   = {Equivariant index formulas for orbifolds},
		journal = {Duke Math. J.},
		year    = {1996},
		volume  = {82},
		number  = {3},
		pages   = {637--652},
		doi     = {10.1215/S0012-7094-96-08226-5},
	}

@article{Verlinde1988,
		author  = {Verlinde, Erik},
		title   = {Fusion rules and modular transformations in {$2$D} conformal field theory},
		journal = {Nuclear Phys. B},
		year    = {1988},
		volume  = {300},
		number  = {3},
		pages   = {360--376},
		doi     = {10.1016/0550-3213(88)90603-7},
	}

@article{Weitsman1992,
		author  = {Weitsman, Jonathan},
		title   = {Real polarization of the moduli space of flat connections on a {R}iemann surface},
		journal = {Comm. Math. Phys.},
		year    = {1992},
		volume  = {145},
		number  = {3},
		pages   = {425--433},
		doi     = {10.1007/BF02099391},
	}

@article{Weitsman2001,
	author  = {Weitsman, Jonathan},
	title   = {Non-Abelian Symplectic Cuts and the Geometric Quantization of Noncompact Manifolds},
	journal = {Lett. Math. Phys.},
	volume  = {56},
	pages   = {31--40},
	year    = {2001}
}

@article{Wu2000,
		author  = {Wu, Siye},
		title   = {A note on higher cohomology groups of {K}{\"a}hler quotients},
		journal = {Ann. Global Anal. Geom.},
		year    = {2000},
		volume  = {18},
		number  = {6},
		pages   = {569--576},
		doi     = {10.1023/A:1006712005064},
	}
